\documentclass[11pt]{amsart}
\usepackage[english]{babel}
\usepackage[margin=1.4in]{geometry}
\usepackage{amsmath}
\usepackage{amssymb}
\usepackage{amsthm}

\usepackage{palatino}

\usepackage{newpxmath}

\usepackage[backref=page]{hyperref}

\usepackage{yfonts}
\usepackage[T1]{fontenc}
\usepackage[utf8x]{inputenc}
\usepackage{enumerate}
\usepackage{enumitem}
\usepackage{verbatim}
\usepackage{graphicx}
\usepackage{faktor}
\usepackage{xcolor}
\usepackage{xfrac}
\usepackage{tikz,tikz-cd,tikz-3dplot}
\usepackage{pgfplots}
\usetikzlibrary{arrows,shadows,positioning, calc, decorations.markings, 
hobby,quotes,angles,decorations.pathreplacing,intersections,shapes}
\usepgflibrary{shapes.geometric}
\usetikzlibrary{fillbetween,backgrounds}

\tikzset{
  arrow/.pic={\path[tips,every arrow/.try,->,>=#1] (0,0) -- +(0,4pt);},
  pics/arrow/.default={triangle 90}
}

\tikzset{->-/.style={decoration={
  markings,
  mark=at position .6 with {\arrow{latex}}},postaction={decorate}}
  }

\tikzset{
  c/.style={every coordinate/.try}
}

\usepackage[all]{xy}
\usepackage{hyperref}
\usepackage[normalem]{ulem}
\usepackage{setspace}

\usepackage{calrsfs}
\DeclareMathAlphabet{\pazocal}{OMS}{zplm}{m}{n}

\newcommand{\tropC}{\scalebox{0.8}[1.3]{$\sqsubset$}}

\newcommand{\tropT}{\scalebox{0.9}[1.2]{$\top$}}
\newcommand{\tropD}{\scalebox{0.9}[1.3]{$\triangle$}}

\newcommand{\trop}{\operatorname{trop}}

\newcommand{\CO}{\mathcal{O}}

\newcommand{\CL}{\mathcal{L}}
\newcommand{\tC}{\widetilde{C}}
\newcommand{\Do}{\pazocal D_1}
\newcommand{\Dt}{\pazocal D_2}
\newcommand{\VZ}{\pazocal{V\!Z}_{2,n}}

\newcommand{\M}[4]{\overline{\pazocal M}_{#1,#2}(#3,#4)}
\newcommand{\PP}{\mathbb P}

\newcommand{\NN}{\mathbb N}
\newcommand{\OO}{\mathcal O}
\newcommand{\ZZ}{\mathbb Z}
\renewcommand{\to}{\rightarrow}

\newcommand{\Aaff}{\mathbb A}

\newcommand{\oM}{\overline{\pazocal M}}

\newcommand{\R}{\operatorname{R}}
\newcommand{\val}{\operatorname{val}}

\newcommand{\dvr}{\varUpsilon}
\newcommand{\cC}{\mathcal C}
\newcommand{\TTC}{\doublewidetilde{\mathcal{C}_{}}}
\newcommand{\TC}{\widetilde{\mathcal C}}
\newcommand{\TTA}{\doublewidetilde{\mathcal A_2}}
\newcommand{\TAt}{\widetilde{\mathcal A_2}}
\newcommand{\TAtn}{\widetilde{\mathcal A_{2,n}}}
\newcommand{\Cp}{\mathcal{C}^\prime}
\newcommand{\Cpp}{\mathcal{C}^{\prime\prime}}
\newcommand{\Tp}{\mathcal{T}^\prime}
\newcommand{\Tpp}{\mathcal{T}^{\prime\prime}}
\newcommand{\Z}{\pazocal Z}
\newcommand{\tZ}{\widetilde{\pazocal Z}}
\newcommand{\Zp}{\pazocal{Z}^\prime}
\newcommand{\Zpp}{\pazocal{Z}^{\prime\prime}}
\newcommand{\Li}{\mathcal L}
\newcommand{\LLi}{\widetilde{\mathcal L}}
\newcommand{\OC}{\overline{\mathcal C}}
\newcommand{\OOC}{\overline{\overline{\mathcal C}}}

\newcommand{\bq}{\begin{equation}}
\newcommand{\eq}{\end{equation}}
\newcommand{\ba}{\begin{aligned}}
\newcommand{\ea}{\end{aligned}}
\newcommand{\be}{\begin{enumerate}}
\newcommand{\ee}{\end{enumerate}}
\newcommand{\on}{\operatorname}
\newcommand{\bsm}{\left(\begin{smallmatrix}}
\newcommand{\esm}{\end{smallmatrix}\right)}                   
\newcommand{\bpm}{\begin{pmatrix}}
\newcommand{\epm}{\end{pmatrix}}
\newcommand{\overliner}{\begin{displaymath}\begin{array}{cccc}}
\newcommand{\earr}{\end{array}\end{displaymath}}
\newcommand{\overlinerl}{\begin{displaymath}\begin{array}{lcl}}
\newcommand{\earrl}{\end{array}\end{displaymath}}
\newcommand{\overlinel}{\begin{displaymath}\begin{array}{l}}
\newcommand{\earl}{\end{array}\end{displaymath}}
\newcommand{\bxym}{ \begin{displaymath}\xymatrix }
\newcommand{\exym}{\end{displaymath}}
\newcommand{\bcd}{\begin{center}\begin{tikzcd}}
\newcommand{\ecd}{\end{tikzcd}\end{center}}

\newcommand{\Hom}{\operatorname{Hom}}

\newcommand{\ev}{\operatorname{ev}}

\theoremstyle{plain}
\newtheorem{thm}{Theorem}[section]

\newtheorem{lem}[thm]{Lemma}
\newtheorem{prop}[thm]{Proposition}

\newtheorem{cor}[thm]{Corollary}
\newtheorem*{thm*}{Theorem}
\newtheorem*{matteo*}{Main Theorem}

\theoremstyle{definition}
\newtheorem{example}[thm]{Example}
\newtheorem{exa}[thm]{Example}
\newtheorem{dfn}[thm]{Definition}

\newtheorem{rem}[thm]{Remark}
\newtheorem{Convention}[thm]{Convention}

\makeatletter
\newcommand{\doublewidetilde}[1]{{%
  \mathpalette\double@widetilde{#1}%
}}
\newcommand{\double@widetilde}[2]{%
  \sbox\z@{$\m@th#1\widetilde{#2}$}%
  \ht\z@=.9\ht\z@
  \widetilde{\box\z@}%
}
\makeatother

\newtheorem{innercustomgeneric}{\customgenericname}
\providecommand{\customgenericname}{}
\newcommand{\newcustomtheorem}[2]{%
  \newenvironment{#1}[1]
  {%
   \renewcommand\customgenericname{#2}%
   \renewcommand\theinnercustomgeneric{##1}%
   \innercustomgeneric
  }
  {\endinnercustomgeneric}
}

\newcustomtheorem{customthm}{Theorem}
\newcustomtheorem{customlemma}{Lemma}
\newcustomtheorem{customrmk}{Remark}

\def\Dk{\tikz[baseline=-3pt]{
\draw (0,0)node[left]{\scriptsize $g=2,d=0$} -- (1,.5) node[right]{\scriptsize $g=0,d_1$} (0,0)--(1,-.5)node[right]{\scriptsize $g=0,d_k$} (0,0) -- (1.2,.2)node[right]{\scriptsize $g=0,d_2$};
\draw (1.2,-.2) node[right]{\scriptsize $\ldots$};
\draw (0,0) circle(2pt)[fill=white];
\fill (1,.5) circle (2pt) (1,-.5) circle (2pt) (1.2,.2) circle (2pt)}
}

\def\Ek{\tikz[baseline=-3pt]{
\draw (-1.5,0) node[left]{\scriptsize $g=1,d_0$} -- (0,0);
\draw (0,0)node[above]{\scriptsize \begin{tabular}{c} $g=1,$ \\ $d=0$ \end{tabular}} -- (1,.5) node[right]{\scriptsize $g=0,d_1$} (0,0)--(1,-.5)node[right]{\scriptsize $g=0,d_k$} (0,0) -- (1.2,.2)node[right]{\scriptsize $g=0,d_2$};
\draw (1.2,-.2) node[right]{\scriptsize $\ldots$};
\draw (0,0) circle(2pt)[fill=white];
\fill (-1.5,0) circle (2pt) (1,.5) circle (2pt) (1,-.5) circle (2pt) (1.2,.2) circle (2pt)}
}

\def\Ekk{\tikz[baseline=-3pt]{
\draw (-3,0)node[above]{\scriptsize \begin{tabular}{c} $g=1,$ \\ $d=0$ \end{tabular}} -- (-4,.5) node[left]{\scriptsize $g=0,d_{1,1}$} (-3,0)--(-4,-.5)node[left]{\scriptsize $g=0,d_{1,k_1}$} (-3,0) -- (-4.2,.2)node[left]{\scriptsize $g=0,d_{1,2}$};
\draw (-3,0) -- (-1.5,0) node[below]{\scriptsize $g=0,d_0$} -- (0,0);
\draw (0,0)node[above]{\scriptsize \begin{tabular}{c} $g=1,$ \\ $d=0$ \end{tabular}} -- (1,.5) node[right]{\scriptsize $g=0,d_{2,1}$} (0,0)--(1,-.5)node[right]{\scriptsize $g=0,d_{2,k_2}$} (0,0) -- (1.2,.2)node[right]{\scriptsize $g=0,d_{2,2}$};
\draw (1.2,-.2) node[right]{\scriptsize $\ldots$} (-4.2,-.2) node[right]{\scriptsize $\ldots$};
\draw (0,0) circle(2pt)[fill=white] (-3,0) circle(2pt)[fill=white];
\fill (-1.5,0) circle (2pt) (1,.5) circle (2pt) (1,-.5) circle (2pt) (1.2,.2) circle (2pt) (-4,.5) circle (2pt) (-4,-.5) circle (2pt) (-4.2,.2) circle (2pt)}
}

\def\pacman{\tikz[baseline=-3pt]{
\draw (1,0) -- (2,.5) node[right]{\scriptsize $g=0, d_1$} (1,0)--(2,-.5)node[right]{\scriptsize $g=0, d_k$} (1,0) -- (2.2,.2)node[right]{\scriptsize $g=0, d_2$};
\draw (2.2,-.2) node[right]{\scriptsize $\ldots$};
\draw (0,0)node[left]{\scriptsize $g=0,d_0$} to[out=20,in=160] (1,0) node[above]{\scriptsize \begin{tabular}{c} $g=1,$ \\ $d=0$ \end{tabular}};
\draw (0,0) to[out=-20,in=-160] (1,0);
\fill (0,0) circle (2pt);
\draw[fill=white] (1,0) circle (2pt);
\fill  (2,.5) circle (2pt) (2,-.5) circle (2pt) (2.2,.2) circle (2pt)}
}

\def\hypell{\tikz[baseline=-3pt]{
\draw (0,0)node[left]{\scriptsize $g=2,d_0=2$} -- (1,.5) node[right]{\scriptsize $g=0, d_1$} (0,0)--(1,-.5)node[right]{\scriptsize $g=0, d_k$} (0,0) -- (1.2,.2)node[right]{\scriptsize $g=0, d_2$};
\draw (1.2,-.2) node[right]{\scriptsize $\ldots$};
\draw (0,0) circle(2pt)[fill=gray!50];
\fill (1,.5) circle (2pt) (1,-.5) circle (2pt) (1.2,.2) circle (2pt)}
}

\newcommand{\thismonth}{\ifcase\month 
  \or January\or February\or March\or April\or May\or June%
  \or July\or August\or September\or October\or November%
  \or December\fi}

\title[A modular desingularisation of $\overline{\pazocal{M}}_2(\PP^r,d)^\text{main}$: log geometry \& singularities]{A smooth compactification of the space \\ of genus two curves in projective space \\ \tiny{via logarithmic geometry and Gorenstein curves}}
\author{Luca Battistella and Francesca Carocci}

\begin{document}

{\setstretch{1.3}
\begin{abstract}
We construct a modular desingularisation of $\overline{\pazocal{M}}_{2,n}(\PP^r,d)^\text{main}$. The geometry of Gorenstein singularities of genus two leads us to consider maps from prestable admissible covers: with this enhanced logarithmic structure, it is possible to desingularise the main component by means of a logarithmic modification. Both isolated and non-reduced singularities appear naturally. Our construction gives rise to a notion of reduced Gromov-Witten invariants in genus two.
\end{abstract}}

\maketitle
\tableofcontents

\vspace{-1cm}

\setstretch{1.08}

\section{Introduction}
Modern enumerative geometry is based on a series of compactifications of the moduli space of smooth embedded curves of genus $g$ and curve class $\beta\in H^+_2(X,\ZZ)$ in a smooth projective variety $X$, and on their virtual intersection theory. The boundary of M. Kontsevich's space of stable maps represents maps from nodal curves, including multiple covers and contracted components. The genus zero theory of projective space provides a smooth compactification with normal crossing boundary. In higher genus, instead, contracted subcurves and finite covers may give rise to boundary components of excess dimension. Even though this moduli space satisfies R. Vakil's Murphy's law, a desingularisation of the main component certainly exists after H. Hironaka's work on resolution of singularities. The main application of the methods developed in this paper is an \emph{explicit modular desingularisation} of the main component of the moduli space of stable maps to projective space in genus two and degree $d\geq3$.

 \begin{thm*}
 There exists a logarithmically smooth and proper DM stack $\VZ(\PP^r,d)$ over $\mathbb C$, with locally free logarithmic structure (and therefore smooth) and a birational forgetful morphism to $\M{2}{n}{\PP^r}{d}^{\text{main}}$, parametrising the following data:
 \begin{itemize}[leftmargin=.5cm]
  \item a pointed admissible hyperelliptic cover \[\psi\colon (C,D_R,x_1,\ldots,x_n)\to(T,D_B,\psi(x_1),\ldots,\psi(x_n))\]
  where $C$ is a prestable curve of arithmetic genus two, $D_R$ and $D_B$ are length six (ramification and branch) divisors, and $(T,D_B,\psi(\mathbf{x}))$ is a stable rational tree;
  \item a map $f\colon C\to \PP^r$;
  \item a bubbling destabilisation $C\leftarrow\widetilde{C}$ and a contraction $\widetilde{C}\to \overline{C}$ to a curve with Gorenstein singularities, such that $f$ factors through $\bar{f}\colon\overline{C}\to \PP^r$, and the latter is not special on any subcurve of $\overline{C}$.
 \end{itemize}
\end{thm*}

More generally, for any smooth projective variety $X$, we construct a proper moduli space $\VZ(X,\beta)$ admitting a perfect obstruction theory and defining \emph{reduced Gromov-Witten invariants}. When $X$ is a projective complete intersection, they satisfy the \emph{quantum Lefschetz} hyperplane principle. In general, we expect them to have a simpler enumerative content compared to standard Gromov-Witten invariants.

\subsection{Main and boundary components}\label{sec:mainandboundary}
The moduli space of stable maps to $X$ (see \cite{Kontsevich}) represents $f\colon C\to X$, where $C$ is a \emph{prestable} (nodal and reduced) curve such that every rational component of $C$ contracted by $f$ contains at least three \emph{special points} (markings and nodes).

When $X=\PP^r$, $f$ is equivalently determined by the data of a line bundle $L$ on $C$, and $r+1$ sections of $L$ that do not vanish simultaneously on $C$. Forgetting the sections we obtain a morphism to the universal Picard stack, parametrising pairs $(C,L)$ of a prestable curve and a line bundle on it:
\[\M{g}{n}{\PP^r}{d}\to\mathfrak{P}ic_{\mathcal C_{g,n}/\mathfrak M_{g,n}}.\]
Obstructions lie in $H^1(C,L)$ \cite{WangJie}. This implies that $\M{0}{n}{\PP^r}{d}$ is smooth. On the other hand, for $g\geq 1$, boundary components may arise where $L$ restricts to a special line bundle on a subcurve of $C$. When $d>2g-2$ we can identify a \emph{main component}: the closure of the locus of maps from a smooth source curve.

In genus one, the generic point of a boundary component has a contracted subcurve of genus one, together with $k\leq r$ tails of positive degree \cite{Vre}.

In genus two, two types of boundary phenomena occur:
\begin{itemize}[leftmargin=0.5cm]
 \item a subcurve of positive genus is contracted, or
 \item the restriction of $f$ to the minimal subcurve of genus two (\emph{core}) is the hyperelliptic cover of a line in $\PP^r$.
\end{itemize}
As an example of the second phenomenon, the main component of $\overline{\pazocal M}_2(\PP^3,5)$ has dimension $20$; but the locus of maps from a reducible curve $C=Z\bigcup_q R$, with $Z$ of genus two covering a line two-to-one, and $R\simeq\PP^1$ parametrising a twisted cubic, has dimension $21$.

We stress that $\overline{\pazocal M}_{2,n}(\PP^r,d)^{\text{main}}$ has no natural modular interpretation, since it is obtained by taking the closure of the nice locus. Moreover, its singularities along the boundary can be nasty \cite{VakilMurphy}. In \cite{BCquartics} we draw the consequences of the present construction, analysing the locus of smoothable maps in $\overline{\pazocal M}_{2}(\PP^2,4)$, a moduli space with more than twenty irreducible components.

\subsection{Good factorisation through a Gorenstein curve}
Our approach to the desingularisation of the main component takes off from a simple observation: a line bundle of degree at least $2g-1$ on a \emph{minimal Gorenstein} curve of genus $g$ has vanishing $h^1$; we have already used this implicitly for nodal curves in the discussion of the irreducible components of $\M{g}{n}{\PP^r}{d}$. The Gorenstein assumption makes this fact a straightforward consequence of Serre duality. Minimality is a weaker notion than irreducibility, but it is needed to ensure that the line bundle has sufficiently positive degree on every subcurve.

This observation raises a natural question: is it possible to replace every $f\colon C\to\PP^r$ with a ``more positive''/``less obstructed'' $\overline f\colon\overline C\to \PP^r$ by ``contracting'' the higher genus subcurves on which $f^*\OO_{\PP^r}(1)$ is special? The answer is ``no'', since every such map is \emph{smoothable}, i.e. it can be deformed into a map from a smooth curve.\footnote{The answer could also be ``yes'': it has been shown by M. Viscardi in genus one that maps from Gorenstein curves satisfying certain stability conditions give rise to irreducible compactifications of $\pazocal{M}_{1,n}(\PP^r,d)$ \cite{VISC}; yet, their deformation theory is hard to grasp, because such is the deformation theory of the singularities that are involved \cite{SMY2}.}

This, on the other hand, gives us a strategy to study the main component:
\begin{enumerate}[leftmargin=.5cm]
 \item\label{item:contraction} Construct a universal contraction $\mathcal C\to\OC$ to a Gorenstein curve, by collapsing the subcurves of $C$ on which $f$ has low degree. This Gorenstein curve will not be ready available on $\M{g}{n}{\PP^r}{d}$, because the contraction map $\mathcal C\to\OC$ has moduli (called \emph{crimping spaces} in \cite{vdW}, and \emph{moduli of attaching data} in \cite{SMY2}); we first need to introduce a modification of $\M{g}{n}{\PP^r}{d}$ along the boundary.
 \item\label{item:factorisation} Take only those maps that admit a factorisation through $\OC$: as we have discussed, these are all smoothable, so the moduli space is at the very least irreducible. It provides a birational model of $\M{g}{n}{\PP^r}{d}^{\text{main}}$, with the advantage of admitting a modular interpretation.
\end{enumerate}
Moreover, this space is unobstructed over a base that can be assumed to be smooth in the low-genus examples that we have at our disposal: \cite{RSPW1}, where the base is a logarithmic modification of the moduli space of prestable curves; and this paper, where the base is a logarithmic modification of the moduli space of hyperelliptic admissible covers. Once such a moduli space is constructed, the proof of smoothness for target $\PP^r$ is entirely conceptual; furthermore, the same methods can be employed to approach the study of different targets (e.g. products of projective spaces, toric varieties, flag varieties) and stability conditions (e.g. quasimaps).

While point \eqref{item:factorisation} was essentially established for us by D. Ranganathan, K. Santos-Parker, and J. Wise, point \eqref{item:contraction} is at all open for $g\geq 2$. In the present work, we make a hopefully meaningful step in this direction.

\subsection{Logarithmic geometry \& singularities}
The moduli space of prestable curves has a natural logarithmic structure induced by its boundary divisor, thus keeping track of the nodes and their smoothing parameters \cite{KatoF}. This induces a logarithmic structure on the moduli space of maps as well. It is a natural question whether the desingularisation of the main component can be achieved by means of a logarithmic modification; it is indeed the case in genus one \cite{RSPW1}, but not in our construction.  We think that the following are truly high-genus phenomena. 

\subsubsection{Augmenting the logarithmic structure: admissible covers} Our first most relevant finding is that, instead, in genus two it is necessary to enrich the logarithmic structure of the base by passing to a moduli space of admissible covers \cite{HarrisMumford}.

Every smooth curve of genus two is hyperelliptic, i.e. the canonical class induces a two-to-one cover of a line branched along six points. The hyperelliptic cover is essentially unique, but, when the curve becomes nodal and in the presence of markings, the uniqueness is lost. It can be restored by making it part of the moduli functor. The resulting space of admissible covers is as nice as the moduli space of curves, but it has the advantage of encoding the Brill-Noether theory of the curve in the logarithmic structure \cite{Mochizuki}. The necessity of such enrichment can be understood already by looking at the isolated Gorenstein singularities of genus two \cite{B}: there are two families of these, basically corresponding to the choice of either a Weierstrass point or a conjugate pair in the semistable model. In order to tell these two cases apart logarithmically, in the construction of $\mathcal C\to\OC$, we start by considering stable maps from the source of an admissible cover.

\subsubsection{Realisable tropical canonical divisors} This points to our second finding. The line bundle giving the contraction $\mathcal C\to\Cp$ (the first step towards $\OC$) is a vertical twist of the relative dualising line bundle. The twist is indeed logarithmic, and the piecewise-linear function on $\operatorname{trop}(\mathcal C)$ determining it is nothing but a \emph{tropical canonical divisor} satisfying certain requirements of compatibility with the admissible cover, and in particular realisable \cite{MUW}. Pulling back from the target of the admissible cover, which is a rational curve, allows us to simplify several computations. The subdivision of the tropical moduli space according to the domain of linearity of such function should therefore be thought of as parametrising all possible Gorenstein contractions of $\mathcal C$ compatible with the degree of $f$.


\subsubsection{Non-isolated singularities} Finally, we underline that non-reduced curves appear as fibres of $\OC$, which requires a careful study of what we call \emph{tailed ribbons}, and further distinguishes our work from its genus one ancestor. This should not come as a surprise: first, non-reduced curves can still be Gorenstein if the nilpotent structure is supported along one-dimensional components, rather than isolated points; second, \emph{ribbons} were introduced in the '90s as limits of smooth canonical curves when the latter tend to the hyperelliptic locus in moduli \cite{Fong}. They appear naturally in our work at the intersection of the main component with the hyperelliptic ones. It is possibly less expected that they show up as well when the core is contracted, as a way of interpolating between isolated singularities whose special branches differ. The construction of $\OC$ is concluded by gluing in $\Cp$ the hyperelliptic cover of a genus two subcurve supported along a boundary divisor.

As a note for future investigation, we remark that we perform both the contraction $\mathcal C\to\mathcal C^{\prime}$ and the pushout $\mathcal C^{\prime}\to\OC$ in sufficiently general families. We wonder whether a pointwise construction might be possible, realising both steps as special instances of a more general pushout of logarithmic subcurves, in the spirit of S. Bozlee's \cite{Bozlee}.

\subsection{Relation to other work and further directions}

\subsubsection{Local equations and resolution by blowing up}\label{sec:desingbyblowingup}
It is always possible to find a local embedding of $\M{g}{n}{\PP^r}{d}$ in a smooth ambient space by looking at:
\[\operatorname{Tot}(\pi_*L)^{\oplus r+1}\subseteq\operatorname{Tot}(\pi_*L\otimes\OO_C(n))^{\oplus r+1}\]
over the Picard stack, where $\OO_C(1)$ is a relative polarisation for the universal curve $\pi\colon\mathcal C\to\mathfrak{P}ic_{g,n}$ \cite[\S 3.2]{ToricQuasimaps}. If the polarisation is chosen carefully, the embedding is given by $\approx\! g$ local equations, repeated $r+1$ times. An approach to the desingularisation of the main component is by blowing up according to these local equations \cite{HL,HLderivedresolution,HLN}. Possibly, the hardest task is to provide a compatible global description of the blow-up centres and their sequence. This method is close in spirit to the original construction of R. Vakil and A. Zinger \cite{ZingerSharp,VZ} (in particular, it involves an iterated blow-up procedure and a good deal of book-keeping), and it has the advantage of simultaneously ``desingularising'' the sheaves $\pi_*\mathcal L^{\otimes k}, k\geq 1,$ on the main component, making them into vector bundles. This in turn makes the theory of projective complete intersections accessible via torus localisation \cite{ZingerredCYhyp,PoparedCYci}.
 
Y.Hu, J.Li, and J.Niu \cite{HLN} carry out this strategy in genus two. For the reader's benefit, we provide a brief conversion chart: in Round $3$, Phase $3$ of their routine, they blow up loci in the Picard stack, where the line bundle is required to restrict to the canonical bundle on the core; this part is subsumed in our work by starting from a birational model of the moduli space of weighted curves, namely the space of admissible covers. In Round $1$ they blow up loci in $\mathfrak{M}_2^\text{wt}$ that correspond explicitly to the boundary phenomena that we have summarised in \S \ref{sec:mainandboundary}. In Round $2$ and $3$ they further blow up loci contained in the exceptional divisor of the previous rounds; these loci often depend on the attaching points of some rational tails being Weierstrass, or conjugate. In the language we have adopted, all of the blow-up centres are encoded in \emph{the domain of linearity of a certain piecewise linear function} - a tropical canonical divisor -, while the Brill-Noether theory of the curve is built in the moduli space at the level of admissible covers already. Our construction is thus more intrinsic. We remark, for example, that the genus one phenomena (corresponding to Round $1$, Phase $2-4$ in \cite{HLN}) are addressed uniformly within this language, with the tropical canonical divisor resembling the function ``distance from the core'' that is at the base of \cite{RSPW1}. The desingularisation of Hu, Li, and Niu seems more efficient than ours, in the (vague) sense that the information encoded by the admissible cover \emph{far from the special branch} is useless for the rest of our construction. Their procedure appears to select a ``special branch'' to start with in Round $1$, and then continue blowing up accordingly only ``around that branch''; the price to be paid seems to be a certain non-canonicity of the blow-up procedure, with subdivisions appearing that are not quite natural from our perspective (when the core is reducible, see Example \ref{exa:reduciblecore}).

In \cite{HuNiuTwisted}, Hu and Niu reconstruct the blow-up by gluing projective bundle strata indexed by \emph{treelike structures} and \emph{level trees}; these data are reminiscent of piecewise-linear functions on tropical curves, though missing both the slope and the metric data. These authors have already pointed out the similarities between their indexing set and the combinatorial data appearing in the moduli space of multi-scale differentials \cite{BCGGM}. This relation certainly deserves further attention: we think that canonical divisors on tropical curves could provide the right language to talk about it, and the geometry of Gorenstein curves could be the informing principle of further investigations. We expect a crucial ingredient for an all-genus desingularisation will be a good compactification of canonical curves.

\subsubsection{Computations in Gromov-Witten theory} Naive computations of what we may now expect to coincide with our reduced invariants have made seldom appearances in the literature: for example, with Zinger's enumeration of genus two curves with a fixed complex structure in $\PP^2$ and $\PP^3$ \cite{Zinger2}, and the computation of characteristic numbers of plane curves due to T. Graber, J. Kock, and R. Pandharipande \cite{KockGraberPandha}. To make the relation with the latter work precise, we should first extend our methods to the analysis of relative and logarithmic stable maps (compare with \cite{BNR} and \cite{RSPW2} in genus one).

$\VZ(X,\beta)$ is only the main component of a moduli space of aligned admissible maps $\mathcal{A}_{2,n}(X,\beta)$, which dominates $\overline{\pazocal{M}}_{2,n}(X,\beta)$ and is virtually birational to it. The virtual class of $\mathcal{A}_{2,n}(X,\beta)$ is expected to split as the sum of its main and boundary components; the contribution of the latter should be expressible in terms of genus zero and reduced genus one invariants via virtual pushforward \cite{ManolachePush}. This would deliver an extension of the Li-Zinger's formula to genus two (see \cite{LZ,ZingerStdVSRed,CoatesManolache}). Together with a localisation computation of reduced invariants, as mentioned in the previous section, this would provide an alternative proof of the genus two mirror theorem for the quintic threefold \cite{GJR,CGL}. Along the same lines, it would ease the computation of genus two Gromov-Witten invariants of Fano and Calabi-Yau complete intersections in projective space.

\subsubsection{Logarithmic maps to toric varieties and tropical realisability} In \cite{RSPW2}, Ranganathan, Santos-Parker, and Wise apply their techniques to desingularising the space of genus one logarithmic maps to a toric variety with respect to its boundary. As an application, they are able to solve the \emph{realisability problem} for tropical maps of genus one \cite{Speyer}: which of them arise as the tropicalization of a map from a smooth genus one curve to an algebraic torus? With Ranganathan, we are working towards a similar result in genus two - where, as far as we know, there is at the moment no reasonable guess as to what the full answer should be, although a clear understanding of the moduli space of tropical curves has been obtained \cite{CuetoMarkwig}.

\subsubsection{Birational geometry of moduli spaces of curves} In \cite{B}, the first author produces a sequence of alternative compactifications of $\pazocal M_{2,n}$ based on replacing genus one and two subcurves with few special points by isolated Gorenstein singularities. Although we do not discuss it here, the techniques developed in this paper also provide a resolution of the rational maps among these spaces. Moreover, the universal Gorenstein curve constructed in this paper unveils the possibility of defining new birational models of $\overline{\pazocal M}_{2,n}$ by including non-reduced curves as well - contrary to $\overline{\pazocal M}_{2,n}(m)$, these models could respect the $\mathfrak S_n$-symmetry in the markings. It would be interesting to compare them with \cite{PolishchukJohnson}, and to establish their position in the Hassett-Keel program \cite{SMY2,AFS1}. Recently, S. Bozlee, B. Kuo, and A. Neff have classified all the compactifications of $\pazocal{M}_{1,n}$ in the stack of Gorenstein curves with distinct markings \cite{BKN} - it turns out that there are many more than envisioned by Smyth, although the numerosity arises more from combinatorial than geometric complications. It would be interesting if all compactifications of $\pazocal{M}_{2,n}$ could be classified by a mixture of our techniques.

\subsection{Plan of the paper} In \S \ref{sec:preliminaries} we establish some language and background material concerning \emph{logarithmic curves}, their tropicalization, and the use of piecewise linear functions; \emph{admissible covers}, and their logarithmic structure; and \emph{Gorenstein curves}, with both isolated singularities and non-reduced structures, including a number of useful properties and formulae, 
the classification of isolated singularities of genus one (due to Smyth) and two (due to the first author), and a novel study of the non-isolated singularities of genus two with non-negative canonical class.

In \S \ref{sec:alignedadmissiblecovers} we introduce the key player: a subdivision of the tropical moduli space of weighted admissible covers based on aligning (ordering) the vertices of the tropical curve with respect to a piecewise linear function constructed from tropical canonical divisors. This subdivision induces a logarithmically \'etale model of the moduli space of weighted admissible covers. It is on this model that we are able to construct a universal family of Gorenstein curves. This process consists of two steps: first, a birational contraction $\mathcal C\to \Cp$; then, a pushout/normalisation $\Cp\to\OC$.

In \S \ref{sec:maps} we apply these methods to desingularise $\overline{\pazocal{M}}_{2,n}(\PP^r,d)^\text{main}$. To be more precise, we repeat the contraction/pushout process twice: first, to ensure that any curve of positive genus has weight at least one; then, to exclude the possibility that the map is hyperelliptic on a genus two subcurve. The intermediate space, although not being smooth, is already interesting in that its invariants satisfy quantum Lefschetz \cite{YPLee}.

\subsection{Notations and conventions} We work throughout over $\mathbb C$.

By \emph{trait} we mean the spectrum $\dvr$ of a discrete valuation ring. It only has two points: the closed one is often denoted by $s$ (or $0$), and the generic one by $\eta$.

Curves will always be projective and $S1$, i.e. without embedded points, but they may be non-reduced. Subcurves are not supposed to be irreducible, but they are usually connected. We call \emph{core} the minimal subcurve of genus two. We may refer to subcurves $C^\prime\subseteq C$ that intersect the rest of $C$ and the markings in only one (resp. two) at worst nodal point(s) - or markings - as \emph{tails} (resp. \emph{bridges}).

The \emph{dual graph} $\Gamma$ of a nodal curve $C$ has a vertex for every irreducible component, an edge for every node, and a leg for every marking (labelled or not); it is endowed with a genus function $g\colon V(\Gamma)\to\ZZ$, and the genus of the graph is given by:
\[g(\Gamma)=h^1(\Gamma)+\sum_{v\in V(\Gamma)}g(v).\]
The graph $\Gamma$ may be further weighted with a function $w\colon V(\Gamma)\to\ZZ$, that should be thought of as recording the degree of a line bundle on $C$.

A graph with no loops is called a \emph{tree}; its valence one vertices are called \emph{leaves} (sometimes, one of them plays a different role and it is therefore named the \emph{root}). A tree with only two leaves is a \emph{chain}.

A \emph{tropical curve} is a graph metrised in a monoid $\overline M$ (most basically, $\mathbb R_{\geq 0}$), i.e. a graph $\Gamma$ as above together with a length function $\ell\colon E(\Gamma)\to\overline M$ on the set of edges (legs are considered to be infinite, instead).

A \emph{logarithmic space} is denoted by $X=(\underline{X},M_X)$ with logarithmic structure $\alpha\colon M_X\to(\OO_X,\cdot)$ (also simply indicated by exponential/logarithmic notation). $M_X$ is a sheaf in the \'{e}tale topology of $\underline{X}$. We assume that logarithmic structures are \emph{fs}, i.e. fine (they admit an atlas of finitely generated and integral charts) and saturated. Let $\overline M_X$ denote the characteristic sheaf $M_X/\alpha^{-1}(\OO^*_X)$. This is a constructible sheaf of abelian groups.

A logarithmic space $X$ gives rise to a cone stack $\operatorname{trop}(X)$ via \emph{tropicalization}.
In the case of logarithmically smooth curves, the tropicalisation is a tropical curve in the above sense.
We abide to the rule that the tropicalization of $C$ should be denoted by the corresponding piecewise-linear character $\tropC$.

We have quite a few families of curves; usually we adopt the following notation:
\begin{itemize}[leftmargin=.5cm]
 \item $\pi\colon C\to S$ will denote a prestable curve, $\widetilde{C}$ (with $\tilde\pi$) a partial destabilisation of $C$, $C^\prime$ (with $\pi'$) a (not necessarily Gorenstein) contraction of $\widetilde{C}$, and $\overline{C}$ (with $\bar\pi$) a (not necessarily reduced) Gorenstein curve dominated by $C^\prime$;
 \item $f\colon C\to X$ will denote a (stable) map to (a smooth and projective) $X$;
 \item $\psi\colon C\to T$ will denote an admissible cover from a genus two to a rational curve, $\psi'\colon C^\prime \to T^\prime$ a double cover (where the curves are not necessarily prestable).
\end{itemize}

\subsection{Acknowledgements} We are grateful to Dhruv Ranganathan and Jonathan Wise for nurturing this project. We would also like to thank Daniele Agostini, Fabio Bernasconi, Sebastian Bozlee, David Holmes, Yi Hu, Navid Nabijou, and Jingchen Niu for  helpful conversations, and Zak Tur\v{c}inovi\'{c} for his support with the pictures. We thank Cristina Manolache and Richard Thomas for pushing us in this direction during our doctoral studies. We thank the anonymous referee for their quick and detailed response, that was of great help in improving the overall quality of the paper. L.B. is supported by the Deutsche Forschungsgemeinschaft (DFG, German Research Foundation) under Germany’s Excellence Strategy EXC-2181/1 - 390900948 (the Heidelberg STRUCTURES Cluster of Excellence). F.C. is supported by the starter grant “Categorified Donaldson–Thomas theory” No. 759967 of the European Research Council.
This work was initiated at the Max-Planck-Institut f\"ur Mathematik in Bonn; it was partly carried out during the workshop on ``Enumerative logarithmic geometry and mirror symmetry'' at the Mathematisches Forschungsinstitut Oberwolfach, and during a visit of the first author to the University of Edinburgh - we gratefully acknowledge all these institutions for their hospitality.

\section{Preliminaries and background material}\label{sec:preliminaries}
\subsection{Logarithmically smooth curves and their tropicalizations} Let $(S,M_S)$ be a logarithmic scheme. A family of \emph{logarithmically smooth curves} over $S$ is a proper and logarithmically smooth morphism $\pi\colon (C,M_C)\to (S,M_S)$ with connected one-dimensional (geometric) fibres, such that $\pi$ is integral and saturated. These hypotheses guarantee flatness and that the fibres are reduced. Logarithmically smooth curves naturally provide a compactification of the moduli space of smooth curves.

\subsubsection{Characterisation of logarithmically smooth curves}
F. Kato proved in \cite{KatoF} that logarithmically smooth curves have at worst nodal singularities; moreover, he provided the following local description: for every geometric point $p\in C$, there exists an \'etale local neighbourhood of $p$ in $C$ with a strict \'etale morphism to
\begin{description}
 \item[{smooth point}] $\Aaff^1_S$ with the log structure pulled back form the base;
 \item[{marking}] $\Aaff^1_S$ with the log structure generated by the zero section and $\pi^*M_S$;
 \item[{node}] $\OO_S[x,y]/(xy=t)$ for some $t\in\OO_S$, with semistable log structure induced by the multiplication map $\Aaff^2_S\to\Aaff^1_S$ and $t\colon S\to \Aaff^1$.
\end{description}

\noindent In the last case, the class of $\log(t)$ in $\overline{M}_S$ is called a \emph{smoothing parameter} for the node.

\noindent At times, we may have to consider more general logarithmic orbicurves \cite{OlssonLogCurves}.

\subsubsection{Minimal logarithmic structures}\label{sec:minimality}
For every prestable curve $\underline{\pi}\colon\underline{C}\to\underline{S}$, there is a \emph{minimal} logarithmic structure $M_S^\text{can}$ on $\underline{S}$, together with a logarithmically smooth enhancement $\pi\colon(\underline{C},M_C^\text{can})\to(\underline{S},M_S^\text{can})$, such that any other logarithmically smooth enhancement of $\underline{\pi}$ is pulled back from the minimal one. If $S=\{s\}$ is a geometric point, the characteristic sheaf of the minimal structure is freely generated by the smoothing parameters of the nodes of $C_s$: \[\overline{M}^{\text{can}}_S=\mathbb N^{\#E(\Gamma(C_s))}.\]

More generally, the concept of minimality - introduced by W.D. Gillam in \cite{Gillam} - serves the purpose of describing which stacks $\mathcal X$ over $(\text{LogSch})$ are induced by ordinary stacks $\pazocal X$ over $(\text{Sch})$ endowed with a logarithmic structure $\alpha\colon M_{\pazocal X}\to(\OO_X,\cdot)$. Gillam defines a \emph{minimal object} $x$ of $\mathcal X$ to be one such that for every solid diagram
\bcd
y^\prime\ar[rr,dashed] & & x\\
& y\ar[ul,"i"]\ar[ur,"j"below] &
\ecd
with $i$ and $j$ over the identity of underlying schemes, there exists a unique dashed arrow in $\mathcal X$ making the diagram commutative. He then shows that $\mathcal X$ is induced by an $(\pazocal X,M_{\pazocal X})$ if and only if:
\begin{enumerate}
 \item for every object $z$ of $\mathcal X$, there exist a minimal $x$, and $z\to x$ covering the identity of underlying schemes;
 \item for every $w\to x$ with $x$ minimal, the corresponding morphism of logarithmic schemes is strict if and only if $w$ is minimal as well.
\end{enumerate}

From this discussion, it follows that the stack of logarithmically smooth curves over $(\text{LogSch})$ is induced by a logarithmic stack over $(\text{Sch})$. It can be identified as the Artin stack of prestable (marked) curves endowed with the divisorial logarithmic structure corresponding to its normal crossings boundary.

\subsubsection{Tropicalization}
The combinatorial structure of an fs logarithmic space $X$ is encoded by its \emph{Artin fan} $\pazocal A_X$ \cite{ACMW,AWbirational}. An Artin fan is a logarithmic stack that  looks like the quotient of a toric variety by its dense torus $[V/T]$ locally in the strict \'etale topology. The morphism $\underline{X}\to\operatorname{Log}$ classifying the logarithmic structure on $X$ \cite{OlssonLogStacks} factors through $X\to\pazocal A_X$. The $2$-category of Artin fans is equivalent to the $2$-category of \emph{cone stacks}, i.e. stacks over the category of rational polyhedral cone complexes \cite[Theorem 6.11]{CCUW}. The latter should be thought of as a generalisation of the category of rational polyhedral cone complexes, where cones are allowed to be glued to themselves via automorphisms. In this way, we can associate to the logarithmic stack $X$ a cone stack $\operatorname{trop}(X)$ known as its \emph{tropicalization}, which is nothing but another incarnation of the Artin fan.
 When $X$ is a scheme with Zariski logarithmic structure (without monodromy), it is the generalised (standard) cone complex:
\[\operatorname{trop}(X)=\left(\coprod_{x\in X}\Hom(\overline{M}_{X,x},\mathbb R_{\geq 0})\right)\bigg/\sim\]
where $x$ runs through the schematic points of $X$, and gluing takes place along the face inclusions induced by the specialisation morphisms \cite[Appendix B]{GS}. The language developed in \cite{CCUW} allows us to talk about \emph{the tropicalization map} within  the category of logarithmic stacks. See also \cite{ACMUW,Ulirsch,Ulirsch2}.

For a logarithmically smooth curve $(C,M_C)$ over a geometric point $(S,M_S)$, the tropicalization $\tropC$ consists of its dual graph $\Gamma(C)$ metrised in $\overline{M}_S$; the length function $\ell\colon E(\Gamma(C))\to\overline{M}_S$ associates to each edge the smoothing parameter of the corresponding node. More precisely, vertices correspond to irreducible components weighted by their geometric genus, and there are (labelled) infinite legs corresponding to the markings. The tropicalization of $C$ is thus a family of classical (i.e. metrised in $\mathbb R$) tropical curves over the cone $\Hom(\overline{M}_S,\mathbb R_{\geq0})$. The construction can be generalised to more general base logarithmic schemes \cite[\S 7]{CCUW}. It induces a (strict, smooth, surjective) morphism:
\[\mathfrak M_{g,n}^\text{log}\to\widetilde{\mathfrak M}_{g,n}^\text{trop},\]
where the latter is the lift of the stack of tropical curves to the category of logarithmic schemes through the tropicalization map.

\subsubsection{Piecewise-linear functions and line bundles}\label{sec:PLfun}
The characteristic monoid at a node $q\in C$ is the amalgamated sum: \[\overline M_{C,q}=\overline M_{S,\pi(q)}\oplus_\NN \NN^{\oplus2},\]
where the map $\NN\to\overline M_{S,\pi(q)}$ is the smoothing parameter $\delta_q$, and $\NN\to\NN^{\oplus 2}$ is the diagonal. It has been noticed in \cite{GS} that:
\[\overline M_{C,q}^\text{gp}\simeq\{(\lambda_1,\lambda_2)\in\overline M_{S,\pi(q)}^{\text{gp},\oplus2}|\lambda_2-\lambda_1\in\ZZ\delta_q\}.\]

If $\tropC$ is a tropical curve metrised in $\overline M$, a \emph{piecewise-linear} (PL) \emph{function} on $\tropC$ with values in $\overline M^\text{gp}$ is a function $\lambda$ from the vertex set of $\tropC$ to $\overline M^\text{gp}$ such that, for any two adjacent vertices $v_1$ and $v_2$ connected by an edge $e_q$, we can write \[\lambda(v_2)-\lambda(v_1)=s(\lambda,e_q)\delta_q\] for some $s(\lambda,e_q)\in\mathbb Z$ (called the \emph{slope} of $\lambda$ along $e_q$, outgoing from $v_1$), where $\delta_q$ is the smoothing parameter of the node $q$ in $\overline M$.

Over a geometric point (or, more generally, if $\overline M_S$ and $\tropC$ are constant over $S$ - in particular, if no nodes are smoothed out) the previous observation, together with a similar analysis at the markings, shows (see \cite[Remark 7.3]{CCUW}) that
\[H^0(C,\overline M_C^\text{gp})=\{\text{PL functions on }\tropC\text{ with values in }\overline M_S^\text{gp}\}.\]

Similarly, every section of $\pi_*\overline M_C/\overline M_S$ is described by the collection of slopes of a rational function on $\tropC$.

The exact sequence:
\[0\to \OO_C^\times\to M_C^\text{gp}\to\overline M_C^\text{gp}\to 0,\]
shows that to every section $\lambda\in\Gamma(C,\overline M_C^\text{gp})$ there is an associated $\OO_C^\times$-torsor (or, equivalently, a line bundle) of lifts of $\lambda$ to $M_C^\text{gp}$. Under the assumption that $\overline M_S$ and $\tropC$ are constant over $S$, every vertex $v$ of $\tropC$ determines an irreducible component of $C$, and the restriction of $\OO_C(-\lambda)$ to $C_v$ is given explicitly by \cite[Proposition 2.4.1]{RSPW1}:
\[\OO_{C_v}(\lambda)\simeq \OO_{C_v}\left(\sum s(\lambda,e_q)q\right)\otimes\pi^*\OO_S(\lambda(v)),\]                                                                                                                                                                                                                                                                                                                                                                                                                                                                                                                
where $s(\lambda,e_q)$ denotes the outgoing slope of $\lambda$ along the edge corresponding to $q$ (either a marking or a node of $C$).   

If we started from $\lambda\in\Gamma(C,\overline M_C)$,
the associated line bundle $\OO_C(-\lambda)$ would come with a cosection $\OO_C(-\lambda)\to \OO_C$ (induced by the logarithmic structure). Such a cosection is not always injective, but when it is it defines an effective Cartier divisor on $C$; when it is not, it will vanish along components of $C$. Nonetheless, the association of this \emph{generalised effective Cartier divisor} to $\lambda$ behaves well under pullbacks, and in fact a useful point of view on logarithmic structures is to consider them as a functorial system of generalised effective Cartier divisors indexed by  $\overline M_C$ (see \cite{BV} for the details). See also \cite[p.9]{Bozlee} for a description in local coordinates.

Finally, we briefly recall the theory of divisors on tropical curves \cite{HaaseMusikerYu}. For this, let $\tropC$ be a tropical curve over $\mathbb R_{\geq 0}$, thought of as its metric realisation. A \emph{divisor} on $\tropC$ is a finite $\ZZ$-linear combination of points on $\tropC$; it is said to be \emph{effective} if all the coefficients are non-negative. The \emph{tropical canonical divisor} is:
\[K_{\tropC}=\sum_{v\in\tropC}(2g(v)-2+\on{val}(v))v.\]
It is effective unless the curve has rational tails. To a piecewise-linear function $\lambda$ on (a subdivision of) $\tropC$, we can associate a divisor supported on its the non-linearity locus:
\[\on{div}(\lambda)=\sum_{v\in\tropC}\left(\sum_{e\in\on{Star}(v)}s(\lambda,e)\right)v.\]
Notice that the dependence on $\lambda$ is only up to a global translation by $\mathbb R$. Divisors of the form $\on{div}(\lambda)$ are called \emph{principal}. Two divisors are \emph{linearly equivalent} if their difference is principal. We can define the \emph{divisor class group} of $\tropC$ by taking equivalence classes of divisors on $\tropC$ modulo linear equivalence. The \emph{linear system} of a divisor is given by:
\[\lvert D\rvert=\{D^\prime\sim_{\text{lin}}D|D^\prime\geq 0\}.\]
 We will loosely refer to a member of $\lvert K_{\tropC}\rvert$ as a tropical canonical divisor.
 
 Notice that a tropical linear system is \emph{tropically convex}, i.e. it is closed under taking $\on{max}$ and real translations (see \cite[Lemma 4]{HaaseMusikerYu}).
 
 Everything here can be restated for tropical curves over a more general base. Subdivisions, though, may only make sense after enlarging the base monoid, often corresponding to a subdivision of its dual cone.

If $\psi\colon\tropC\to\tropT$ is a harmonic morphism of tropical curves \cite{ABBR}, both divisors and piecewise-linear functions can be pulled back form $\tropT$ to $\tropC$ - in both cases, a coefficient accounting for the expansion factor of $\psi$ along edges must be included.

\subsubsection{Alignments and blow-ups}\label{sec:logblowups}
The monoid $\overline M$ induces a partial order on $\overline M^{\text{gp}}$: \[m_1\leq m_2 \qquad \Leftrightarrow \qquad m_2-m_1\in\overline M.\]
Given a logarithmic scheme $(\underline S,M_S)$ with a logarithmic ideal $K=(m_1\ldots,m_h)$, the logarithmic subfunctor of $S$ defined by requiring that there is always a minimum among the $m_i,\ i=1,\ldots,h,$ - i.e. that the ideal is locally principal - is represented by the blow-up of $S$ in the ideal $\alpha(K)$ \cite[\S 3.4]{Santos-Parker}. Tropically, it corresponds to a subdivision of the cone $\Hom(\overline{M}_S,\mathbb R_{\geq0})$, and vice versa.

This simple observation has had many fruitful applications in moduli theory, starting from \cite{MarcusWise}.

\subsection{Admissible covers and their logarithmic structure}\label{sec:adm}

Admissible covers have been introduced by J. Harris and D. Mumford in \cite{HarrisMumford} as a compactification of Hurwitz' spaces. A fully-fledged moduli theory for them has been developed only later by S. Mochizuki with the introduction of logarithmic techniques \cite{Mochizuki}, and by D. Abramovich, A. Corti, and A. Vistoli with the introduction of twisted curves \cite{ACV}. They have also been generalised by B. Kim in \cite[\S\S 5.2 and 7.2]{KimLog} (he calls them log stable $\mu$-ramified maps). We will only be concerned with the case of double covers of $\PP^1$, i.e. hyperelliptic curves.

In order to motivate the construction, we observe that every smooth curve of genus two is hyperelliptic (under the canonical map) with six Weierstrass points (i.e. ramification points for the hyperelliptic cover; this follows from the Riemann-Hurwitz formula). When the curve is allowed to become nodal, there still exists a degree two map to the projective line, but it is no longer finite. In this case, it is appropriate to define Weierstrass points to be limit of Weierstrass points in a smoothing family; the drawback is that whether a point on a (contracted) rational subcurve is Weierstrass or not may depend on the choice of the smoothing family. To resolve this ambiguity within the definition of the moduli problem, a solution is to expand the target, so as to consider the space of \emph{finite} morphisms from nodal curves of genus two to rational \emph{trees}. Here is a formal definition.

\begin{dfn}\label{def:admissiblecover}
 A family of \emph{admissible hyperelliptic covers} over $S$ is a finite morphism $\psi\colon (C,D_R)\to (T,D_B)$ over $S$ such that:
 \begin{enumerate}
  \item $(C,D_R)$ and $(T,D_B)$ are prestable curves with (unlabelled) smooth disjoint multisections $D_R$ and $D_B$ of length $2g+2$, $C$ has arithmetic genus $g$, and $(T,D_B)$ is a stable rational tree;
  \item $\psi$ is a double cover on an open $U\subseteq T$ dense over $S$;
  \item\label{point:logetale} $\psi$ is \'etale on $C^\text{sm}\setminus D_R$, it maps $D_R$ to $D_B$ with simple ramification, and it maps nodes of $C$ to nodes of $T$ so that in local\footnote{We refer the reader to \cite{Mochizuki}, and in particular Remark 2 of \S3.9, for a detailed discussion.} coordinates:
  \[ \psi^\#\colon \OO_S[u,v]/(uv-s)\to\OO_S[x,y]/(xy-t)\]
  maps $u\mapsto x^i,v\mapsto y^i, s\mapsto t^i$ for $i=1$ or $2$.
 \end{enumerate}
\end{dfn}
Mochizuki shows that condition \eqref{point:logetale} can be replaced by requiring that $\psi$ lifts to a \emph{logarithmically \'etale} morphisms of logarithmic schemes $(C,M_C)\to(T,M_T)$ over $(S,M_S)$, so that the image of a smoothing parameter of $T$ is either a smoothing parameter of $C$ or its double. Moreover there is a minimal logarithmic structure over $S$ that makes this possible: if $p$ is a node of $C$ that $\psi$ maps to the node $q$ of $T$, \[M_S^{\psi-\text{can}}=M_S^{C-\text{can}}\oplus_{\OO_S^*}M_S^{T-\text{can}}/\sim,\qquad (0,\delta_q)\sim(i\delta_p,0),\]
where $i$ is the local multiplicity defined at the end of \eqref{point:logetale}. In the case of double covers, this simply means that the minimal logarithmic structure for the admissible cover is the same as that of the source curve, except that the smoothing parameters of two nodes have been identified if $\psi$ matches them both with the same node of the target. We thus obtain a logarithmic Deligne-Mumford stack.
\begin{thm}[{\cite[\S 3.22]{Mochizuki}}] The moduli stack of admissible hyperelliptic covers $\mathcal A_{g,0,2}$ is a logarithmically smooth with locally free logarithmic structure (and therefore smooth), proper Deligne-Mumford stack, with a logarithmically \'etale morphism to $\overline{\pazocal M}_{0,2g+2}/\mathfrak S_{2g+2}$.
\end{thm}
To an admissible hyperelliptic cover we can associate a harmonic morphism of degree two to a metric tree $\psi\colon\tropC\to\tropT$. The tropical geometry of hyperelliptic and admissible covers has been analysed, for instance, in \cite{ChanHyperelliptic,CMR}.

\subsection{Gorenstein curves}
\subsubsection{Tools and formulae}
Curves shall always be assumed Cohen-Macaulay, i.e. $S1$; they might still be non-reduced along some subcurve, but they have no embedded points.

\begin{dfn}
 A curve $X$ is Gorenstein if its dualising sheaf $\omega_X$ is a line bundle.
\end{dfn}

A fundamental role in the study of singularities is played by the conductor ideal.
\begin{dfn}
 Let $\nu\colon\widetilde X\to X$ be a finite morphism. The \emph{conductor ideal} of $\nu$ is: \[\mathfrak c = \operatorname{Ann}_{\OO_X}\left(\nu_*\OO_{\widetilde X}/\OO_X\right)=\mathcal{H}om_{X}(\nu_*\OO_{\widetilde X},\OO_X).\]
\end{dfn}
The conductor is the largest $\OO_X-$ideal sheaf that is also a $\nu_*\OO_{\widetilde X}-$ideal sheaf. When $X$ is a reduced curve with finite normalisation - it is always the case over a field -, blowing up the conductor ideal recovers the normalisation \cite{Wilson}. The conductor ideal admits a further characterisation in terms of duality theory as:
\[\mathfrak c =\nu^!\OO_X=\omega_\nu.\]

\begin{rem}\label{rem:normalisation}
Let us restrict to schemes of finite type over a field. For such an $X$, the normalisation $X^\nu\to X$ has two universal properties (see e.g. \cite[\href{https://stacks.math.columbia.edu/tag/035Q}{Tag 035Q}]{stacks-project}):
\begin{enumerate}
 \item It is final in the category of dominant morphisms $Y\to X$ with $Y$ normal.
 \item\label{point:normalisation2} For every finite and birational morphisms $Y\to X$, the normalisation of $X$ factors uniquely through a morphism $X^\nu\to Y$.
\end{enumerate}
Now, suppose that $X$ is a Cohen-Macaulay curve. The smooth (possibly disconnected) curve $(X_\text{red})^\nu$ is what is normally referred to as the normalisation of $X$. Notice that $(X_\text{red})^\nu$ can be obtained as a (sequence of) blow-up(s) along a subscheme (points) $\delta$ of $X_\text{red}$. Set $X^\prime=\on{Bl}_{\delta}X$, and consider the following diagram:
\bcd
X & X^\prime\ar[l,"\nu^\prime" above]& (X^\prime_\text{red})^\nu\ar[l] \\
X_\text{red}\ar[u,hook] & (X_\text{red})^\nu\ar[l]\ar[ul,"f"]\ar[u,dashed,"\alpha"]\ar[ur,dashed,"\beta"] &
\ecd
The arrow $\alpha$ exists by the universal property of blowing up (because $f^{-1}\pazocal I_{\delta/X}$ is principal in $\OO_{(X_\text{red})^\nu}$). Since $(X_\text{red})^\nu$ is reduced and normal, $\alpha$ factors through $(X^\prime_\text{red})^\nu$, hence the arrow $\beta$ exists. Since it is a finite and birational morphism of normal curves, $\beta$ is an isomorphism. Abusing notation, we are going to call $X^\prime$ the \emph{normalisation} of $X$. Note that the underlying reduced curve of every irreducible component is regular. Moreover, the sheaf $\nu^\prime_*\OO_{X^\prime}/\OO_X$ is supported along the singular locus of $X_\text{red}$. Suppose that $X^\prime$ contains a ribbon (see Definition  \ref{def:ribbon}; or, more generally, a \emph{multiple curve} in the sense of \cite{Drezet-parametrisation}); the reader should be aware that, if we only ask of a morphism $g\colon Y\to X$ that the sheaf $g_*\OO_Y/\OO_X$ be supported at points, we can construct infinitely many more $Y$ by blowing up $X^\prime$ successively at a number of smooth points of $X_\text{red}$. Indeed, for ribbons, these are essentially all the morphisms that restrict to the identity on the underlying $\PP^1$, yet altering the square-zero ideal that defines the double structure (see \cite[Theorem 1.10]{Bayer-Eisenbud}).

\end{rem}

It is known that for a nodal curve $X$ the conductor ideal relates the dualising line bundle of $X$ to that of its normalisation $X^\nu$; more generally \cite[Proposition 1.2]{Catanese}:

\begin{prop}[Noether's formula]\label{prop:Noether}
 Let $\nu\colon\widetilde X\to X$ be a finite birational morphism of Gorenstein curves. Then, viewing  $\mathfrak c$ as an ideal sheaf on $\widetilde X$: \[\omega_{\widetilde X}=\nu^*\omega_X\otimes\mathfrak c.\]
\end{prop}

\begin{cor}\label{cor:principalconductor}
 Let $X$ be a Gorenstein curve with Gorenstein normalisation $\nu\colon\widetilde X\to X$. The conductor is a principal ideal on $\widetilde X$.
\end{cor}

\begin{dfn}
 For a finite birational morphism $\nu\colon\widetilde{X}\to X$ of curves, the coherent $\OO_X-$module $\nu_*\OO_{\widetilde X}/\OO_X$ has finite support; its length is called the \emph{$\delta$-invariant} of $\nu$. When $\nu$ is the normalisation, we simply call it the $\delta$-invariant of $X$. When $X$ is reduced, it is the sum of the local contributions of all the isolated singularities of $X$.
\end{dfn}

We review a result of J.P. Serre \cite[\S~4, Proposition~7]{serreAG} for possibly non-reduced curves.

\begin{lem}\label{lem:n=2delta}
 If $X$ is a Gorenstein curve and $\widetilde X$ its normalisation, then 
 \begin{equation}\label{eqn:Gorenstein}
  \operatorname{length}\left(\OO_{\widetilde X}/\mathfrak c\right)=2\delta.
 \end{equation}
\end{lem}
\begin{proof}
Since the normalisation is finite and birational, we may work locally around the support of $\nu_*\OO_{\widetilde X}/\OO_X$, which we may therefore assume is a single point. Let $X=\operatorname{Spec}(R)$ and $\widetilde X=\operatorname{Spec}(\widetilde R)$. Using $\mathfrak c\subseteq R\subseteq \widetilde R$ we see that \eqref{eqn:Gorenstein} is equivalent to:
\[\dim_{\mathbb C}(R/\mathfrak c)=\delta.\]

Let $\omega$ be a local generator of $\Omega_R$; given $f\in R$, Noether's formula implies that $f\omega\in \Omega_{\widetilde{R}}$ if and only if $f\in \mathfrak c$.
We therefore get: \[R/\mathfrak c\hookrightarrow \Omega_R/\Omega_{\widetilde{R}}.\] Since $\Omega_R$ is free on rank one and $\omega$ a generator, the above map is surjective as well. Finally, applying $\operatorname{Hom}_R(-,\Omega_R)$ to the normalisation exact sequence, we obtain $\Omega_R/\Omega_{\widetilde{R}}=\operatorname{Ext}^1_R(\widetilde{R}/R,\Omega_R)$, and therefore $\dim_{\mathbb C}\Omega_R/\Omega_{\widetilde{R}}=\delta$.
\end{proof}

\begin{rem} The converse is true when $X$ is reduced, but may fail otherwise; an example of which is provided by the subalgebra of $\mathbb C[\![s,\epsilon]\!]\times \mathbb C[\![t]\!]$ generated by $\langle s,\epsilon,t\rangle$, i.e. the transverse union of a ribbon with a line, which, though satisfying \eqref{eqn:Gorenstein}, is not Gorenstein as a consequence of the following lemma. What seems to be lacking in the non-reduced case is a good theory of dualising sheaves as rational differential forms on the normalisation satisfying some residue condition.
\end{rem}

\begin{dfn}
 A curve is \emph{decomposable} if it is obtained by gluing two curves along a (reduced closed) point.
\end{dfn}
\begin{rem}
 This is equivalent to the definition of \cite[Definition 2.1]{Stevens}. Indeed, the fiber product over $\mathbb C$ of $\mathbb C[x_1,\ldots,x_m]/I_1$ and $\mathbb C[y_1,\ldots,y_n]/I_2$ is isomorphic to
 \begin{equation}\label{eqn:deomposable}
  \underbrace{\mathbb C[x_1,\ldots,x_m,y_1,\ldots,y_n]}_S/\underbrace{(I_1(\mathbf x),y_1,\ldots,y_n)}_{J_1}\cap\underbrace{(x_1,\ldots,x_m,I_2(\mathbf y))}_{J_2}.
 \end{equation}
On the other hand, with notation as in \eqref{eqn:deomposable}, there is a short exact sequence:
\[0\to S/J_1\cap J_2\to S/J_1\times S/J_2\to S/J_1+J_2\simeq\mathbb C\to 0,\]
which is exact in the middle because every element there can be represented as a pair of polynomials in $\mathbf x$ only (resp. in $\mathbf y$ only), which is easy to lift.
\end{rem}

We review \cite[Proposition 2.1]{AFSGm} for not necessarily reduced curves.

\begin{lem}\label{lem:indecomposable}
 A decomposable curve may be Gorenstein only if it is a node.
\end{lem}
\begin{proof}
Let us assume that $X$ is the decomposable union of $X_1$ and $X_2$.

Notice first that if $X$ is Gorenstein, then so are $X_1$ and $X_2$.
Indeed, say we can find $(a_1,\ldots,a_n)\subseteq\mathfrak m_X$ a regular sequence such that $\OO_{X^0}=\OO_X/(a_1,\dots,a_n)$ is Gorenstein of dimension zero, i.e.  $\operatorname{Hom}_{X^0}(\mathbb C,\OO_{X^0})\cong \mathbb C$ (and all higher Ext groups vanish, which is automatic by \cite[Theorem~18.1]{matsumura1989}).
Since $\mathfrak m_X= \mathfrak m_{X_1}\oplus\mathfrak m_{X_2}$, this defines a regular sequence in the latter as well, and  $\OO_{X_1}\oplus \OO_{X_2}/(a_1,\dots,a_n)\cong \OO_{X_1^0}\oplus \OO_{X_2^0}$ is an extension of $k$ with $\OO_{X^0}$. From this we see that $\dim_{\mathbb C} \operatorname{Hom}_{X_i^0}(\mathbb C, \OO_{X_i^0})=1,\ i=1,2,$ which is to say that $X_1$ and $X_2$ are Gorenstein as well.

Now, $X$ and $X_1\sqcup X_2$ have the same normalisation $\widetilde X$, and
\begin{equation}\label{eqn:deltacomparison}
\delta_{X}=\delta_{X_1}+\delta_{X_2} +1.  
\end{equation}
On the other hand, $\OO_{\widetilde X}/\OO_X\twoheadrightarrow \OO_{\widetilde X}/\OO_{X_1\sqcup X_2}$ implies that $\mathfrak c_X\subseteq\mathfrak c_{X_1\sqcup X_2}$ (as ideals of $\OO_{\widetilde X}$).
We claim that, if $X$ is Gorenstein, the reduced curves underlying $X_1$ and $X_2$ are both smooth.
\begin{itemize}[leftmargin=.5cm]
 \item If they are both singular, then $\mathfrak c_{X_1\sqcup X_2}\subseteq \mathfrak m_{X_1}\oplus\mathfrak m_{X_2}=\mathfrak m_X$, so $\mathfrak c_X=\mathfrak c_{X_1\sqcup X_2}$. By Lemma \ref{lem:n=2delta}, we obtain:
 $2\delta_X=\dim_{\mathbb C}(\OO_{\widetilde X}/\mathfrak c_X)=\dim_{\mathbb C}(\OO_{\widetilde X}/\mathfrak c_{X_1\sqcup X_2})=2\delta_{X_1}+2\delta_{X_2}$, which contradicts \eqref{eqn:deltacomparison}.
 \item If $X_1$ is singular and $X_2$ smooth, \eqref{eqn:deltacomparison} reduces to $\delta_X=\delta_{X_1}+1$. On the other hand, $\mathfrak c_{X_2}=\OO_{X_2}$, and the contradiction comes from:
 $2\delta_X=\dim_{\mathbb C}(\OO_{\widetilde X}/\mathfrak c_X)=\dim_{\mathbb C}(\OO_{\widetilde X}/\mathfrak c_{X_1})=2\delta_{X_1}$.
\end{itemize}

Finally, since both underlying curves are smooth, we see from the definition of decomposability that:
\[\widehat{\OO}_X=\mathbb C[\![x,\epsilon_i,y,\eta_j]\!]/(xy,x\eta_j,\epsilon_i y,\epsilon_i\eta_j,\epsilon_i^{m_i},\eta_j^{n_j})_{\substack{i=1,\ldots,h\\j=1,\ldots,k}}.\]
If $X_1$ and $X_2$ are reduced, we recover the node. In all other cases, one can verify that $\dim_{\mathbb C}\operatorname{Hom}(\mathbb{C},\OO_X)\geq 2$ (it is generated by $xy,x\epsilon_1^{m_1-1},\ldots$).

\end{proof}

\subsubsection{Isolated singularities}

Let $(X,x)$ be (the germ of) a reduced curve with a unique singular point $x$, with normalisation $\nu\colon\widetilde{X}\to X$. The following is a measure of how much of the arithmetic genus of a projective curve is hiding in its singularities.
\begin{dfn}\label{def:genus}\cite{SMY1}
If $X$ has $m$ branches at $x$, the \emph{genus} of $(X,x)$ is:
\[g=\delta-m+1.\] 
\end{dfn}

The classification of isolated Gorenstein singularities of genus one has been carried out by D.I. Smyth in \cite[A.3]{SMY1}.

\begin{prop}
 An $(X,x)$ of genus one with $m$ branches is locally isomorphic to:
\begin{description}
 \item[$m=1$] the cusp, $V(y^2-x^3)\subseteq\Aaff^2_{x,y}$;
 \item[$m=2$] the tacnode, $V(y^2-yx^2)\subseteq\Aaff^2_{x,y}$;
 \item[$m\geq 3$] the union of $m$ general lines through the origin of $\Aaff^{m-1}$.
\end{description}
\end{prop}

All of these singularities are smoothable. Choosing a one-parameter smoothing and passing to a regular semistable model, the \emph{semistable tail}, i.e. the subcurve contracted to the singularity, admits a simple description: if we mark the semistable tail by the intersection with the rest of the central fibre, then it is a \emph{balanced} nodal curve of genus one. This means that it consists of a genus one core - which can be either smooth, or a circle of rational curves - together with some rational trees supporting the markings, and the distance between a marking and the core - i.e. the length of the corresponding rational chain - is independent of the chosen marking \cite[Proposition 2.12]{SMY1}.

\smallskip

The classification of isolated Gorenstein singularities of genus two has been carried out by the first author in \cite[\S2]{B}.

\begin{prop}
 The unique unibranch Gorenstein singularity of genus two is the \emph{ramphoidal cusp} or $A_4$-singularity, $V(y^2-x^5)\subseteq\Aaff^2_{x,y}$. For every $m\geq 2$, there are exactly two  isomorphism classes of germs of isolated Gorenstein singularities of genus two:
\end{prop}

 \begin{tabular}{c||c|c}
  & type I & type II \\
  \hline
  \hline
  Parametr. & 
  \begin{minipage}{.37\textwidth}{
 \begin{align*}
 \begin{split}
  x_1= & t_1\oplus0\oplus\ldots\oplus t_m^3\\
  x_2= & 0\oplus t_2\oplus\ldots\oplus t_m^3\\
  &\ldots\\
  x_{m-1}= & 0\oplus\ldots\oplus t_{m-1}\oplus t_m^3\\
  x_m= & 0\oplus\ldots\oplus0\oplus t_m^2 
  \end{split}
 \end{align*}}
  \end{minipage}
&
  \begin{minipage}{.39\textwidth}{
   \begin{align*}
 \begin{split}
  x_1= & t_1\oplus0\oplus\ldots\oplus t_m\\
  x_2= & 0\oplus t_2\oplus\ldots\oplus t_m^2\\
  &\ldots\\
  x_{m-1}= & 0\oplus\ldots\oplus t_{m-1}\oplus t_m^2\\
  (y= & 0\oplus t_2^3 \quad\text{if } m=2)
 \end{split}
 \end{align*}}
  \end{minipage}
\\
\hline
  Equations & 
  \begin{minipage}{.37\textwidth}{
  
  \begin{description}
  \item[{\textcolor{gray}{$m=1$}}] \textcolor{gray}{$x^5-y^2$ ($A_4$);}
  \item[$m=2$] $x_2(x_2^3-x_1^2)$ ($D_5$);
  \item[$m=3$] $\langle x_3(x_1-x_2),x_3^3-x_1x_2\rangle$;
  \item[$m\geq 4$] $\langle x_i(x_j-x_k),\\ x_m(x_i-x_j),\\ x_m^3-x_1x_2\rangle_{i\neq j\neq k\in\{1,\ldots,m-1\}}$
   \end{description}} 
  \end{minipage}
&
  \begin{minipage}{.39\textwidth}{
  
  \begin{description}
  \item[$m=2$] $y(y-x_1^3)$ ($A_5$);
  \item[$m=3$] $x_1x_2(x_2-x_1^2)$ ($D_6$);
  \item[$m\geq 4$] $\langle x_3(x_1^2-x_2),\\ x_i(x_j-x_k)\rangle_{\substack{1\leq i<j<k\leq m-1 \text{ or }\\1<j<k<i\leq m-1}}$
   \end{description}} 
  \end{minipage}
 \end{tabular}

\smallskip 
 
It is not hard to see that, for $m\geq 3$, every type I singularity is the union of a cusp with $m-1$ lines, and every type II singularity is the union of a tacnode with $m-2$ lines; in each case, we refer to the components of the genus one subcurve as the \emph{special branches}.

The description of semistable tails is a bit more cumbersome, but enlightening: assuming that the genus two core is smooth, the special branches are closer to the core and attached to a Weierstrass point (respectively, two conjugate points) in type I (resp. II); all other branches are equidistant and further away from the core. In fact, the ratio between the length of the special rational chain and the others is fixed to $\frac{1}{3}$ (resp. $\frac{1}{2}$) in type I (resp. II). For a more detailed statement in case the core is not smooth see \cite[Propositions 4.3-4.6]{B}.

\subsubsection{Non-reduced structures}

Multiple curves were investigated in the '90s in connection to Green's conjecture \cite{Bayer-Eisenbud}. Ribbons, in particular, were understood to arise as limits (in the Hilbert scheme of $\PP^{g-1}$) of canonical curves, as the curve becomes hyperelliptic \cite{Fong}.

\begin{dfn}\label{def:ribbon}
 A \emph{ribbon} is a double structure on $\PP^1$, i.e. a non-reduced curve $R$ with $R_{\text{red}}=\PP^1$ defined by a square-zero ideal $\pazocal I_{R_{\text{red}}/R}$ that is a line bundle on $R_{\text{red}}$.
\end{dfn}

\begin{exa}
 There is only one ribbon of genus two, $R_2$, up to isomorphism: it is the first infinitesimal neighbourhood of the zero section in $\operatorname{Tot}_{\PP^1}(\OO(3))$. The short exact sequence
 \[0\to \mathcal I_{\PP^1/R_2}\simeq\OO_{\PP^1}(-3)\to \OO_{R_2}\to\OO_{\PP^1}\to 0\]
 is split by restricting the projection $\operatorname{Tot}_{\PP^1}(\OO(3))\to\PP^1$ to $R_2$. $R_2$ is therefore called a \emph{split ribbon}. In fact, all ribbons of genus at most two are split.
\end{exa}

Automorphisms and moduli of multiple curves have been studied in \cite{Drezet-parametrisation}. We are going to encounter non-reduced structures along singular curves as well.

\begin{exa}\label{exa:1-tailedribbon}
 The simplest example of a Gorenstein non-reduced structure with singular underlying curve is given by the union of a ribbon and a line \emph{along a double point}. Local equations are $\mathbb C[\![x,y]\!]/(x^2y)$. It is easy to see that, for such a curve to have genus two, the ribbon needs to have ideal $\mathcal I_{\PP^1/R}\simeq\OO_{\PP^1}(-2)$. Such a curve can be realised in the linear system $\lvert 2D_+ + F\rvert$ on the Hirzebruch surface $\mathbb F_2$, where $D_+$ denotes the class of the positive section, and $F$ the class of a fiber.
\end{exa}

\begin{exa}\label{exa:tailedribbons}
 More generally, the union of a ribbon, $R$, with a rational $k$-fold point along a double point (representing a generic tangent vector to the $k$-fold point) is Gorenstein. Local equations are given by
 \[\mathbb C[\![x_1,\ldots,x_k,y]\!]/(x_ix_j,(x_i-x_j)y)_{1\leq i <j\leq k}.\]
 To see that this is Gorenstein, it is enough to find a regular element $\xi$ such that the quotient be Gorenstein of dimension zero. Let $\xi=\sum_{i=1}^kx_i-y$. The quotient $A$ is a graded finite-dimensional algebra, with $A_0=\mathbb C$, $A_1=\mathbb C\langle x_1,\ldots,x_k\rangle$, and $A_2=\mathbb C\langle x_1^2=\ldots=x_k^2\rangle$, having one-dimensional socle $A_2$.
 
 Now, we can obtain a Gorenstein projective curve of genus two, $C$, by gluing a ribbon together with some $k_i$-fold points ($i=1,\ldots,r$) at distinct (closed) points of the ribbon, by iterating the local construction above. From the short exact sequence
 \[0\to \OO_C\to \OO_R\oplus\bigoplus_{i=1}^r\OO_{\PP^1}^{\oplus k_i}\to \bigoplus_{i=1}^r(\mathbb C^{k_i-1}\oplus\mathbb C[\epsilon])\to 0,\]
 we can compute the Euler characteristic of $R$, whose structure sheaf satisfies thus
 \[0\to\OO_{\PP^1}(r-3)\to\OO_R\to\OO_{\PP^1}\to0,\]
 depending only on the total number of ``noded'' points $r$, and not on the number of branches of each $k_i$-fold point.
\end{exa}
\begin{dfn}
 We call $C$ as in Example \ref{exa:tailedribbons} a \emph{$(k_1,\ldots,k_r)$-tailed ribbon} of genus two.
\end{dfn}

\begin{rem}\label{rem:Noetheradjunction}
 We can employ Noether's formula (Proposition \ref{prop:Noether}) and adjunction (on a surface containing the ribbon, i.e. $\operatorname{Tot}(\OO_{\PP^1}(3-r))$) to compute the restriction of the dualising sheaf to every component of a $(k_1,\ldots,k_r)$-tailed ribbon $C$. First, the ``normalisation'' is given by:
 \begin{flalign*}
  \mathbb C[\![x_1,\ldots,x_k,y]\!]/(x_ix_j,(x_i-x_j)y)_{1\leq i <j\leq k} &\to&  \mathbb C[\![s]\!][\epsilon]/(\epsilon^2)\times \mathbb C[\![t_1]\!]\times\ldots\times \mathbb C[\![t_k]\!] &\\
  x_i &\mapsto&(\epsilon,0,\ldots,t_i,\ldots,0),i=1,\ldots,k\\
  y&\mapsto&(s,0,\ldots,0)\qquad
 \end{flalign*}
From this we compute the conductor $\mathfrak c=\langle x_1^2,\ldots,x_k^2,y\rangle$. The restriction to a tail is:
\[\omega_{C|T}\simeq\omega_T\otimes\OO(2q)\simeq\OO_T.\]
The restriction to the $\PP^1$ underlying the ribbon is:
\[\omega_{C|R_\text{red}}\simeq \omega_{R_\text{red}}\otimes \mathcal I_{R_\text{red}/R}^\vee\otimes\mathfrak c_{R_\text{red}}^\vee\simeq\OO_{R_\text{red}}(-2+(3-r)+r)\simeq\OO_{R_\text{red}}(1).\]
\end{rem}

For a nodal curve $C$ with core $Z$, we say that a component $D$ of $C$ \emph{cleaves} to $q\in Z$ if $D$ is joined to $Z$ at $q$ through a chain of rational curves.

\begin{lem}\label{lem:ssmodelsofribbons}
 Regular semistable models of a $(k_1,\ldots,k_r)$-tailed ribbon can be classified. If the core is a smooth curve of genus two $Z$, the sets of tails $\{T_1^i,\ldots,T_{k_i}^i\}_{i=1,\ldots,r}$ cleave to the same (or to conjugate) point(s) of $Z$, and are all equidistant from the core (independently of $i$ and $j$ in $T^i_j$). The configuration of attaching points on $Z$ is such that the hyperelliptic cover maps it to the corresponding configuration of noded points on $R_\text{red}$, up to reparametrisation.
\end{lem}

\begin{proof}
 A word-by-word repetition of the argument of \cite[Proposition 4.3]{B}.
\end{proof}

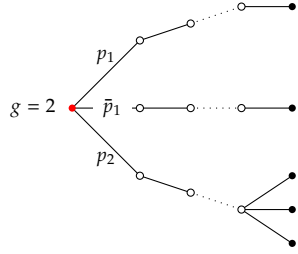
\begin{figure}[htb]
	\centering
	\begin{minipage}[t]{0.5\textwidth}
		\centering
		\begin{tikzpicture} [scale=.45]
		\tikzstyle{every node}=[font=\normalsize]
		\tikzset{arrow/.style={latex-latex}}
		\def\S{2.825cm} 
	\coordinate (O) at (0,0,0);
		  \coordinate (A) at (2,2,0);
		 \coordinate (B) at (2,0,0);
		 \coordinate (C) at (2,-2,0);
		 
		  \coordinate (A1) at (3.5,2.5,0);
		 \coordinate (B1) at (3.5,0,0);
		 \coordinate (C1) at (3.5,-2.5,0);
		 
		  \coordinate (A2) at (5,3,0);
		 \coordinate (B2) at (5,0,0);
		 \coordinate (C2) at (5,-3,0);
		 
		 \coordinate (a) at (6.5,3,0);
		 \coordinate (b) at (6.5,0,0);
		 \coordinate (c1) at (6.5,-2,0);
		 \coordinate (c2) at (6.5,-3,0);
		 \coordinate (c3) at (6.5,-4,0);

    \node at (-0.2,0,0) [left] {\tiny{$g=2$}};

	\draw (O)--node[above]{\tiny{$p_1$}}(A);
	\draw (O)--node[fill=white,right=-.18cm]{\tiny{$\bar{p}_1$}}(B);
	\draw (O)--node[below]{\tiny{$p_2$}}(C);
	
	\draw (A1)--(A);
	\draw[dotted] (A1)--(A2);
	\draw (A2)--(a);	 
	
	\draw (B1)--(B);
	\draw[dotted] (B1)--(B2);
	\draw (B2)--(b);

		\draw (C1)--(C);
	\draw[dotted] (C1)--(C2);
	\draw (C2)--(c1);	 
	\draw (C2)--(c2);
	\draw (C2)--(c3);
	
	\foreach \x in {a,b,c1,c2,c3}
   \fill (\x) circle (3pt);
    \fill[red] (O) circle (3pt);
 \foreach \x in {A,B,C,A1,A2,B1,B2,C1,C2}
   \draw[fill=white] (\x) circle (3pt);
     	
\end{tikzpicture}
	\end{minipage}
\caption{Dual graph of the central fibre of a regular semistable model for a $(2,3)$-tailed ribbon with $p_1,\overline{p}_1$ conjugate.}
\label{fig:exa_ribbon_ss_tail}
\end{figure}

See Figure \ref{fig:exa_ribbon_ss_tail} for an example of a semistable tail of a tailed ribbon.

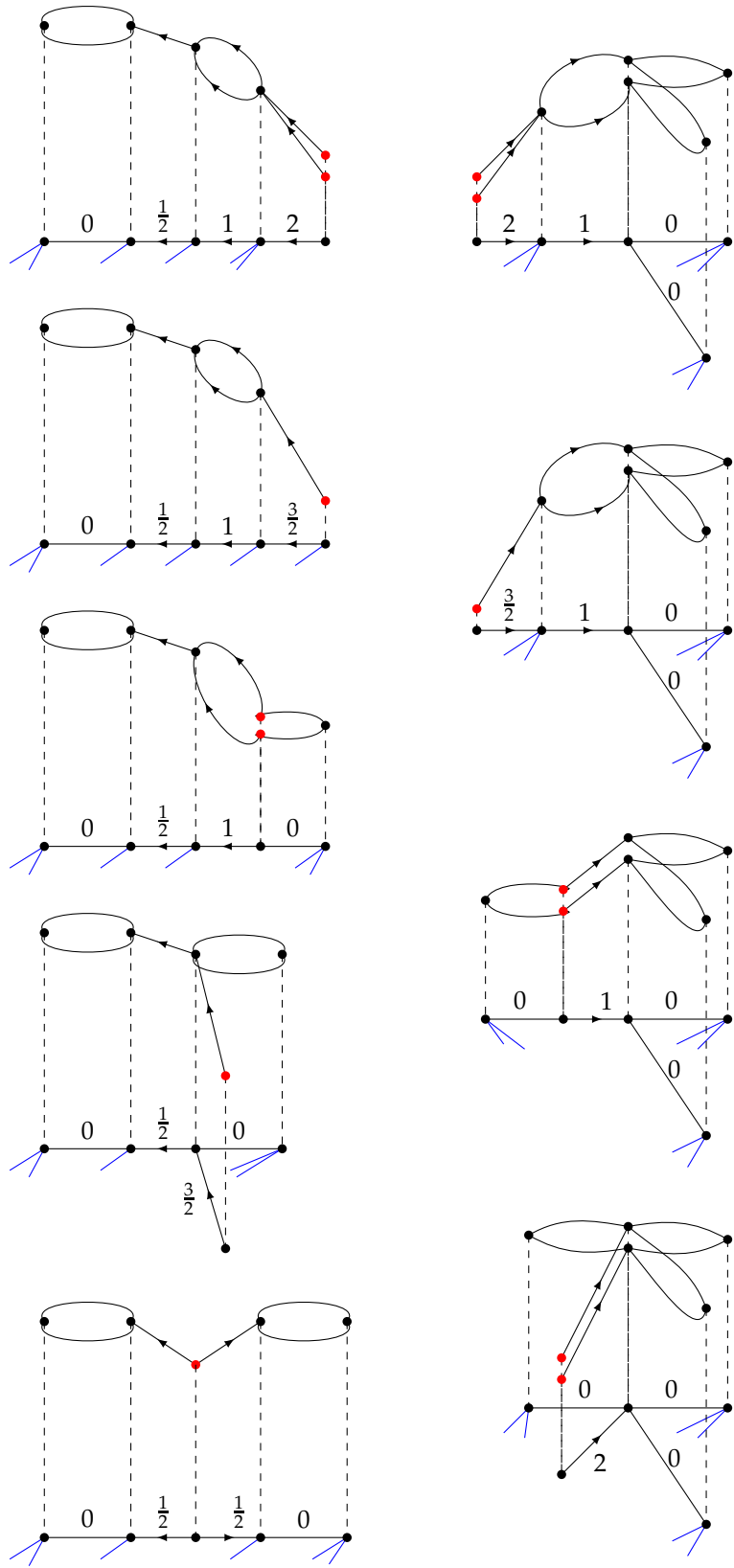
\begin{figure}
	\centering
		\begin{tikzpicture} [scale=.6]
		\tikzstyle{every node}=[font=\normalsize]
		\tikzset{arrow/.style={latex-latex}}
		\def\S{2.825cm} 

		
		\coordinate (A) at (0,6,0);
		  \coordinate (B) at (2,6,0);
		 \coordinate (C) at (3.5,5.5,0);
		 
		  \coordinate (CM) at (3.5,5,0);;
		 	\coordinate (AM) at (5,6,0);
		  \coordinate (BM) at (7,6,0);;

		 \coordinate (D) at  (5,4.5,0); 
		 \coordinate (Dd) at (5.5,5.5,0); 
		\coordinate (D1) at  (5,4,0);
		\coordinate (D2) at  (5,3.6,0);
		
		\coordinate (Ew) at (6.5,5,6);
		\coordinate (Ee) at (6.5,3.8,0);
		\coordinate (E) at (6.5,2,0); 
         \coordinate (E1) at (6.5,2.5,0); 
         \coordinate (E2) at ( 6.5,3,0); 
      	\coordinate (At) at (0,1,0);
		  \coordinate (Bt) at (2,1,0);
		 \coordinate (Ct) at (3.5,1,0);
		\coordinate (Dt) at  (5,1,0);
		\coordinate (Et) at (6.5,1,0); 
		\coordinate (Ddt) at (5.5,1,0); 
		\coordinate (Ewt) at (6.5,1,6);); 
		\coordinate (BMt) at (7,1,0);;

		\coordinate (P1) at (-.3,1,1.3);
		\coordinate (P2) at (.3,1,1.7);
		\coordinate (P3) at (1.8,1,1.3);
		\coordinate (P4) at (3.3,1,1.3);
		\coordinate (P5) at (4.8,1,1.3);
		\coordinate (P6) at (6.3,1,1.3);
		\coordinate (Pp6) at (6.8,1,1.7); 
		\coordinate (Q6) at (5.1,1,1.7); 
        \coordinate (T6) at (6.5,1,6);
        
        \coordinate (QM6) at (6.7,1,1.3); 
        \coordinate (TM6) at (7.2,1,1.6); 
      
      
      \draw (A) to [out=130,in=50] (B);
      \draw (A) to [out=-130,in=-50] (B);
      \draw[->-] (C)--(B);
      \draw[->-] (D) to [out=70,in=70](C);
      \draw[->-] (D) to [out=-110,in=-110](C);
       \draw[->-] (E)--(D);
         
        \draw[->-] (Et)--(Dt);
        \draw[->-] (Dt)--(Ct);
        \draw[->-] (Ct)--(Bt);
        \draw (Bt)--(At);

		\draw[blue] (At)--(P1);
		\draw[blue] (At)--(P2);
		\draw[blue] (Bt)--(P3);
		\draw[blue] (Ct)--(P4);
		\draw[blue] (Dt)--(P5);
		\draw[blue] (Et)--(P6);
		
\node at (1,1,0) [above] {\small{0}};
\node at (2.75,1,0) [above] {\small{$\frac{1}{2}$}};
\node at (4.25,1,0) [above] {\small{$1$}};
\node at (5.75,1,0) [above] {\small{$\frac{3}{2}$}};

       
       \draw[dashed] (A)--(At);
        \draw[dashed] (B)--(Bt);
       \draw[dashed] (C)--(Ct);
       \draw[dashed] (D)--(Dt);
        \draw[dashed] (E)--(Et);

      \foreach \x in {A,B,C,D,At,Bt,Ct,Dt,Et}
   \fill (\x) circle (3pt);
     	\fill[red] (E) circle (3pt);



\begin{scope}[every coordinate/.style={shift={(0,7,0)}}]

       \draw ([c]A) to [out=130,in=50] ([c]B);
      \draw ([c]A) to [out=-130,in=-50] ([c]B);
       \draw[->-] ([c]C)--([c]B);
      \draw[->-] ([c]D) to [out=70,in=70]([c]C);
      \draw[->-] ([c]D) to [out=-110,in=-110]([c]C);
       \draw[->-] ([c]E1)--([c]D);
       \draw[->-] ([c]E2)--([c]D);
         
        \draw[->-] ([c]Et)--([c]Dt);
        \draw[->-] ([c]Dt)--([c]Ct);
        \draw[->-] ([c]Ct)--([c]Bt);
        \draw ([c]Bt)--([c]At);

		\draw[blue] ([c]At)--([c]P1);
		\draw[blue] ([c]At)--([c]P2);
		\draw[blue] ([c]Bt)--([c]P3);
		\draw[blue] ([c]Ct)--([c]P4);
		\draw[blue] ([c]Dt)--([c]P5);
		\draw[blue] ([c]Dt)--([c]Q6);
		
\node at (1,8,0) [above] {\small{0}};
\node at (2.75,8,0) [above] {\small{$\frac{1}{2}$}};
\node at (4.25,8,0) [above] {\small{$1$}};
\node at (5.75,8,0) [above] {\small{$2$}};

       
       \draw[dashed] ([c]A)--([c]At);
        \draw[dashed] ([c]B)--([c]Bt);
       \draw[dashed] ([c]C)--([c]Ct);
       \draw[dashed] ([c]D)--([c]Dt);
        \draw[dashed] ([c]E1)--([c]Et);
        \draw[dashed] ([c]E2)--([c]Et);

      \foreach \x in {A,B,C,D,At,Bt,Ct,Dt,Et}
   \fill ([c]\x) circle (3pt);
     	\fill[red] ([c]E1) circle (3pt);
     		\fill[red] ([c]E2) circle (3pt);
\end{scope}


\begin{scope}[every coordinate/.style={shift={(0,-7,0)}}]

     \draw ([c]A) to [out=130,in=50] ([c]B);
      \draw ([c]A) to [out=-130,in=-50] ([c]B);
      \draw[->-] ([c]C)--([c]B);
      \draw[->-] ([c]D1) to [out=70,in=70]([c]C);
      \draw[->-] ([c]D2) to [out=-110,in=-110]([c]C);
      \draw ([c]Ee) to [out=110,in=180]([c]D1);
      \draw ([c]Ee) to [out=-110,in=-180]([c]D2);

        \draw ([c]Et)--([c]Dt);
        \draw[->-] ([c]Dt)--([c]Ct);
        \draw[->-] ([c]Ct)--([c]Bt);
        \draw ([c]Bt)--([c]At);

		\draw[blue] ([c]At)--([c]P1);
		\draw[blue] ([c]At)--([c]P2);
		\draw[blue] ([c]Bt)--([c]P3);
		\draw[blue] ([c]Ct)--([c]P4);
		\draw[blue] ([c]Et)--([c]P6);
	    \draw[blue] ([c]Et)--([c]Pp6);
		
\node at (1,-6,0) [above] {\small{0}};
\node at (2.75,-6,0) [above] {\small{$\frac{1}{2}$}};
\node at (4.25,-6,0) [above] {\small{$1$}};
\node at (5.75,-6,0) [above] {\small{$0$}};

       
       \draw[dashed] ([c]A)--([c]At);
        \draw[dashed] ([c]B)--([c]Bt);
       \draw[dashed] ([c]C)--([c]Ct);
       \draw[dashed] ([c]D1)--([c]Dt);
        \draw[dashed] ([c]D2)--([c]Dt);
        \draw[dashed] ([c]Ee)--([c]Et);

      \foreach \x in {A,B,C,Ee,At,Bt,Ct,Dt,Et}
   \fill ([c]\x) circle (3pt);
     	\fill[red] ([c]D1) circle (3pt);
     		\fill[red] ([c]D2) circle (3pt);
\end{scope}


\begin{scope}[every coordinate/.style={shift={(0,-14,0)}}]

     \draw ([c]A) to [out=130,in=50] ([c]B);
      \draw ([c]A) to [out=-130,in=-50] ([c]B);
    \draw[->-] ([c]C)--([c]B);
      \draw ([c]C) to [out=130,in=50]([c]Dd);
      \draw([c]C) to [out=-130,in=-50]([c]Dd);
      \draw[->-] ([c]Ew) -- ([c]C);

        \draw [->-] ([c]Ewt)--([c]Ct);
        \draw ([c]Ddt)--([c]Ct);
        \draw[->-] ([c]Ct)--([c]Bt);
        \draw ([c]Bt)--([c]At);

		\draw[blue] ([c]At)--([c]P1);
		\draw[blue] ([c]At)--([c]P2);
		\draw[blue] ([c]Bt)--([c]P3);
		\draw[blue] ([c]Ddt)--([c]P5);
		\draw[blue] ([c]Ddt)--([c]Q6);
		\draw[blue] ([c]Ewt)--([c]T6);

\node at (1,-13,0) [above] {\small{0}};
\node at (2.75,-13,0) [above] {\small{$\frac{1}{2}$}};
\node at (4.5,-13,0) [above] {\small{$0$}};
\node at (4.9,-13,3) [left] {\small{$\frac{3}{2}$}};

       
       \draw[dashed] ([c]A)--([c]At);
        \draw[dashed] ([c]B)--([c]Bt);
       \draw[dashed] ([c]C)--([c]Ct);
       \draw[dashed] ([c]Dd)--([c]Ddt);
     \draw[dashed] ([c]Ew)--([c]Ewt);

      \foreach \x in {A,B,C,Dd,At,Bt,Ct,Ddt,Ewt}
   \fill ([c]\x) circle (3pt);
     	\fill[red] ([c]Ew) circle (3pt);
     		
\end{scope}

 
 
 \begin{scope}[every coordinate/.style={shift={(0,-23,0)}}]

       \draw ([c]A) to [out=130,in=50] ([c]B);
      \draw ([c]A) to [out=-130,in=-50] ([c]B);
       \draw[->-] ([c]CM)--([c]B);
        \draw[->-] ([c]CM)--([c]AM);
      \draw ([c]AM) to [out=130,in=50] ([c]BM);
      \draw ([c]AM) to [out=-130,in=-50] ([c]BM);
      
        \draw ([c]At)--([c]Bt);
        \draw[->-] ([c]Ct)--([c]Dt);
        \draw[->-] ([c]Ct)--([c]Bt);
        \draw ([c]Dt)--([c]BMt);

		\draw[blue] ([c]At)--([c]P1);
		\draw[blue] ([c]At)--([c]P2);
		\draw[blue] ([c]Bt)--([c]P3);
		\draw[blue] ([c]Dt)--([c]P5);
		\draw[blue] ([c]BMt)--([c]QM6);
		\draw[blue] ([c]BMt)--([c]TM6);

        \foreach \x in {A,B,AM,BM,At,Bt,Ct,Dt,BMt}
   \fill ([c]\x) circle (3pt);
     	\fill[red] ([c]CM) circle (3pt);
                
       
       \draw[dashed] ([c]A)--([c]At);
        \draw[dashed] ([c]B)--([c]Bt);
       \draw[dashed] ([c]CM)--([c]Ct);
       \draw[dashed] ([c]AM)--([c]Dt);
     \draw[dashed] ([c]BM)--([c]BMt);

\node at (1,-22,0) [above] {\small{0}};
\node at (2.75,-22,0) [above] {\small{$\frac{1}{2}$}};
\node at (4.5,-22,0) [above] {\small{$\frac{1}{2}$}};
\node at (6,-22,0) [above] {\small{$0$}};

    \end{scope}



\coordinate (Cc) at (0,0.5,0); 
\coordinate (C1) at (0,1,0);
\coordinate (C2) at (0,1.5,0);
\coordinate (A1) at (1.5,3,0);

\coordinate (G) at (0.2,2.75,0);
\coordinate (F1) at (2,3,0);
\coordinate(F2) at (2,2.5,0);
\coordinate (A14) at (1.2,4,0);
\coordinate (H1) at (3.5,2.7,4);
\coordinate (H2) at (3.5,2.2,4);

\coordinate (B1) at (3.5,4.2,0);
\coordinate (B2) at (3.5,3.7,0);
\coordinate (A2) at (8,5,7);
\coordinate (A3) at (5.8,3.9,0);


\coordinate (C1t) at (0,0,0);
\coordinate (Gt) at (0.2,0,0);
\coordinate (Ft) at (2,0,0);


\coordinate (A14t) at (1.2,0,0);
\coordinate (H1t) at (3.5,0,4);

\coordinate (A1t) at (1.5,0,0);
\coordinate (B1t) at (3.5,0,0);
\coordinate (A2t) at (8,0,7);
\coordinate (A3t) at (5.8,0,0); 

\coordinate (Q1) at ( 1.2,0,1.5);
\coordinate (Q2) at ( 1.8,0,1.8);
\coordinate (Q3) at ( 8.3,0,8.9);
\coordinate (Q4) at (7.6,0,8);
\coordinate (Q5) at ( 5.2,0,1.5);
\coordinate (Q61) at ( 5.8,0,1.8);

\begin{scope}[every coordinate/.style={shift={(10,8,0)}}]


\draw[->-] ([c]C1)--([c]A1);
\draw[->-] ([c]C2)--([c]A1);
\draw[->-] ([c]A1) to [out=110,in=150] ([c]B1);
\draw[->-] ([c]A1) to [out=-70,in=-70] ([c]B2);
\draw ([c]B1) to [out=10,in=150] ([c]A3);
\draw ([c]B2) to [out=-10,in=-150] ([c]A3);
\draw ([c]B1) to [out=-40,in=110] ([c]A2);
\draw ([c]B2) to [out=-50,in=-100] ([c]A2);
\draw[->-] ([c]C1t)--([c]A1t);
\draw[->-] ([c]A1t)--([c]B1t);
\draw ([c]B1t)--([c]A2t);
\draw ([c]B1t)--([c]A3t);


		\draw[blue] ([c]A1t)--([c]Q1);
		\draw[blue] ([c]A1t)--([c]Q2);
		\draw[blue] ([c]A2t)--([c]Q3);
		\draw[blue] ([c]A2t)--([c]Q4);
		\draw[blue] ([c]A3t)--([c]Q5);
		\draw[blue] ([c]A3t)--([c]Q61);

\foreach \x in {A1,A2,A3,B1,B2,C1t,A1t,B1t,A2t,A3t}
   \fill ([c]\x) circle (3pt);

\node at (10.75,8,0) [above] {\small{2}};
\node at (12.5,8,0) [above] {\small{$1$}};
\node at (15.35,8,3) [right] {\small{$0$}};
\node at (14.5,8,0) [above] {\small{$0$}};

       
       \draw[dashed] ([c]C1)--([c]C1t);
        \draw[dashed] ([c]C2)--([c]C1t);
       \draw[dashed] ([c]A1)--([c]A1t);
       \draw[dashed] ([c]B1)--([c]B1t);
        \draw[dashed] ([c]B2)--([c]B1t);
        \draw[dashed] ([c]A2)--([c]A2t);
        \draw[dashed] ([c]A3)--([c]A3t);

  \fill[red] ([c]C1) circle (3pt);
    \fill[red] ([c]C2) circle (3pt);

\end{scope}

\begin{scope}[every coordinate/.style={shift={(10,-1,0)}}]


\draw[->-] ([c]Cc)--([c]A1);
\draw[->-] ([c]A1) to [out=110,in=150] ([c]B1);
\draw[->-] ([c]A1) to [out=-70,in=-70] ([c]B2);

\draw ([c]B1) to [out=10,in=150] ([c]A3);
\draw ([c]B2) to [out=-10,in=-150] ([c]A3);
\draw ([c]B1) to [out=-40,in=110] ([c]A2);
\draw ([c]B2) to [out=-50,in=-100] ([c]A2);
\draw[->-] ([c]C1t)--([c]A1t);
\draw[->-] ([c]A1t)--([c]B1t);
\draw ([c]B1t)--([c]A2t);
\draw ([c]B1t)--([c]A3t);


		\draw[blue] ([c]A1t)--([c]Q1);
		\draw[blue] ([c]A1t)--([c]Q2);
		\draw[blue] ([c]A2t)--([c]Q3);
		\draw[blue] ([c]A2t)--([c]Q4);
		\draw[blue] ([c]A3t)--([c]Q5);
		\draw[blue] ([c]A3t)--([c]Q61);

\foreach \x in {A1,A2,A3,B1,B2,C1t,A1t,B1t,A2t,A3t}
   \fill ([c]\x) circle (3pt);
  

\node at (10.75,-1,0) [above] {\small{$\frac{3}{2}$}};
\node at (12.5,-1,0) [above] {\small{$1$}};
\node at (15.35,-1,3) [right] {\small{$0$}};
\node at (14.5,-1,0) [above] {\small{$0$}};

       
       \draw[dashed] ([c]Cc)--([c]C1t);
       \draw[dashed] ([c]A1)--([c]A1t);
       \draw[dashed] ([c]B1)--([c]B1t);
        \draw[dashed] ([c]B2)--([c]B1t);
        \draw[dashed] ([c]A2)--([c]A2t);
        \draw[dashed] ([c]A3)--([c]A3t);
        
         \fill[red] ([c]Cc) circle (3pt);
\end{scope}

\begin{scope}[every coordinate/.style={shift={(10,-10,0)}}]


\draw ([c]F1) to [out=0,in=70] ([c]G);
\draw ([c]F2) to [out=0,in=-70] ([c]G);
\draw[->-] ([c]F1) to  ([c]B1);
\draw[->-] ([c]F2) to  ([c]B2);
\draw ([c]B1) to [out=10,in=150] ([c]A3);
\draw ([c]B2) to [out=-10,in=-150] ([c]A3);

\draw ([c]B1) to [out=-40,in=110] ([c]A2);
\draw ([c]B2) to [out=-50,in=-100] ([c]A2);

\draw ([c]Gt)--([c]Ft);
\draw[->-] ([c]Ft)--([c]B1t);
\draw ([c]B1t)--([c]A2t);
\draw ([c]B1t)--([c]A3t);


		\draw[blue] ([c]Gt)--([c]Q1);
		\draw[blue] ([c]Gt)--([c]Q2);
		\draw[blue] ([c]A2t)--([c]Q3);
		\draw[blue] ([c]A2t)--([c]Q4);
		\draw[blue] ([c]A3t)--([c]Q5);
		\draw[blue] ([c]A3t)--([c]Q61);

\foreach \x in {G,A2,A3,B1,B2,Gt,Ft,B1t,A2t,A3t}
   \fill ([c]\x) circle (3pt);

\node at (11,-10,0) [above] {\small{$0$}};
\node at (13,-10,0) [above] {\small{$1$}};
\node at (15.35,-10,3) [right] {\small{$0$}};
\node at (14.5,-10,0) [above] {\small{$0$}};

       
       \draw[dashed] ([c]G)--([c]Gt);
       \draw[dashed] ([c]F1)--([c]Ft);
       \draw[dashed] ([c]F2)--([c]Ft);
        \draw[dashed] ([c]B2)--([c]B1t);
        \draw[dashed] ([c]A2)--([c]A2t);
        \draw[dashed] ([c]A3)--([c]A3t);

   \fill[red] ([c]F1) circle (3pt);
   \fill[red] ([c]F2) circle (3pt);
\end{scope}

\begin{scope}[every coordinate/.style={shift={(10,-19,0)}}]


\draw ([c]B1) to [out=170,in=30] ([c]A14);
\draw ([c]B2) to [out=190,in=-30] ([c]A14);

\draw[->-] ([c]H1) to ([c]B1);
\draw[->-] ([c]H2) to ([c]B2);
\draw ([c]B1) to [out=10,in=150] ([c]A3);
\draw ([c]B2) to [out=-10,in=-150] ([c]A3);

\draw ([c]B1) to [out=-40,in=110] ([c]A2);
\draw ([c]B2) to [out=-50,in=-100] ([c]A2);

\draw ([c]A14t)--([c]B1t);
\draw[->-] ([c]H1t)--([c]B1t);
\draw ([c]B1t)--([c]A2t);
\draw ([c]B1t)--([c]A3t);


		\draw[blue] ([c]A14t)--([c]Q1);
		\draw[blue] ([c]A14t)--([c]Q2);
		\draw[blue] ([c]A2t)--([c]Q3);
		\draw[blue] ([c]A2t)--([c]Q4);
		\draw[blue] ([c]A3t)--([c]Q5);
		\draw[blue] ([c]A3t)--([c]Q61);

\foreach \x in {A14,A2,A3,B1,B2,H1t, A14t,B1t,A2t,A3t}
   \fill ([c]\x) circle (3pt);
  

\node at (12.5,-19,0) [above] {\small{$0$}};
\node at (13.25,-19.5,2) [right] {\small{$2$}};
\node at (15.35,-19,3) [right] {\small{$0$}};
\node at (14.5,-19,0) [above] {\small{$0$}};

       
       \draw[dashed] ([c]A14)--([c]A14t);
       \draw[dashed] ([c]H1)--([c]H1t);
       \draw[dashed] ([c]H2)--([c]H1t);
     \draw[dashed] ([c]B1)--([c]B1t);
        \draw[dashed] ([c]B2)--([c]B1t);
        \draw[dashed] ([c]A2)--([c]A2t);
        \draw[dashed] ([c]A3)--([c]A3t);
        
 \fill[red] ([c]H1) circle (3pt);
   \fill[red] ([c]H2) circle (3pt);
\end{scope}
\end{tikzpicture}
\caption{Admissible functions on maximally degenerated tropical hyperelliptic curves: the dumbbell (l), and the theta graph (r).}
\label{fig:admredcores}
\end{figure}

\begin{rem}
  The description of the semistable models of a $(k_1,\ldots,k_r)$-tailed ribbon (Lemma \ref{lem:ssmodelsofribbons}) in case the former has reducible core is more cumbersome. Indeed, the length of a tree depends on which component of the core it is attached to. See Proposition \ref{prop:tropical_ss_tails} and Figure \ref{fig:admredcores} below for a precise description in terms of piecewise-linear functions on the tropicalization. A ribbon arises when the special component (that we denote pictorially by a red dot) belongs to the core. The piecewise-linear function has then constant slope $1$ on any tree emanating from the core, and in this sense it is reminiscent of the genus one situation. 
  
\end{rem}

\begin{exa}
 A semistable model of the $1$-tailed ribbon $C$ (Example \ref{exa:1-tailedribbon}) can be computed by taking the pencil \[u_{+}^2f_1-tu_{-}^2p_5(f_1,f_2)=0\] on $\mathbb F_2\times\Aaff^1_t$, with $p_5$ a generic homogeneous polynomial of degree $5$ in two variables. The resulting smoothing family is singular at six points along the central fiber, including the node of $C$. After blowing them up and normalising, we obtain a genus two curve covering the ribbon two-to-one; the tail of $C$ is attached to a Weierstrass point. This appears to be an accident due to the restriction that the smoothing comes from a pencil on $\mathbb F_2$.
\end{exa}


In order to describe semistable tails in complete generality, there is no better way than expressing the problem in terms of existence of a certain piecewise linear function on the dual graph. We refer the reader to Sections \ref{sec:PLfun} and \ref{sec:adm} for background and notation.

\begin{prop}\label{prop:tropical_ss_tails}
 Let $C$ be a Gorenstein curve of arithmetic genus two. Let $\mathcal C\to\dvr$ be a one-parameter smoothing of $C$. Let $\phi\colon \mathcal C^{\text{ss}}\to\mathcal C$ be a semistable model with exceptional locus $Z$. Then, up to destabilising $\mathcal C^{\text{ss}}$, there exists a hyperelliptic admissible cover $\psi\colon \mathcal C^{\text{ss}}\to \mathcal T$, and a tropical canonical divisor $\lambda\colon\tropC\to \mathbb R$ pulled back along $\on{trop}(\psi)$, such that $Z$ is cut out by $\OO_{\mathcal C^{\text{ss}}}(-\lambda)$, and we have $\omega_{\mathcal C^{\text{ss}}}=\phi^*\omega_{\mathcal C}(\lambda)$.
\end{prop}
\begin{proof}
 The proof in the isolated case has appeared in \cite[\S 4.4]{B}. The argument is insensitive to rational tails away from the core, so we may as well assume that $C$ is minimal (see Definition \ref{def:minimal} below). In this case, the canonical linear series gives a two-to-one map $\bar{\psi}\colon \mathcal C\to \PP_\dvr(\pi_*\omega_{\mathcal C})=:\PP$. The general fibre is the hyperelliptic cover of a smooth curve of genus two. Applying semistable reduction to the target and branch divisor, we may lift $\bar{\psi}$ to a map of nodal curves, thanks to the properness of the moduli space of admissible covers:
 \bcd
 \mathcal C^{\text{ss}}\ar[d,"\phi"]\ar[r,"\psi"] & \mathcal T\ar[d,"\phi_T"]\\
 \mathcal C\ar[r,"\bar{\psi}"] & \PP
 \ecd
 The line bundles $\phi_T^*\OO_{\PP}(1)$ and $\omega_{\mathcal T}(\frac{1}{2}D_B)$ have the same total degree on $\mathcal T$, so, since the latter is a tree, their difference is associated to a piecewise linear function $\lambda_T$ on $\tropT$. Pulling back via $\mathcal C$, on the other hand, we find:\[\phi^*\omega_{\mathcal C}=\omega_{\mathcal C^{\text{ss}}}(-\on{trop}(\psi)^*\lambda_T).\]
\end{proof}

For the reader's benefit, Figure \ref{fig:admredcores} provides a pictorial description of such piecewise linear functions in case the core is a genus two configuration of rational curves, including the \emph{slopes} along the nodes (or rational chains). This is where the \emph{admissible functions} of Section \ref{sec:admfun} originate from. The two columns correspond to the two maximally degenerate (stable) tropical curves of genus two. Here, the red vertex(ices) correspond(s) to the component(s) mapping to the special component(s) of the singularity/ribbon. Blue legs represent branch points of the admissible cover. Tails are not drawn, but can be understood by a tropical modification: the corresponding vertices  would lie on trees of slope $1$ towards the core, all at the same height, and lower than the red vertex(ices). See also the conventions set out in Section \ref{sec:admfun}.

\begin{rem}
 The non-reduced structures of genus two are numerous; indeed, as opposed to the case of isolated singularities, where the $\delta$-invariant is always related to the number of branches via the genus, in the non-reduced case a very large $\delta$-invariant can always be compensated by lowering the genus of the ribbon (that can be negative). 
\end{rem}
Yet, if we require $\omega_C\geq0$, we find a characterisation of tailed ribbons.
\begin{prop}
 Let $C$ be a Gorenstein curve of genus two, consisting of a ribbon glued \emph{in some way} to a number of (possibly singular) rational curves. Assume that the dualising sheaf is non-negative, and positive on the ribbon. Then $C$ is a $(k_1,\ldots,k_r)$-tailed ribbon.
\end{prop}
\begin{proof}
 We may assume that the normalisation sequence is
 \[ 0 \to \OO_C\to \OO_R\oplus\bigoplus_{i=1}^r\OO_{\PP^1}^{\oplus k_i}\to\bigoplus_{i=1}^r\mathcal F_i\to 0,\]
 where the $\mathcal F_i$ are sheaves supported at the gluing points; call $\delta_i$ their lengths. We claim that $\delta_i\geq k_i$; indeed, this is the minimum number of conditions that we need to impose to ensure that the values of the functions on different branches agree at the preimage of the singular point. But in fact we may even assume $\delta_i\geq k_i+1$, by the Gorenstein assumption and Lemma \ref{lem:indecomposable}.
 
 Now, it is easy to compute that for the genus of $C$ to be two we need
 \[\mathcal I_{R_\text{red}/R}\simeq \OO_{\PP^1}\left(\sum_{i=1}^r(\delta_i-k_i)-3\right).\]
 So, by applying Noether's formula and adjunction as in Remark \ref{rem:Noetheradjunction}, we find
 \[\omega_{C|_{R_{\text{red}}}}\simeq\omega_{\PP^1}\otimes\mathcal I_{R_\text{red}/R}^\vee\otimes(\mathfrak c_{|R_{\text{red}}})^\vee\simeq \OO_{\PP^1}\left(1+\sum_{i=1}^r(k_i-\deg(\mathfrak c_{i|R_{\text{red}}})-\delta_i)\right).\]
 We study the local contributions $k_i-\deg(\mathfrak c_{i|R_{\text{red}}})-\delta_i$ around each singular point. We are then looking for Gorenstein subalgebras $A$ of $\widetilde A=\mathbb C[\![s]\!][\epsilon](\epsilon^2)\times \mathbb C[\![t_1]\!]\times\ldots\times \mathbb C[\![t_k]\!]$. By Corollary \ref{cor:principalconductor} we know that $\mathfrak c$ is principal as an ideal of $\widetilde A$; up to rescaling the generators, we may assume that $\mathfrak c$ is generated by the element
 \[(s^c+\epsilon s^d,t_1^{c_1},\ldots,t_k^{c_k}).\]
 Indeed, under our assumptions, the degree of $\omega_C$ has to be one on $R_{\text{red}}$, and zero on every tail; thus, it is easy to see that $c_1=\ldots=c_k=2$.
 Now, if $d\geq c$, then $s^c\in \mathfrak c$; if $d<c$, then $s^{2c-d}\in\mathfrak c$. In the second case,
 \begin{align*}
  \widetilde A/\mathfrak c=&\langle 1, s,\epsilon,\ldots,\epsilon s^d(=-s^c),s^{d+1},\ldots, s^{c-1},\epsilon s^{c-1}(=s^{2c-d-1});\\
  &1,t_1;\\
  &\ldots\\
  &1,t_k\rangle
 \end{align*}
In any case we find $2\delta\overbrace{=}^{\text{Lemma \ref{lem:n=2delta}}}\dim_\mathbb{C}\widetilde A/\mathfrak c=2c+2k$, and $\delta_i=k_i-\deg(\mathfrak c_{i|R_{\text{red}}})$ for all $i$.

By assumption (that the reduced subcurve underlying $R$ is smooth), we must have a generator $y$ of $A$ whose linear part contains $s+\epsilon p(s)$. By Lemma \ref{lem:indecomposable} and $t_i^2\in\mathfrak c$, there must be generators with non-trivial linear part in $t_i$. Up to taking polynomial combinations of these generators, we may assume they take the form:
\begin{align*}
 y=&\ s+\epsilon p(s)&\\
 x_i=&\ \epsilon q_i(s)\oplus t_i, &\quad i=1,\ldots,k
\end{align*}
where $p$ and $q_i$ are monomial. Since $1+\epsilon \tilde{p}(s)$ is invertible in $\widetilde A$, we may assume that $y=s$. Furthermore, by looking at the conductor ideal, we see that for $x_i$ we may assume $q_i(s)=q_is^{c-1}$, with $q_i\in\mathbb C^\times$. Up to blowing $\widetilde A$ down (see Remark \ref{rem:normalisation}), we may reduce to the case $c=1$. Thus, we find the singularities of Example \ref{exa:tailedribbons}, with local equations:
\[ \mathbb{C}[\![x_1,\ldots,x_k,y]\!]/(x_ix_j,x_iy-x_jy)_{1\leq i<j\leq k},\]
and ideal sheaf $\mathcal I_{R_\text{red}/R}=\OO_{R_\text{red}}(r-3)$.
\end{proof}

\subsection{$h^1-$vanishing}
First we describe all the possible cores. Recall the following:
\begin{dfn}\label{def:minimal}
 A curve is \emph{minimal} if there is no proper subcurve of the same arithmetic genus.
\end{dfn}

\begin{rem}
 Up to replacing nodes with rational bridges, minimal curves of genus one are either irreducible or elliptic $m$-fold points with rational branches ($m\geq 2$).
\end{rem}

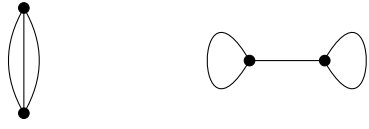
\begin{figure}
 \begin{tikzpicture}
  \draw[fill] (0,.7) circle (2pt) (0,-.7) circle (2pt) (3,0) circle (2pt) (4,0) circle (2pt);
  \path (0,.7) edge (0,-.7) (0,.7) edge [bend left] (0,-.7) (0,.7) edge [bend right] (0,-.7);
  \draw[scale=3] (1,0) node{} to [out=120,in=240,loop] ();
  \draw[scale=3] (1.33,0) node{} to [out=60,in=-60,loop] ();
  \draw (3,0)--(4,0);
 \end{tikzpicture}
\caption{Minimal theta and dumbbell configurations (dual graphs).}
\label{fig:thetadumbbell}
\end{figure}

\begin{prop}
 Up to replacing nodes with rational bridges, a \emph{minimal} Gorenstein curve of genus two with isolated singularities is one of the following:
 \begin{enumerate}
  \item irreducible;
  \item the nodal union of two minimal curves of genus one (including a dumbbell configuration of $\PP^1$s, see Figure \ref{fig:thetadumbbell});
  \item a theta configuration of $\PP^1$s (see Figure \ref{fig:thetadumbbell});
  \item a genus one singularity with a minimal genus one branch;
  \item a type I singularity with rational branches;
  \item a type II singularity with rational branches.
 \end{enumerate}
\end{prop}

\begin{lem}
 With numbering as in the previous Proposition, the dualising line bundle of a minimal Gorenstein curve of genus two has:
 \begin{enumerate}
 \setcounter{enumi}{3}
  \item degree $2$ on the genus one branch, or on the common branch of the two elliptic $m$-fold points;
  \item degree $2$ on the special branch;
  \item multidegree $(1,1)$ on the special branches.
 \end{enumerate}
\end{lem}
\begin{proof}
 Follows from Noether's formula (Proposition \ref{prop:Noether}) and an explicit calculation.
\end{proof}

Now we describe sufficient conditions for the vanishing of higher cohomology of a line bundle on a minimal curve of genus two. This will allow us to prove the unobstructedness of the $\lambda-$aligned admissible maps satisfying the factorisation property. The proof of the following lemmas boils down to a simple albeit tedious application of the normalisation exact sequence.

\begin{lem}\label{lem:h1van1}
 Let $C$ be a minimal Gorenstein curve of genus two with isolated singularities. A line bundle $L$ on $C$ having non-negative multidegree, positive degree on every subcurve of genus one, degree at least two, and $L\neq\omega_C$, has vanishing $h^1$.
\end{lem}

\begin{lem}\label{lem:h1van2}
 Let $C$ be a $(k_1,\ldots,k_r)$-tailed ribbon, and $L$ a non-negative line bundle on it. Then $h^1(C,L)=0$ if:
 \begin{itemize}
  \item $L$ has positive degree on at least two $k$-fold points; or
  \item $L$ restricts to $\OO_{\PP^1}(1)$ on $R_\text{red}$ (where $R$ denotes the component of multiplicity two), and it has positive degree on at least one $k$-fold point.
 \end{itemize}
\end{lem}

\section{Gorenstein curves over aligned admissible covers}\label{sec:alignedadmissiblecovers}

\begin{dfn}
 A weighted admissible cover consists of a hyperelliptic admissible cover \[\psi\colon(C,D_R)\to(T,D_B)\] with prestable target, together with a weight function $w\colon V(\tropC)\to\NN$, such that, with the induced weight function $w_T\colon V(\tropT)\to\NN$ given by
 \[w_T(v)=\sum_{v^\prime\in\psi^{-1}(v)}w(v^\prime),\]
 the rational tree $T$ is \emph{weighted-stable} (every weight-zero component has at least three nodes or branch points).
 We denote by $\mathcal A_2^{\text{wt}}$ the smooth Artin stack of weighted admissible covers where $C$ has arithmetic genus two.
\end{dfn}

\subsection{The intuition and strategy} Our goal is to produce a morphism to a Gorenstein curve $\cC\to\OC$ such that a line bundle on $\OC$ of degree as prescribed by the weight function $w$ would have vanishing higher cohomology. This would ensure the unobstructedness of the space of maps to projective space.

Classically, we would look for a vertical divisor $\pazocal Z$, supported on the exceptional locus of $\cC\to\OC$, such that the contraction is associated to a line bundle of the form $\omega_{\cC}(\pazocal Z)$ - this guarantees that $\OC$ is Gorenstein. Tropically, $\pazocal Z$ is replaced by a PL function $\lambda$ on $\tropC$. 
The relevant information encoded by $\lambda$ is a collection of slopes along the edges of $\tropC$. These are the objects of the stack \textbf{Div}, introduced in \cite[\S 4.2]{MarcusWise}.

Finding $\lambda$ such that $\omega_{\cC}(\lambda)$ is trivial on the exceptional locus of $\cC\to\OC$ reduces to a simple calculation of degrees - as opposed to a more complicated equality in the Picard group of a genus two curve - because we impose that $\lambda$ is pulled back from the target of the admissible cover.

When the weight of the core is zero, basically any of the positive-weight tails can be elected as the special branch of the Gorenstein singularity $\overline{C}$, be it isolated or a ribbon. The choice of $\lambda$ can be thought of as the choice of a degree-one divisor (i.e. a point) on the target of the admissible cover. The support of this divisor corresponds to the special component of the singularity.

Given a standard tropical weighted admissible cover over $\mathbb R_{\geq 0}$, the function $\lambda$ is uniquely determined. It is a member of the tropical canonical linear series on the maximal subcurve $\tropD\subseteq\tropC$ (making it into a \emph{level graph} in the sense of \cite{BCGGM}) such that:
\begin{itemize}[leftmargin=.5cm]
 \item the interior $\tropD^\circ$ contains every subcurve of arithmetic genus $g$ and $w\leq2g-2$;
 \item every connected component $\Omega$ of $\tropD^\circ$ has weight at most $2p_a(\Omega)-2$, and $\overline{\Omega}$ has weight at least $2p_a(\Omega)-1$.
\end{itemize}
For a more general family of tropical curves, this determines a (polyhedral and simplicial) subdivision of the base cone. In order to prove this, we actually define $\lambda$ by interpolating among finitely many piecewise linear functions, which we call \emph{admissible}. 
Our definition of admissible functions has been inspired by that of S. Bozlee's \emph{mesa curves} \cite[\S 3]{Bozlee}. It could be said that it constitutes a Copernican revolution with respect to \cite{RSPW1}.

In case the special component has weight $2$, if this represents the degree of a map, the only possibility is that the map restrict on the core to the hyperelliptic cover of a line. In this case, factorisation through the Gorenstein singularity is not enough to ensure smoothability. If we again replace the core by a ribbon, imposing a second factorisation requires the map to satisfy some ramification condition on the nearby branches, which turns out to be sufficient for the obstructions to vanish.

In order to avoid complications (and nasty singularities) arising from non-factoring situations - i.e. when the core has weight $2$, but it cannot possibly be the weight of the dualising bundle of a Gorenstein singularity - we proceed in two steps: first, in \S \ref{sec:admfun}, we make sure that every subcurve of positive genus has positive weight. Then, in \S \ref{sec:secalign}, we discard the locus where the weight of the core is $1$ or $1+1$. After doing this, we can safely replace any weight-two core with a ribbon. We sum up the construction of this and the next sections in the following diagram of moduli spaces:

\bcd[cramped]
\VZ(X)\ar[r,hook] &[-.7cm] \doublewidetilde{\mathcal A_2}(X)\ar[d]\ar[r]\ar[dr,phantom,"\Box"] & \TAt(X)^{\text{fact}}\ar[d]\ar[r,hook] & \TAt(X)\ar[d]\ar[r]\ar[dr,phantom,"\Box"] & \oM_2(X)\ar[d] \\
& \doublewidetilde{\mathcal A_2}\ar[r] & \TAt^\circ \ar[r,hook] & \TAt \ar[r] &\mathfrak M_2^{\text{wt}} \\
& \mathcal{C}\leftarrow\TTC\to\Cpp\to\OOC\ar[u] & &\mathcal{C}\leftarrow\TC\to\Cp\to\OC\ar[u] & \mathcal{C}\ar[u]
\ecd

\subsection{Admissible functions and aligned admissible covers}\label{sec:admfun} 
Let $\psi\colon(C,D_R)\to(T,D_B)$ be a weighted admissible cover over $S$, endowed with the minimal logarithmic structure. Let $\tropT$ denote the tropicalization of the target. $\tropT$ is therefore a weighted \emph{tree}, metrised in the monoid $\overline M_S$. The metric structure can be disregarded for the minute.

It follows from the Riemann-Hurwitz formula that the following $\frac{1}{2}\mathbb{Z}$-divisor $D^\prime$, supported on the vertices of $\tropT$, pulls back to the canonical divisor of $\tropC$ (see \S \ref{sec:PLfun}):
 \begin{equation}\label{eqn:Dv}
 D^\prime(v)=\val(v)-2+\frac{1}{2}\on{deg}(D_B|_v),
 \end{equation}
 where $\val(v)$ denotes the valence (number of adjacent edges) of the vertex $v$, and $D_B$ is the branch divisor of the cover. The degree of $D^\prime$ on $\tropT$ is $1$.

 Let $D$ be another $\frac{1}{2}\mathbb{Z}$-divisor of degree $1$ on $\tropT$. Since $\tropT$ is a tree, its Jacobian is trivial, so any two divisors of the same degree are linearly equivalent (this carries over to $\frac{1}{2}\mathbb{Z}$-divisors). Therefore, there exists a unique collection $\overline{\lambda}_T$ of half-integral slopes on the edges of $\tropT$
  such that
 \[D=D^\prime+\on{div}(\overline\lambda_T).\] If $D$ has integral coefficients, we observe that half-integral slopes may only appear along edges over which $\on{trop}(\psi)$ has expansion factor two, so it is possible to think of $\overline{\lambda}_T$ as a piecewise-linear function on $\tropT$ with values in $\overline M_S^\text{gp}$, up to a global translation of $\overline M_S^\text{gp}$. In fact, it may be more accurate to think of $\overline\lambda_T$ as a PL function on the tropicalization of the orbicurve $[C/\mathfrak S_2]$.

 If $D$ is effective, the pullback $\overline\lambda$ of $\overline\lambda_T$ along $\on{trop}(\psi)$ - or rather the divisor $K_{\tropC}+\on{div}(\overline\lambda)$ - is an element of the tropical canonical linear system.
 


\begin{dfn}\label{def:admissiblefunction}
 Let $S$ be a geometric point.
 An \emph{admissible function} $\overline\lambda$ on $\tropC$ is a piecewise-linear function with integral slopes (with values in $\overline M_S^\text{gp}$, but defined only up to a global translation by $\overline M_S^\text{gp}$) such that:
 \begin{itemize}[leftmargin=.5cm]
  \item $\overline\lambda$ is the pullback along $\on{trop}(\psi)$ of a (possibly half-integral) $\overline\lambda_T$ on $\tropT$;
  \item $K_{\tropC}+\operatorname{div}(\overline\lambda)\geq 0$ defines an element of the tropical canonical system, which is the pullback of an effective, degree 1 divisor $D$ supported on \emph{exactly one} vertex of $\tropT$.
 \end{itemize}
\end{dfn}

We are going to call \emph{special} the vertex of $\tropT$ supporting $D$, or its preimage(s) in $\tropC$. With an eye to the future, it will correspond to the special branch(es) of the Gorenstein curve of genus two.

\begin{rem}\label{rmk:finiteness}
 There are only \emph{finitely many} admissible functions compatible with a given $\psi$. Indeed, they are in bijection with the vertices of $\tropT$.
\end{rem}

When the core is maximally degenerate, i.e. a configuration of rational curves, all the possible admissible functions are depicted in Figure \ref{fig:admredcores}. The red vertex(ices) represent(s) the one(s) supporting $K_{\tropC}+\operatorname{div}(\overline\lambda)$. Missing from the picture is what happens outside the core, but this is easily explained: there may be any number of rational trees, on which $\overline\lambda$ has constant slope $1$ towards the core.

\smallskip

For the sake of concreteness, we now look at some examples, adopting the:

\begin{Convention}
With the application to stable maps in mind, we think of the source of an admissible cover as the destabilisation of a weighted-stable curve. In the following, we represent a component of positive weight by a black circle, and a component of weight zero by a white one, unless it is unstable (either a rational tail introduced as the conjugate of a rational tail of positive weight, or a rational bridge introduced by slicing $\tropC$), in which case it is represented by a cross, and the corresponding edge is dotted. A red vertex represents the component supporting the divisor $D$, or its preimages in $\tropC$. The blue legs (\emph{B-legs}) represent the branching divisor $D_B$ of $\psi$ (see Definition \ref{def:admissiblecover}): the number of B-legs attached to a vertex $v$ of $\tropT$ is the degree of $D_B\cap T_v$, where $T_v$ is the irreducible component of $T$ corresponding to $v$.
\end{Convention}

\begin{exa}\label{ese:lambdabar}
Assume the core is irreducible of weight zero, and there are two tails of (large) positive weight: one of them ($R$) is attached to a Weierstrass point of the core; the other one ($L$) is attached to a general point, so that a weight zero tail ($\bar L$) must be sprouted from the conjugate point in order for the admissible cover to exist (Figure \ref{fig:exa1}). The tropical cover has expansion factor $2$ (in gray) along the rightmost edge, corresponding to the ramification order of the algebro-geometric map at each branch of the node.
\begin{figure}[hbt]
\begin{tikzpicture} [scale=.4]
\tikzstyle{every node}=[font=\normalsize]
		\tikzset{arrow/.style={latex-latex}}
		
	\tikzset{cross/.style={cross out, draw, thick,
         minimum size=2*(#1-\pgflinewidth), 
         inner sep=1.2pt, outer sep=1.2pt}}


		\coordinate (Oc) at (0,3,0);
		\coordinate (T) at (-2,4,0);
		\coordinate (bT) at (-2,2,0);
		\coordinate (W) at (3,1.5,0);
		\coordinate (Ot) at (0,0,0);
		\coordinate (Tt) at (-2,0,0);
		\coordinate (Wt) at (3,0,0);
		
		\coordinate (P1) at (-.5,-1,0);
		\coordinate (P2) at (0,-1,0);
		\coordinate (P3) at (.5,-1,0);
		\coordinate (P4) at (1,-1,0);
		\coordinate (P5) at (1.5,-1,0);
    	\coordinate (P6) at (3.2,-.8,0);

		\begin{scope}[every coordinate/.style={shift={(0,1,0)}}]
		
		
	\draw  ([c]T) node [cross]{};
			
	\foreach \x in {bT,W,Tt,Wt} 
          \fill ([c]\x) circle (3pt);
     	
     	\draw ([c]Oc)node[right=.3cm]{\small $Z$}-- ([c]bT) node[left]{\small$L$};
     	\draw[dotted] ([c]Oc)-- ([c]T) node [left]{\small$\bar L$};
     	\draw ([c]Oc)--node[below,gray]{$2$} ([c]W) node [right]{\small$R$};
     	\draw ([c]Ot)-- ([c]Tt);
     	\draw ([c]Ot)-- ([c]Wt);
     	
     	\draw[blue] ([c]Ot)--([c]P1);
     	\draw[blue] ([c]Ot)--([c]P2);
     	\draw[blue] ([c]Ot)--([c]P3);
     	\draw[blue] ([c]Ot)--([c]P4);
     	\draw[blue] ([c]Ot)--([c]P5);
     	\draw[blue] ([c]Wt)--([c]P6);

       \draw[dashed] ([c]Oc)--([c]Ot);
       \draw[dashed] ([c]T)--([c]Tt);
       \draw[dashed] ([c]W)--([c]Wt);
       
       \draw[fill=white]  ([c]Ot) circle (3pt);
       \draw[fill=white]  ([c]Oc) node[above]{\tiny $g=2$} circle (3pt);
         \end{scope}
       \end{tikzpicture}
\caption{A weighted admissible cover with weight-zero core.}
\label{fig:exa1}
\end{figure}
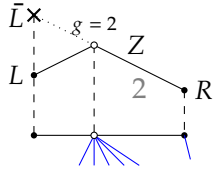

There are three admissible functions $\overline{\lambda}_i$ compatible with the given weighted admissible cover; see Figure \ref{fig:exa1adm}.
 
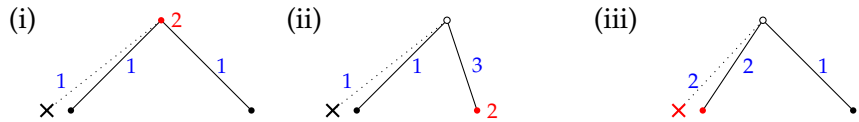
\begin{figure}[h]
\begin{tikzpicture} [scale=.4]
\tikzstyle{every node}=[font=\normalsize]
		\tikzset{arrow/.style={latex-latex}}
		
	\tikzset{cross/.style={cross out, draw, thick,
         minimum size=2*(#1-\pgflinewidth), 
         inner sep=1.2pt, outer sep=1.2pt}}
         
\coordinate (lab) at (-3.8,0,0);
\coordinate (O) at (0,0,0);
\coordinate (T1) at (-3,-3,0);
\coordinate (T1B) at (-3.8,-3,0);	
\coordinate (W) at (3,-3,0);
  \coordinate (W1) at (1,-3,0);
  
  \coordinate (T2) at (-2,-3,0);
\coordinate (T2B) at (-2.8,-3,0);


\coordinate (s1) at (-1.5,-1.5,0);
\coordinate (s1B) at (-2.8,-2,0);	
\coordinate (sW) at (1.5,-1.5,0);	  
\coordinate (sW1) at (.5,-1.5,0);	 		  
	\coordinate (s2) at (-1,-1.5,0);
\coordinate (s2B) at (-1.8,-2,0);


		  \begin{scope}[every coordinate/.style={shift={(-8,0,0)}}]
		  
	  \draw ([c]lab) node [left]{(i)};
          
          \draw ([c]O)-- ([c]T1);
          \draw [dotted] ([c]O)-- ([c]T1B);
               \draw ([c]O)-- ([c]W);
        
            \node at ([c]s1) [right] {\begin{color}{blue}\scriptsize{$1$}\end{color}};   
          \node at ([c]s1B) [left] {\begin{color}{blue}\scriptsize{$1$}\end{color}}; 
           \node at ([c]sW) [right] {\begin{color}{blue}\scriptsize{$1$}\end{color}}; 
       
          \fill[red]  ([c]O) circle (3pt) node[right]{\color{red}\scriptsize{$2$}};
	  \draw ([c] T1B)   node [cross] {};
	  \fill ([c]W) circle (3pt);
          \fill ([c]T1) circle (3pt);

         \end{scope}
		
		  \begin{scope}[every coordinate/.style={shift={(1.5,0,0)}}]
	
        \draw ([c]lab)node[left]{(ii)};
          
          \draw ([c]O)-- ([c]T1);
          \draw [dotted] ([c]O)-- ([c]T1B);
          \draw ([c]O)-- ([c]W1);
        
          \node at ([c]s1) [right] {\begin{color}{blue}\scriptsize{$1$}\end{color}};   
          \node at ([c]s1B) [left] {\begin{color}{blue}\scriptsize{$1$}\end{color}}; 
          \node at ([c]sW1) [right] {\begin{color}{blue}\scriptsize{$3$}\end{color}};

	   \draw[fill=white]  ([c]O) circle (3pt);
	   \draw ([c] T1B)   node [cross] {};
	   \fill[red] ([c]W1) circle (3pt) node[right]{\color{red}\scriptsize{$2$}};
           \fill ([c]T1) circle (3pt);
       
         \end{scope}
		
 \begin{scope}[every coordinate/.style={shift={(12,0,0)}}]
          
          \draw ([c]lab)node[left]{(iii)};
          
          \draw ([c]O)-- ([c]T2);
          \draw[dotted]  ([c]O)-- ([c]T2B);
          \draw ([c]O)-- ([c]W);
        
          \node at ([c]s2) [right] {\begin{color}{blue}\scriptsize{$2$}\end{color}};   
          \node at ([c]s2B) [left] {\begin{color}{blue}\scriptsize{$2$}\end{color}}; 
          \node at ([c]sW) [right] {\begin{color}{blue}\scriptsize{$1$}\end{color}};

	    \draw[fill=white] ([c]O) circle (3pt);
	    \draw ([c] T2B)  node [cross, red] {};
  	    \fill ([c]W) circle (3pt);
            \fill[red] ([c]T2) circle (3pt);

         \end{scope}

	\end{tikzpicture}
	\caption{Admissible functions compatible with Figure \ref{fig:exa1}.}
	\label{fig:exa1adm}
\end{figure}

Note that the pullback of $D$ always has degree $2$ - although the tropical cover is injective on the rightmost edge, the expansion factor $2$ provides the correct multiplicity. It may be instructive to compute the multi-degree of various tropical divisors on $\tropC$ - we represent them as vectors in $\ZZ^4$ by ordering the components from left to right $(\bar L,L,Z,R)$. The canonical $K_{\tropC}$ gives $(-1,-1,5,-1)$.

\begin{centering} 
 \begin{tabular}{c||c|c}
  & $\on{div}(\overline{\lambda}_i)$ &  $\on{trop}(\psi)^*D=K_{\tropC}+\on{div}(\overline{\lambda}_i)$ \\
  \hline
  \hline
  (i) & $(1,1,-3,1)$ & $(0,0,2,0)$ \\
  (ii) & $(1,1,-5,3)$ & $(0,0,0,2)$ \\
  (iii) & $(2,2,-3,1)$ & $(1,1,0,0)$
 \end{tabular}

\end{centering}
\end{exa}


\begin{dfn}\label{dfn:alignmentsubdivision}
 A \emph{pre-aligned} admissible cover over a logarithmic scheme $(S,M_S)$ is one for which the values
 \begin{equation}\label{eqn:alignment}\{\overline\lambda(v)-\overline\lambda(v^\prime)|v,v^\prime\in V(\tropC)\}\subseteq (\overline M_S,\geq)\end{equation}
 are comparable, for every admissible function $\overline\lambda$ compatible with $\psi$ (Definition \ref{def:admissiblefunction}).
\end{dfn}

Pre-aligned admissible covers form a subfunctor of weighted admissible covers over $(\text{LogSch})$. The next lemma follows from the definitions and Section \ref{sec:minimality}.

\begin{lem}
 The minimal logarithmic structure of a pre-aligned admissible cover is obtained, starting from the minimal logarithmic structure of the admissible cover, by adding in the elements of $\overline M_S^{\text{gp}}$ indicated in \eqref{eqn:alignment} and sharpifying.
\end{lem}

\begin{rem}
 Aligning determines a subdivision $\Sigma$ of the tropical moduli space $\sigma=\Hom(\overline M_S,\mathbb R_{\geq0})$, or, equivalently, a logarithmic modification of $\mathcal A_{2}^\text{wt}$, which will be denoted by $\mathcal{A}^{\text{pre}}_{2}$. See Section \ref{sec:logblowups}.
\end{rem}

 Let $\psi$ be a pre-aligned admissible cover, and $\{\overline\lambda_1,\ldots,\overline\lambda_s\}$ the set of admissible functions compatible with $\psi$. We are going to choose genuine PL functions $\lambda_i$ on $\tropC$ with values in $\overline M_S^\text{gp}$ lifting the $\overline\lambda_i$. A lift determines and is determined by the set of vertices mapping to $0$ in $\overline M_S^{\text{gp}}$; let us denote by $\tropC_{\geq0}$ (resp. $\tropC_{>0}$) the set of vertices with values in $\overline M_S\subseteq \overline M_S^{\text{gp}}$ (resp. $\overline M_S\setminus \overline M^*_S$; note that we may assume $\overline M^*_S=\{0\}$). In the end, we are going to define the object of interest as an interpolation/truncation of these lifts.
 
\begin{dfn}\label{def:lambdalift}
 Let $S$ be a geometric point. Define a lift $\lambda_i\in\Gamma(S,\pi_*\overline M_C^{\text{gp}})$ of $\overline\lambda_i$ 
 by requiring that $w(\tropC_{>0})\leq 0$ and $w(\tropC_{\geq 0})\geq 1$.
\end{dfn}

\begin{dfn}
 Define a PL function on a subdivision $\widetilde{\tropC}$ of $\tropC$ by:
 \[\lambda=\max\{0,\lambda_1,\ldots,\lambda_s\}.\]
 
We denote by $\tropD^\circ$ the support of $\lambda$ (appeared as $\tropC_{>0}$ before), and by $\tropD$ the minimal subcurve of arithmetic genus two in $\widetilde{\tropC}$ containing the closure of $\tropD^\circ$.
\end{dfn}

\begin{rem}
 Parallel to the above definition, we denote by $\tropD_T^\circ$ the support of $\lambda_T$, and by $\tropD_T$ the image of $\tropD$ under $\trop(\psi)$.
 
 Since tropical linear systems are tropically convex (see \S \ref{sec:PLfun}), it is still true that $\lambda$ is the pullback of a PL function $\lambda_T$ (with slopes in $\frac{1}{2}\mathbb{Z}$) on $\tropT$, such that $D^\prime+\on{div}(\lambda_T)$ is an effective divisor of degree $1$ on $\tropD_T$ (notice that here, in the definition of $D'$, the valence of a vertex is the one in $\tropD_T$ and not the one in $\tropT$). This is clear away from $0$, where $\lambda$ coincides with one of the admissible functions above; at $0$, on the other hand, all the slopes of $\lambda_T$ are nonnegative, so the only doubt is for a vertex of valence $1$ supporting exactly one $B$-leg; but such a vertex corresponds to a rational tail in $\tropD$, so it can be included in $\tropD_T$ only if the slope of $\lambda_T$ along the uniqe adjacent edge is strictly positive.
 
 The divisor $D=D^\prime+\on{div}(\lambda_T)$  needs not  be supported on a vertex of $\tropT$, and could be the sum of two half-points in the locus over which $\on{trop}(\psi)$ has expansion factor $2$.
 
 Since by Riemann-Hurwitz we have $\trop(\psi)^*(D')=K_{\tropD}$, the following holds:
 \begin{equation}\label{eq:fundamental}
  \trop(\psi)^*(D)=K_{\tropD}+\on{div}(\lambda).
 \end{equation}

\end{rem}

\begin{lem}
 On every cone of $\on{trop}(\mathcal{A}^{\text{pre}}_{2})$, $\lambda$ is well-defined as a PL function on a subdivision of $\tropC$.
\end{lem}

\begin{proof}

For $\lambda$ to be well-defined, we need to argue that, for every vertex $v$ of $\tropC$, the values $\{\lambda_i(v)|i=1,\ldots,s\}$ are comparable, and there is a unique way of subdividing $\tropC$ in order to make $\lambda$ piecewise-linear. Notice that, a priori, we only know that the values of a single $\lambda_i$ at the vertices of $\tropC$ are ordered.

It is convenient to work on $\tropT$. The leaves of $\tropT$ have positive weight due to stability, so $\lambda_T$ necessarily takes the value $0$ on them.

In order to prove that the maximum is well defined at every vertex, we can proceed inductively from the leaves and run through all of $\tropT$. Assume that $\lambda_T(v_1)$ has been determined. We need to establish the behaviour of $\lambda_T$ on the edge $e$ between $v_1$ and $v_2$, and its value on $v_2.$ Upon relabelling, we may suppose that $\lambda_1$ is the function with maximal slope along $e$ among the ones with $\lambda_i(v_1)=\lambda_T(v_1)$; so $\lambda_T$ coincides with $\lambda_1$ in a neighbourhood of $v_1$. 
If  $\lambda_1$ is also the function of maximal slope along $e$ among all the admissible functions, then $\lambda_T(v_2)=\lambda_1(v_2)$ and we are done.

If not, there must be a function $\lambda_2$ with slope along $e$ greater than that of $\lambda_1$ but $\lambda_2(v_1)\leq\lambda_1(v_1).$
First, assume that $\overline\lambda_1$ and $\overline\lambda_2$ differ only along $e$, i.e. $D(\lambda_1,v_2)=1$ and $D(\lambda_2,v_1)=1$. Then, there is a unique way to interpolate between $\lambda_1$ and $\lambda_2$ along $e$, namely by solving the equations:
 \[
 \left\{
 \begin{array}{ r @{{}={}} l }
  \ell(e)&\ell(e^\prime)+\ell(e^{\prime\prime})\\
  \lambda_2(v_2)-\lambda_1(v_1) & s(\lambda_1,e)\ell(e^\prime)+s(\lambda_2,e)\ell(e^{\prime\prime}).
 \end{array}
 \right.
\]
This uniquely determines the subdivision of $e$ on which $\lambda_T$ is piecewise-linear, and the values of $\lambda_T$ on both $v_2$ and the newly introduced vertex $v_{1,2}.$

More generally, $\lambda_1$ and $\lambda_2$ will differ along the chain $R$ connecting the support of $D_1$ to that of $D_2$, and containing the edge $e$.
 Then $\lambda_2-\lambda_1=D_2-D_1$ is a rational function with constant slope $1$ along $R$.
 In order to define $\lambda_T$ we need to interpolate between $\lambda_1$ and $\lambda_2$ along $R$, see Figure~\ref{fig:lambdacomparison}.

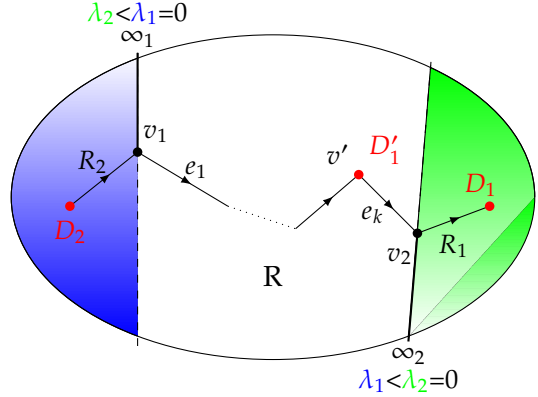
\begin{figure}[htb]
	\centering
		\begin{tikzpicture} [scale=.6]
		\tikzstyle{every node}=[font=\normalsize]
		\tikzset{arrow/.style={latex-latex}}

\draw[name path=E] ellipse (5.8 and 3.6);
\draw[dashed,name path=L1] (3.5,3) to (3,-3.2);
\draw[dashed,name path=L2] (-3,3.2) to (-3,-3.3);

\begin{scope}[on background layer]
\filldraw[draw=black,bottom color=blue, top color=white]
    [intersection segments={of= L2 and E}];
    
    \filldraw[bottom color=white, top color=green]
   [intersection segments={of=E and L1}];
    \path[name path=line1] (3,-3.1) --  (5.7,0);
    \fill[bottom color=white, top color=green]
[intersection segments={of=line1 and E}];
 \end{scope}
 \coordinate (A) at (-4.5,-0.2);
  \coordinate (B) at (-3,1);
  \coordinate (C) at (-1,-.2);

    \coordinate (C1) at (.5,-.7);
    \coordinate (D) at (1.9,0.5);
    \coordinate(E) at (3.2,-0.8);
    \coordinate (F) at (4.8,-.2);

\draw[->-] (A) to (B);
\draw[->-] (B) to (C);
\draw[dotted] (C) to (C1);

\draw[->-] (C1) to (D);
\draw[->-] (D) to (E);
\draw[->-] (E) to (F);
\draw[thick] (B) to (-3,3.2);
\draw[thick] (E) to (3,-3.2);

\node at (-3,3.6) [above]  {\small{\begin{color}{green}$\lambda_2$\end{color}<\begin{color}{blue}$\lambda_1$\end{color}=0}};
\node at (-3,3.1) [above]  {\small{$\infty_1$}};

\node at (-4.5,-0.3) [below] {\color{red}\small{$D_2$}};
\node at (-4,0.3) [above] {\small{$R_2$}};
\node at (-2.6,1) [above] {\small{$v_1$}};
\node at (-1.7,0.2) [above] {\small{$e_1$}};
\node at (0,-1.3) [below] {R};
\node at (1.9,.5) [above left] {\small{$v'$}};
\node at (1.9,.5) [above right] {\color{red}\small{$D'_1$}};
\node at (2.7,-.4) [left] {\small{$e_k$}};
\node at (2.8,-.9) [below] {\small{$v_2$}};
 \node at (4,-.6) [below] {\small{$R_1$}};
 
 \node at (3,-3.6) [below]  {\small{\begin{color}{blue}$\lambda_1$\end{color}<\begin{color}{green}$\lambda_2$\end{color}=0}};
\node at (3,-3.1) [below]  {\small{$\infty_2$}};
 \node at (4.6,-.2) [above] {\color{red}\small{$D_1$}};

 \foreach \x in {B,E}
	\fill (\x) circle (3pt);
\foreach \x in {A,D,F}
	\fill[red] (\x) circle (3pt);
\end{tikzpicture}
\caption{Comparison of $\lambda_1$ and $\lambda_2$.}
\label{fig:lambdacomparison}
\end{figure}


We will argue that this is possible because we can always inductively reduce to the above case, for which we have an explicit formula for the interpolation.
 
  We may assume that $R_2$ has zero length, because changing $\lambda_2$ on $R_2$ increases it without affecting $\lambda_2\leq\lambda_1$ (in the blue-shaded region), and similarly for $R_1$.
  
  Now let $\overline\lambda_1^\prime$ be the function that differs from $\overline\lambda_1$ only on the last edge $e_k$ of $R$. If $\lambda_1^\prime(\infty_1)=0$, then $\lambda_1^\prime\geq\lambda_1$, so that we can instead compare $\lambda_1^\prime$ and $\lambda_2$, which is possible by the inductive assumption. If instead $\lambda_1^\prime(\infty_1)<0$, then $\lambda_1^\prime(\infty_1^\prime)=0$ for some $\infty_1^\prime$ in the green-shaded region (since $\lambda_1^\prime<\lambda_1$ on the complement). In fact, it is easy to see that $\lambda_1^\prime(\infty_2)=0$, so $\lambda_2(v^\prime)=\lambda_1^\prime(v^\prime)<\lambda_1(v^\prime)$, and we conclude that $\lambda_T$ interpolates between $\lambda_1$ and $\lambda_2$ (equivalently $\lambda_1^\prime$) only along $e_k$ - which we have already solved.
 \end{proof}
 
 \begin{rem}
  Over a geometric point, there is no obstruction to defining the lifts $\lambda_i$ of $\overline{\lambda}_i$. On the other hand, they are stable under generisation (edge contraction), because the weight of $\tropC_{>0}$ stays the same, and the weight of $\tropC_{\geq 0}$ can only go up. So both the $\lambda_i$ and $\lambda$ define sections of $\pi_*\overline M_C$ over any base $S$.
 \end{rem}
 
 \begin{dfn}\label{def:semistablereduction}
  Let $\Sigma^\prime$ be the subdivision of $\sigma=\on{trop}(\mathcal A_{2}^\text{wt})$ determined by applying universal semistable reduction \cite[Theorem 2.4.2]{Molcho} to the morphism $\widetilde{\tropC}\to\sigma$.
  
  Let $a\colon\TAt\to\mathcal{A}^{\text{wt}}_{2}$ be the corresponding logarithmically \'etale model.
 \end{dfn}

 It follows from the construction and from \cite{Molcho,AK00} that there is a logarithmically smooth curve $\widetilde{\mathcal C}\to \TAt$, which is a partial destabilisation of $a^*\mathcal{C}$, such that $\widetilde{\tropC}$ is its tropicalisation, $\lambda\in\Gamma(\widetilde{\mathcal C},\overline M_{\widetilde{\mathcal C}})$, and all of its values at vertices of $\widetilde{\tropC}$ are comparable.
 
\begin{rem}
Notice that $\Sigma^\prime$ is neither finer, nor coarser than the subdivision $\Sigma$ previously defined in \ref{dfn:alignmentsubdivision}. It is not finer, because $\lambda$ is insensitive to the comparison of values of $\lambda_i$ below $0$. It is not coarser, because cones of $\widetilde{\tropC}$ do not map surjectively to cones of $\Sigma$ - see Example \ref{exa:reduciblecore}.

Notice also that $a$ is not simply a logarithmic blow-up in general. Indeed, due to the existence of some simplicial but not smooth cones in the subdivision, in order for the lattice map to be surjective it may be necessary to perform a Kummer extension of the base logarithmic structure \cite[\S 4]{BV}:
\[M_{\mathcal A}\subseteq M_{\widetilde{\mathcal A}}\subseteq\frac{1}{2}M_{\mathcal A}^\text{gp},\]
as prescribed by the subdivision $\widetilde{\tropC}$ and the slopes of $\lambda$, see Figure \ref{fig:halfintegers}. This operation should be thought of as a generalised root stack construction, or generalising the presentation of a simplicial affine toric variety as the coarse moduli of $[\Aaff^n/G]$, $G$ the finite abelian group encoding the difference between the two lattices at stake. This is different from the genus one case (it has to do with the appearance of slopes other than $1$ in $\lambda$), but it is nothing new with respect to \cite{Molcho}.

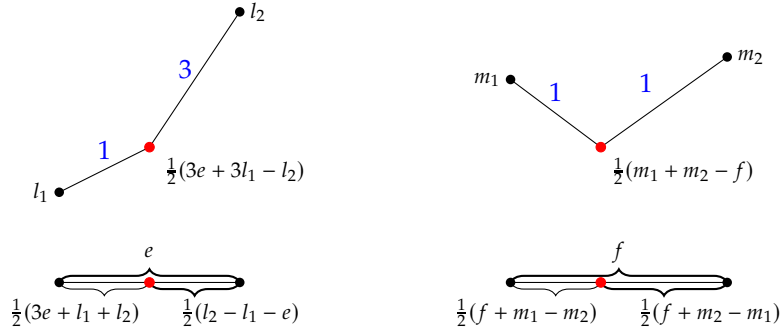
\begin{figure}[htb]
	\centering
		\begin{tikzpicture} [scale=.6]
		\tikzstyle{every node}=[font=\normalsize]
		\tikzset{arrow/.style={latex-latex}}
		
		\tikzset{cross/.style={cross out, draw, thick,
         minimum size=2*(#1-\pgflinewidth), 
         inner sep=1.5pt, outer sep=1.5pt}}
\coordinate (A) at (0,2);
\coordinate (B) at (-2,1);
\coordinate (C) at (2,5);

\coordinate (B1) at (-2,3.5);
\coordinate (C1) at (2.8,4); 

\coordinate (At) at (0,-1);
\coordinate (Bt) at (-2,-1);
\coordinate (Ct) at (2,-1);
\coordinate (C1t) at (2.8,-1); 

 \foreach \x in {B,C,Bt,Ct}
	\fill (\x) circle (3pt);

\draw (A)--(B);
\draw (A)--(C);
\draw (At)--(Bt);
\draw (At)--(Ct);
\node at (B) [left] {\scriptsize{$l_1$}};
\node at (-1,1.5) [above]  {\color{blue}\small{$1$}};
\node at (C) [right] {\scriptsize{$l_2$}};
\node at (0.8,3.3) [above]  {\color{blue}\small{$3$}};
\node at (0.1,2) [below right] {\scriptsize{$\frac{1}{2}(3e+3l_1-l_2)$}};

\draw [decorate,decoration={brace,amplitude=5pt,mirror},xshift=0.4pt,yshift=-0.4pt] ([c]Bt)--([c]At) node[black,midway,xshift=-0.4cm,yshift=-0.4cm] {\scriptsize{$\frac{1}{2}(3e+l_1+l_2)$}};
\draw [thick,decorate,decoration={brace,amplitude=5pt},xshift=0.4pt,yshift=-0.4pt] ([c]Bt)--([c]Ct) node[black,midway,yshift=0.4cm] {\scriptsize{$e$}};
\draw [thick,decorate,decoration={brace,amplitude=5pt,mirror},xshift=0.4pt,yshift=-0.4pt] ([c]At)--([c]Ct) node[black,midway,xshift=0.6cm,yshift=-0.4cm] {\scriptsize{$\frac{1}{2}(l_2-l_1-e)$}};

\draw[red,fill=red] ([c]At) circle (3pt);
\draw[red,fill=red] ([c]A) circle (3pt);

\begin{scope}[every coordinate/.style={shift={(10,0)}}]
\foreach \x in {B1,C1,Bt,C1t}
	\fill ([c]\x) circle (3pt);

\draw ([c]A)--([c]B1);
\draw ([c]A)--([c]C1);
\draw ([c]At)--([c]Bt);
\draw ([c]At)--([c]C1t);
\node at ([c]B1) [left] {\scriptsize{$m_1$}};
\node at (9,2.8) [above]  {\color{blue}\small{$1$}};
\node at ([c]C1) [right] {\scriptsize{$m_2$}};
\node at (11,3) [above]  {\color{blue}\small{$1$}};
\node at ([c]A) [below right] {\scriptsize{$\frac{1}{2}(m_1+m_2-f)$}};
\draw [decorate,decoration={brace,amplitude=5pt,mirror},xshift=0.4pt,yshift=-0.4pt] ([c]Bt)--([c]At) node[black,midway,xshift=-0.4cm,yshift=-0.4cm] {\scriptsize{$\frac{1}{2}(f+m_1-m_2)$}};
\draw [thick,decorate,decoration={brace,amplitude=5pt},xshift=0.4pt,yshift=-0.4pt] ([c]Bt)--([c]C1t) node[black,midway,yshift=0.4cm] {\scriptsize{$f$}};
\draw [thick,decorate,decoration={brace,amplitude=5pt,mirror},xshift=0.4pt,yshift=-0.4pt] ([c]At)--([c]C1t) node[black,midway,xshift=0.6cm,yshift=-0.4cm] {\scriptsize{$\frac{1}{2}(f+m_2-m_1)$}};

\draw[red,fill=red] ([c]At) circle (3pt);
\draw[red,fill=red] ([c]A) circle (3pt);

\end{scope}
\end{tikzpicture}
\caption{The two situations in which half-integral lengths occur.}
\label{fig:halfintegers}
\end{figure}

\end{rem}

\begin{thm}\label{thm:tildeA}
 The moduli space $\TAt$ of aligned admissible covers is a logarithmically smooth stack with locally free logarithmic structure (and therefore smooth) admitting a log \'{e}tale morphism to ${\mathcal A}_{2}^\text{wt}$.
\end{thm}
\begin{proof}
The construction of the minimal logarithmic structure can be traced back to Definition \ref{def:semistablereduction}. The fact that it is locally free can be justified as follows: $\lambda$ provides us with a PL map from $\widetilde{\tropC}$ to a polyhedral subdivision of $\mathbb R_{\geq 0}$. The rank of the logarithmic structure is the number of finite edges in the polyhedral subdivision of $\mathbb R_{\geq 0}$. For every level, there is at most one edge of $\widetilde{\tropC}$ mapping to it with expansion factor strictly greater than $1$: the  smoothing parameter corresponding to the said edge of $\widetilde{\tropC}$ can be taken as a generator of the minimal logarithmic structure.

 Logarithmic smoothness follows from that of $\mathcal A_{2}^\text{wt}$, since $a\colon\TAt\to \mathcal A_2$ is logarithmically \'{e}tale.
 
 Finally, a scheme logarithmically smooth over the trivial log point is smooth if and only if the stalks of its characteristic sheaf are locally free \cite[Lemma 5.2]{Niziol}.
\end{proof}

\subsection{Combinatorial properties of $\lambda$}

 In the following, we work over the standard log point, so $\tropC$ is a graph metrised in $\mathbb{R}$. To an open subcurve $\tropC_*\subseteq\tropC$ there corresponds a proper subcurve $C_*\subseteq C$; for a closed tropical subcurve, we take its algebro-geometric counterpart to be the closure of the preimage under tropicalisation. By abuse of notation, we will refer to $\tropC_*$ or to $C_*$ interchangeably.

\begin{lem}\label{lem:supp_pos_gen}
               Every component $\Omega$ of $\tropD^\circ$ has positive genus.
              \end{lem}
              \begin{proof}
               Restricting equation \eqref{eq:fundamental} to $\Omega$, we find:
               \[0\leq\on{deg}(\trop(\psi)^*(D)_{|\Omega})=2p_a(\Omega)-2+\sum_{i=1}^h (1+s(\lambda,e_i)),\]
               where $e_1,\ldots,e_h$ are the edges outgoing from $\Omega$. Since by assumption $s(\lambda,e_i)\leq-1$ for every $i$, we see that $p_a(\Omega)$ cannot be $0$.
              \end{proof}

\begin{lem}\label{lem:D1_D2}
 If $D=\frac{1}{2}(D_1+D_2)$, then $\lambda(D_1)=\lambda(D_2)=0$.
\end{lem}
\begin{proof}
 Since by definition every $\lambda_i$ has $D$ supported on a single vertex, this situation may only occur after interpolating with the zero function. 
 
 We first explain how to reduce to the path between $D_1$ and $D_2$. Consider an edge $e$ of $\tropD_T$ outside this path, and denote by $\tropD_-$ and $\tropD_+$ the connected components of $\tropD_T\setminus\{e\}$. Assume that $D$ is entirely supported on $\tropD_+$.  Restricting equation \eqref{eq:fundamental} to $\tropD_-$, we find:
 \[0=\val(\tropD_-)-2+\frac{1}{2}\# D_{B|\tropD_-}+\on{div}(\lambda_{T|\tropD_-}).\]
 We notice that $\val(\tropD_-)=1$ and that $\on{div}(\lambda_{T|\tropD_-})$ is simply the slope $s$ of $\lambda_T$ along $e$, oriented away from $\tropD_-$. Then $s=1-\frac{1}{2}\# D_{B|\tropD_-}$. Examining the restriction of \eqref{eq:fundamental} to $\tropD_+$, we notice that it is unaltered by contracting $e$ and all the edges in $\tropD_-$.
 
 Assume by contradiction that $\lambda_T(D_1)<\lambda_T(D_2)$. By the previous reduction, we may assume that $\tropD_T$ consists of the path between $D_1$ and $D_2$ only. By assumption, there must be an edge $e$ on this path along which $\lambda_T$ has non-zero slope. Contract all the rest. Once again by equation \eqref{eq:fundamental} it is easy to see that the slope can be $\frac{1}{2}$ (resp. $1$) from $D_1$ to $D_2$ if and only if the distribution of $B$-legs is $(2,4)$ (resp. $(1,5)$). Now, in order to see the half-points, the slope should be $\frac{1}{2}$; on the other hand, it has to be the slope of one of the $\lambda_i$, and we notice that they have integral slope along $e$ if the distribution is $(2,4)$. This is a contradiction, so we conclude that $\lambda_T(D_1)=\lambda_T(D_2)$. The value has to be $0$ because half-points can only appear after interpolating with the zero function.
\end{proof}

\begin{cor}\label{cor:two_genus_one}
 In the situation of the previous lemma, let $e$ be any edge of $\tropD_T$ between $D_1$ and $D_2$. Then each of the connected components of $\trop(\psi)^{-1}(\tropD_T\setminus e^\circ)$ has genus one.
\end{cor} 
\begin{proof}
Once we know that $\lambda_T$ has slope $0$ between $D_1$ and $D_2$, this follows from a straightforward application of equation \eqref{eq:fundamental}.
\end{proof}

We can now make the following

\begin{dfn}\label{def:rho1}
 Let $\rho_1$ denote the value of $\lambda$ on the special vertex(ices).
\end{dfn}

\begin{lem}
 The subcurve $\tropC_{\geq \rho_1}$ is connected.
\end{lem}
\begin{proof}
 Suppose that there are two connected components $\tropD_1$ and $\tropD_2$. Since $\tropC$ is path-connected, we can find an oriented path $P$ from $\tropD_1$ to $\tropD_2$. The slope of $\lambda$ must be negative at the beginning and positive at the end of $P$, so $\on{div}(\lambda)_{|P}\geq 2$. Restricting equation \eqref{eq:fundamental} to $P$, we find
 \[0=\trop(\psi)^*(D)_{|P}=2p_a(P)-2 +(\on{div}(\lambda)_{|P}+2),\]
 or $p_a(P)\leq -1$, which is a contradiction.
\end{proof}

\begin{lem}\label{lem:ilcorestasopraD}
 The subcurve $\tropC_{\geq \rho_1}$ contains the core.
\end{lem}
\begin{proof}
 Equivalently, $p_a(\tropC_{\geq \rho_1})=2$. This follows from the inspection of Figure \ref{fig:admredcores}, but we give a more formal argument - the following has been suggested by the referee. If $\rho_1=0$ there is nothing to prove, because $\tropC_{\geq 0}=\tropC$, so let us assume that $\rho_1>0$. Let $\tropC_{<\rho_1,\bullet}$ be one of the $h$ connected components of $\tropC\setminus\tropC_{\geq \rho_1}$, let $e_1,\ldots,e_k$ be the edges between  $\tropC_{<\rho_1,\bullet}$ and $\tropC_{\geq \rho_1}$, with positive slope $s_1,\ldots,s_k$. Restricting equation \eqref{eq:fundamental} to $\tropC_{<\rho_1,\bullet}$ we find:
 \[0=2p_a(\tropC_{<\rho_1,\bullet})-2+k+\sum_{i=1}^ks_i,\]
 which holds if and only if $p_a(\tropC_{<\rho_1,\bullet})=0,\ k=1$, and $s_1=1$.
 
 Restricting equation \eqref{eq:fundamental} to $\tropC_{\geq \rho_1}$ we find:
 \[2=2p_a(\tropC_{\geq \rho_1})-2+h-\sum_{j=i}^h1,\]
 which gives the desired result.
\end{proof}

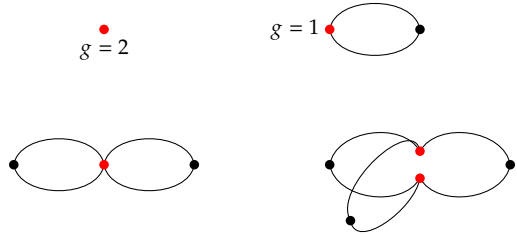
\begin{figure}[htb]
	\centering
		\begin{tikzpicture} [scale=.6]
		\tikzstyle{every node}=[font=\normalsize]
		\tikzset{arrow/.style={latex-latex}}
\coordinate (A) at (0,3,0);
  
 \coordinate (B1) at (-2,3,0);
 \coordinate (B2) at (2,3,0);
  
  \coordinate (A1) at (0,3.3,0);
 \coordinate (A2) at (0,2.7,0);
  \coordinate (B3) at (0,3.3,4);  
    
\fill[red] (A) circle (3pt);
\node at (A) [below] {\scriptsize{$g=2$}};

\begin{scope}[every coordinate/.style={shift={(0,-3,0)}}]
\draw ([c]B1) to [out=80,in=100] ([c]A);
\draw ([c]B1) to [out=-80,in=-100] ([c]A);

\draw ([c]A) to [out=80,in=100] ([c]B2);
\draw ([c]A) to [out=-80,in=-100] ([c]B2);
  
\fill[red] ([c]A) circle (3pt);
\fill ([c]B1) circle (3pt);
\fill ([c]B2) circle (3pt);
\end{scope}

\begin{scope}[every coordinate/.style={shift={(5,0,0)}}]
\draw ([c]A) to [out=80,in=100] ([c]B2);
\draw ([c]A) to [out=-80,in=-100] ([c]B2);
  
\fill[red] ([c]A) circle (3pt);
\fill ([c]B2) circle (3pt);
\node at ([c]A) [below,left] {\scriptsize{$g=1$}};
\end{scope}

\begin{scope}[every coordinate/.style={shift={(7,-3,0)}}]
\draw ([c]B1) to [out=80,in=120] ([c]A1);
\draw ([c]B1) to [out=-80,in=-120] ([c]A2);

\draw ([c]A1) to [out=60,in=100] ([c]B2);
\draw ([c]A2) to [out=-60,in=-100] ([c]B2);
 
 \draw ([c]B3) to [out=120,in=100] ([c]A1);
\draw ([c]B3) to [out=-80,in=-100] ([c]A2);   
  
\fill[red] ([c]A1) circle (3pt);
  \fill[red] ([c]A2) circle (3pt);
 
 \foreach \x in {B1,B2,B3} 
\fill ([c]\x) circle (3pt);

\end{scope}

\end{tikzpicture}
\caption{The top level of $\Delta$ over $\Do\setminus\Dt$.}
\label{top level D1}
\end{figure}

\begin{cor}\label{cor:rho1=rhomax}
 When $\lambda$ is constant on $\tropC_{\geq\rho_1}$, then every positive genus subcurve of $\tropC_{\geq\rho_1}$ supports $D$. The possible shapes of $\tropC_{\geq\rho_1}$ are depicted in Figure \ref{top level D1}.
\end{cor}
\begin{proof}
 Let $F$ be any positive genus subcurve on which $\lambda_{|F}\equiv\rho_1$. Restricting equation \eqref{eq:fundamental} to $F$ we find:
 \[\on{deg}\left(\trop(\psi)^*(D)_{|F}\right)=2p_a(F)-2+\sum_{i=1}^k \left(1+s(\lambda,e_i)\right),\]
 where $e_1,\ldots,e_k$ are the edges from $\tropC\setminus F$ to $F$. It follows from the assumption and the proof above that the slope of $\lambda$ along any $e_i$ can only be $0$ or $-1$. The only possibility for $\on{deg}(\trop(\psi)^*(D)_{|F})$ to be zero remains that of a genus one curve with all outgoing slopes $-1$. But then $\lambda$ would not be constant on $\tropC_{\geq\rho_1}$.
\end{proof}

By construction, there is always at least one vertex of positive weight on the boundary of $\tropD_T^\circ$. We make some stronger statements that we will need below.

\begin{lem}\label{lem:uniquevertexonboundary}
Suppose that $\tropD^\circ$ is connected of genus two. If there is only one vertex of positive weight on the boundary of $\tropD_T^\circ$, then it is the special one.
\end{lem}
\begin{proof}
Let us call  $v_1$ the vertex in the statement. Observe that, by Lemma \ref{lem:D1_D2} and Corollary \ref{cor:two_genus_one}, there has to be a single point of $\tropT$ supporting $D$. We argue by contradiction, showing that if $v_1$ did not support $D$ we could find a different (interpolation of) admissible function(s) greater than $\lambda$; but the latter is defined as the maximum. Extend $\lambda$ with slope $1$ outside $\tropD$. Let $v^\prime$ denote the first vertex of $\tropT$ (not supporting $D$) encountered on the path from $D$ to $v_1$, and let $\ell_1$ denote the distance from $D$ to $v'$. Consider the function $\mu_1$ that is $0$ from $v_1$ to $v'$, has slope $1$ from $v'$ to $D$, and has value $\ell_1$ on $D$ and further away from $v_1$; extend it to $\tropT$ by making it locally constant outside the path from $v_1$ to $D$. The function $\lambda+\mu_1>\lambda$ should still be cut at $\lambda+\mu_1(v_1)=0$, unless one (call it $v_2$) of the vertices of positive weight beyond $D$ had distance $\ell_2<\ell_1$ from $\partial\tropD_T$. In this case, let $\ell_2$ denote the minimal such length. Then consider the function $\mu_2$ that has value $\ell_2$ on $D$ and beyond, and decreases with slope $1$ towards $v_1$, until it reaches $0$ (which happens before $v'$ since $\ell_2<\ell_1$). We observe that $\lambda+\mu_2>\lambda$ should still be cut at $\lambda+\mu_2(v_1)=\lambda+\mu_2(v_2)=0$, and it has two vertices of positive weight on $\partial\tropD_T$ (on two different sides of $D$). This is a contradiction unless $l_1=0$.
\end{proof}

 \begin{example}\label{ese:unique vertex boundary}
 
For the reader's benefit we provide a pictorial proof of Lemma \ref{lem:uniquevertexonboundary} in one example, see Figure \ref{Figure:proofLemmauniquevertex}. In the following pictures, the circle delimits the locus where $\lambda$ is positive. Outside the circle, $\lambda$ can be extended with slope $1$. Conjugate branches are not drawn - we leave it to the reader to complete the picture as necessary.
 
 In the left-most picture - which is the one we start our argument by contradiction with - the values taken by $\lambda$ are the following:
\[
 \lambda(v_1)=0,\;\;\lambda(\text{core})=\ell,\;\;\lambda(v')=\ell-2\ell'=m+2\ell_1,\;\;\lambda(w)=m,\;\;
 \lambda(v_2)=-\ell_2,
 \]
where we have denoted by $w$ the vertex supporting $D$ (it might be arising from the subdivison of an edge), and by $m$ the distance from $w$ to the the circle.
 
 Now we start ``moving'' $D$ towards $v_1$, as shown in the second picture from the left; if $\ell_2$ is long enough, $D$ may reach $v'$ before $v_2$ reaches $\partial\tropD_T$. The value of $\lambda$ stays the same on the core, $v_1$, $v'$, and $v_3$, but it increases otherwise:
 \[\lambda(w)=m+\ell_1,\;\;\;\; \lambda(v_2)=\ell_1-\ell_2.\]

 We can keep pushing $D$ towards $v_1$ until either $D$ reaches $v_1$, or $v_2$ reaches $\partial\tropD_T$, in which case we are no longer in the hypotheses of the Lemma. These are the two possibilities illustrated in the right-most pictures of Figure \ref{Figure:proofLemmauniquevertex}.
 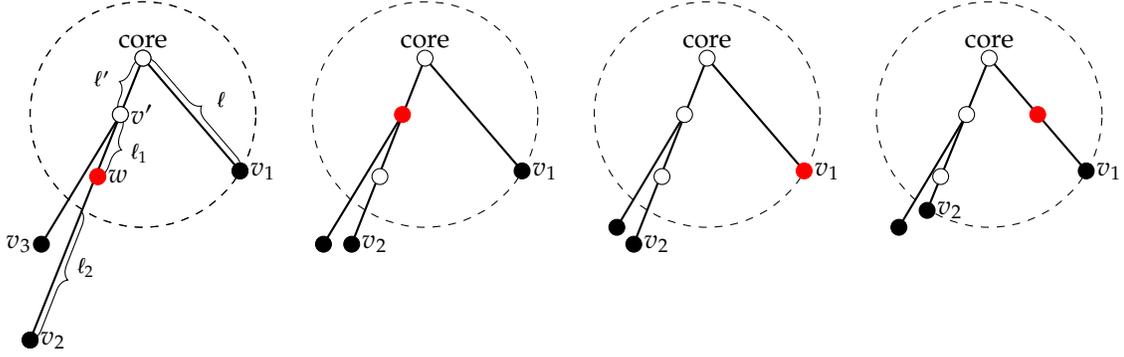
\begin{figure}
  \begin{tikzpicture}[scale=1.5]    
  \tikzset{cross/.style={cross out, draw, thick,
         minimum size=2*(#1-\pgflinewidth), 
         inner sep=1.2pt, outer sep=1.2pt}}

 \coordinate (O) at (0,0);
\coordinate (v1) at (0.86,-0.5);
\coordinate (C) at (0,.5);
\coordinate(v2) at (-1,-2);
\coordinate(v3) at (-.9,-1.15);
\coordinate(D) at (-.4,-.55);
\coordinate(v') at (-.20,0);
\coordinate(DT) at (-.55,-.85);
\coordinate(v31) at (-.8,-1);
\coordinate(v21) at (-.65, -1.15);
\coordinate(v11) at (0.43, 0);
\coordinate(D') at (-.43,-.55);
   \draw[dashed] (0,0)  circle (1cm);
   \draw[fill] (v1)circle (2pt);
     \draw[fill] (v2)circle (2pt);
     \draw[fill] (v3)circle (2pt);
     \draw[fill, red] (D)circle (2pt);
\draw (C) circle (2pt);
    \draw (v') circle (2pt);
     
   \draw[thick] (C)--(v1);
   \draw[thick] (C)--(v2);
   \draw[thick] (v')--(v3);
    \node at (v1) [right] {\small{$v_1$}};
    \node at (v2) [right] {\small{$v_2$}};
    \node at (v3) [left] {\small{$v_3$}};
     \node at (v') [right] {\small{$v'$}};
     \node at (D) [right] {\small{$w$}};
    \node at (C) [above] {\small{core}};
      \draw [decorate,decoration={brace,amplitude=5pt,mirror},xshift=0.2pt,yshift=-0.2pt] (C)--(v') node[black,midway,xshift=-0.4cm,yshift=0.1cm] {\scriptsize{$\ell'$}};
           \draw [decorate,decoration={brace,amplitude=5pt,mirror},xshift=0.2pt,yshift=-0.5pt] (D)--(v') node[black,midway,xshift=0.4cm,yshift=-0.1cm] {\scriptsize{$\ell_1$}};
      
       \draw [decorate,decoration={brace,amplitude=5pt,mirror},xshift=0.2pt,yshift=-0.5pt] (v2)--(DT) node[black,midway,xshift=0.4cm,yshift=0.1cm] {\scriptsize{$\ell_2$}};
               \draw [decorate,decoration={brace,amplitude=5pt,mirror},xshift=0.2pt,yshift=-0.5pt] (v1)--(C) node[black,midway,xshift=0.4cm,yshift=0.2cm] {\scriptsize{$\ell$}};

   \draw[dashed] (0,0)  circle (1cm);
   \draw[fill] (v1)circle (2pt);
     \draw[fill] (v2)circle (2pt);
     \draw[fill, red] (D)circle (2pt);
\draw[fill=white] (C) circle (2pt);
    \draw[fill=white] (v') circle (2pt);
    

\begin{scope}[every coordinate/.style={shift={(2.5,0,0)}}]

   \draw[thick] ([c]C)--([c]v1);
   \draw[thick] ([c]C)--([c]v21);
   \draw[thick] ([c]v')--([c]v3);
    \node at ([c]v1) [right] {\small{$v_1$}};
    \node at ([c]v21) [right] {\small{$v_2$}};
    \node at ([c]C) [above] {\small{core}}; 
    
   \draw[dashed] ([c]O)  circle (1cm);
   \draw[fill] ([c]v1)circle (2pt);
     \draw[fill] ([c]v21)circle (2pt);
     \draw[fill, red] ([c]v')circle (2pt);
\draw[fill=white] ([c]C) circle (2pt);
    \draw[fill=white] ([c]D) circle (2pt);
     \draw[fill] ([c]v3)circle (2pt);
\end{scope}

\begin{scope}[every coordinate/.style={shift={(7.5,0,0)}}]
 
   \draw[thick] ([c]C)--([c]v1);
   \draw[thick] ([c]C)--([c]DT);
   \draw[thick] ([c]v') -- ([c]v31);
    \node at ([c]v1) [right] {\small{$v_1$}};
    \node at ([c]DT) [right] {\small{$v_2$}};
    \node at ([c]C) [above] {\small{core}};
    
          
   \draw[dashed] ([c]O)  circle (1cm);
   \draw[fill] ([c]v1)circle (2pt);
\draw[fill] ([c]v31)circle (2pt);
\draw[fill] ([c]DT)circle (2pt);
     \draw[fill, red] ([c]v11)circle (2pt);
\draw[fill=white] ([c]C) circle (2pt);
\draw[fill=white] ([c]v') circle (2pt);
    \draw[fill=white] ([c]D') circle (2pt);
     
\end{scope}

\begin{scope}[every coordinate/.style={shift={(5,0,0)}}]

   \draw[thick] ([c]C)--([c]v1);
   \draw[thick] ([c]C)--([c]v21);
    \node at ([c]v1) [right] {\small{$v_1$}};
    \node at ([c]v21) [right] {\small{$v_2$}};
    \node at ([c]C) [above] {\small{core}};
      \draw[thick] ([c]v') -- ([c]v31);
      
   \draw[dashed] ([c]O)  circle (1cm);
   \draw[fill=white] ([c]v')circle (2pt);
     \draw[fill] ([c]v21)circle (2pt);
     \draw[fill, red] ([c]v1)circle (2pt);
     \draw[fill] ([c]v31)circle (2pt);
\draw[fill=white] ([c]C) circle (2pt);
    \draw[fill=white] ([c]D) circle (2pt);
\end{scope}
  \end{tikzpicture}
\caption{A graphical proof of Lemma \ref{lem:uniquevertexonboundary}.}
\label{Figure:proofLemmauniquevertex}
 \end{figure}
\end{example}

The following statements are by-products of the proof given above:

\begin{cor}\label{cor:Dpositiveweight}
Suppose that $\tropD^\circ$ is connected of genus two.  If $D$ is supported on the boundary of $\tropD_T^\circ$, the corresponding vertex has positive weight.
\end{cor}

\begin{cor}\label{cor:oneitherside}
Suppose that $\tropD^\circ$ is connected of genus two. If $D$ is supported in $\tropD^\circ_T$, then there is at least one vertex of positive weight on the boundary of $\tropD_T^\circ$ on either side of $D$.
\end{cor}

Finally, we deal with a genus one situation.

\begin{lem}\label{lem:two_genus_one_pos_wt}
 Let $\tropD^\circ_1$ be a connected component of $\tropD^\circ$ of genus one. Then there is at least one vertex of positive weight on its boundary.
\end{lem}
\begin{proof}
 Let $E_1$ denote the genus one core of $\tropD^\circ_1$. There is an (at least one) admissible function $\overline\lambda_i$ behaving like the distance from $E_1$: it is the one with $D$ supported on (any vertex of) the other genus one subcurve $E_2$ of $\tropC$. Then the vertex of positive weight closest to $E_1$ can be at most as far as $E_2$, because otherwise $\tropD^\circ$ would be connected of genus two. Now $\lambda$ looks like $\lambda_i$ in a neighbourhood of $E_1$.
\end{proof}
For the sake of concreteness, we notice that, in the situation of the above lemma, either $E_2$ is on the boundary of $\tropD^\circ_1$, supports $D$, and has positive weight (essentially by Corollary \ref{cor:Dpositiveweight}), or the path from $E_1$ to $E_2$ crosses the boundary of $\tropD^\circ_1$ at a rational vertex supporting $\frac{1}{2}D$, and proceeds with slope $0$ in a neighbourhood.

\subsection{The secondary alignment}\label{sec:secalign}

\begin{lem}
 The locus $V$ where $w(\tropD)$ is $1$, or it is $2$ but supported on two non-conjugate vertices of $\tropC$, is closed in $\TAt$.
\end{lem}
\begin{proof}
 The statement is of a topological nature. The loci where the combinatorial type of $\tropC$ is constant form a constructible stratification of $\TAt$. By Noetherianity, it is enough to check that the described locus is closed under specialisation for any trait, i.e. if we have an edge contraction $\tropC_s\to\tropC_\eta$ such that the latter is in $V$, then the former must be as well. This is easily checked.
\end{proof}

\begin{dfn}\label{def:niceopen}
 Let $\TAt^\circ$ be the open substack of $\TAt$ defined as the complement of the locus $V$ described in the previous Lemma.
\end{dfn}

\begin{dfn}\label{def:radial}
 Let $\Sigma''$ be the subdivision of $\on{trop}(\TAt^\circ)$ defined by
 \begin{itemize}
  \item formally contracting $\tropD$ to a point $z$;
  \item aligning the rest of the curve with respect to the distance from $z$.
 \end{itemize}
Let $(\TAt^\circ)^\text{pre}$ denote the corresponding logarithmic blow-up.
\end{dfn}

On $\TAt^\circ$ the weight of $\tropD$ is at least $2$. If it is $2$, it is supported on a single vertex $v_1$ of $\tropD_T$. In this case, by Lemma \ref{lem:uniquevertexonboundary}, $\lambda$ has $D$ supported on $v_1$. By relabelling, assume that $\overline \lambda_1$ was the admissible function (Definition \ref{def:admissiblefunction}) determined by having $v_1$ as the special vertex. By Lemma \ref{lem:ilcorestasopraD}, the slope of $\overline \lambda_1$ is always $1$ below the level of $v_1$, hence the alignment of the previous definition makes sense. We extend $\lambda$ to reach the next vertex of positive weight.

\begin{dfn}
 Let $S$ be a geometric point of $(\TAt^\circ)^{\text{pre}}$, and suppose that $w(\tropD)=2$.  Let $\lambda_1^\prime\in \Gamma(S,\pi_*\overline M_C)$ be the lift of $\overline \lambda_1$ determined by $w(\tropC_{>0})\leq 2$ and $w(\tropC_{\geq0})\geq 3$. Let \[\widetilde{\lambda}=\max\{0,\lambda_1'\}.\]
\end{dfn}

\begin{lem}
 If a specialisation $\eta \rightsquigarrow s$ in $(\TAt^\circ)^\text{pre}$ induces an edge contraction $\tropC_s\to\tropC_\eta$ such that $w(\tropC_{s,\geq 0})=2$ but $w(\tropC_{\eta,\geq 0})\geq 3$, then $\widetilde{\lambda}$ generises to $\lambda$.
\end{lem}
\begin{proof}
 By Lemma \ref{lem:uniquevertexonboundary}, on the cone of $\on{trop}(\TAt^\circ)$ corresponding to the point $s$, the function $\lambda$ was a lift of $\overline{\lambda}_1$, with a cutoff at the first vertex $v_1$ of positive weight. Extending it to the second vertex $v_2$ of positive weight makes sense thanks to Definition \ref{def:radial}. When an edge contraction makes $w(\tropD)$ grow, it means in particular that $v_2$ and $v_1$ become at the same height, so $\widetilde{\lambda}$ coincides with the original $\lambda$.
\end{proof}

This shows that $\widetilde{\lambda}$ extends $\lambda$ as a well-defined (continuous) PL function on a subdivision $\doublewidetilde{\tropC}$ of $\widetilde{\tropC}$ for any $S$-point of $(\TAt^\circ)^\text{pre}$, which we think of as an extension of $\lambda$ on certain cones of $\on{trop}((\TAt^\circ)^\text{pre})$.

\begin{dfn}\label{def:doubletildeA}
 Let $\doublewidetilde{\mathcal A_2}$ denote the logarithmic blow-up of $\TAt^\circ$ induced by applying universal semistable reduction to $\doublewidetilde{\tropC}\to\on{trop}((\TAt^\circ)^\text{pre})$.
\end{dfn}

\begin{rem}
 The spaces $(\TAt^\circ)^\text{pre}$ and $\doublewidetilde{\mathcal A_2}$ are analogous to, respectively, the space of radially aligned and the space of centrally aligned logarithmic curves in \cite{RSPW1}.
\end{rem}

Similarly to the previous section, we have the following result.

\begin{thm}\label{thm:doubletildeA}
 The moduli space $\doublewidetilde{\mathcal A_2}$ is a logarithmically smooth stack with locally free logarithmic structure (and therefore smooth).
 There is a logarithmically smooth curve $\TTC\to \doublewidetilde{\mathcal A_2}$, which is a partial destabilisation of $\tilde{a}^*\mathcal{C}$, such that $\doublewidetilde{\tropC}$ is its tropicalisation, $\widetilde{\lambda}\in\Gamma(\TTC,\overline M_{\TTC})$, and all of its values at vertices of $\doublewidetilde{\tropC}$ are comparable.
\end{thm}

\subsection{Examples of subdivisions}
In this section we collect a few examples to show how the subdivisions  $\Sigma$ and $\Sigma^\prime$ of a cone $\sigma\in\operatorname{trop}(\mathcal{A}_2^{\text{wt}})$,  the combinatorics of $\lambda$ on the various cones, and the associated singularities in $\overline{\mathcal C}$ look like. The construction of $\OC$ will be carried out in the next section.

\begin{exa}
 The stabilisation of $C$ has smooth core of weight zero, and two rational tails of high weight, one of which attached  to a Weierstrass point (compare with Example \ref{ese:lambdabar}).
In this case, the subdivisions $\Sigma$ and $\Sigma^\prime$ coincide. Notice that there is a simplicial non-smooth cone $\sigma_3$; correspondingly, half edge-lengths occur in a Kummer extension of $M_{\widetilde{A}}$. See Figure \ref{esesub1}.  
\begin{figure}[h]
\begin{tikzpicture} [scale=.4]
\tikzstyle{every node}=[font=\normalsize]
		\tikzset{arrow/.style={latex-latex}}
		
	\tikzset{cross/.style={cross out, draw, thick,
         minimum size=2*(#1-\pgflinewidth), 
         inner sep=1.2pt, outer sep=1.2pt}}


		\coordinate (Oc) at (0,3,0);
		\coordinate (T) at (-2,4,0);
		\coordinate (bT) at (-2,2,0);
		\coordinate (W) at (3,1.5,0);
		\coordinate (Ot) at (0,0,0);
		\coordinate (Tt) at (-2,0,0);
		\coordinate (Wt) at (3,0,0);
		
		\coordinate (P1) at (-.5,-1,0);
		\coordinate (P2) at (0,-1,0);
		\coordinate (P3) at (.5,-1,0);
		\coordinate (P4) at (1,-1,0);
		\coordinate (P5) at (1.5,-1,0);
    	\coordinate (P6) at (3.2,-.8,0);

		\begin{scope}[every coordinate/.style={shift={(2,5,0)}}]

     	
     	\draw[dotted] ([c]Oc)-- ([c]T);
     	\draw ([c]Oc)-- ([c]bT);
     	\draw ([c]Oc)-- ([c]W);
     	\draw ([c]Ot)-- ([c]Tt);
     	\draw ([c]Ot)-- ([c]Wt);
     	
     	\draw[blue] ([c]Ot)--([c]P1);
     	\draw[blue] ([c]Ot)--([c]P2);
     	\draw[blue] ([c]Ot)--([c]P3);
     	\draw[blue] ([c]Ot)--([c]P4);
     	\draw[blue] ([c]Ot)--([c]P5);
     	\draw[blue] ([c]Wt)--([c]P6);

       
       \node at (1,8.5,0) [above] {\scriptsize{$l_1$}};
        \node at (.9,8.4,0) [below] {\scriptsize{$l_1$}};
        \node at (3.5,7.2,0) [above] {\scriptsize{$l_2$}};
        
       \node at (1,5,0) [above] {\scriptsize{$l_1$}};
        \node at (3.5,5,0) [above] {\scriptsize{$2l_2$}};
        
      \node at (0,9,0) [left] {\scriptsize{$\overline{T}$}};
        \node at (0,7,0) [left] {\scriptsize{$T$}};
         \node at (5,6.2,0) [right] {\scriptsize{$W$}};

       \draw[dashed] ([c]Oc)--([c]Ot);
       \draw[dashed] ([c]T)--([c]Tt);
       \draw[dashed] ([c]W)--([c]Wt);
       
		\draw[fill=white]  ([c]Oc) circle (3pt);
		\draw  ([c]T) node [cross]{};
		\draw[fill=white]  ([c]Ot) circle (3pt);
	   \foreach \x in {bT,W,Tt,Wt} 
          \fill ([c]\x) circle (3pt);
       
         \end{scope}
		\coordinate (O) at (0,0,0);
		\coordinate (l2) at (6,0,0);
		\coordinate (l1) at (0,6,0);
		\coordinate (r1) at (3,6,0);
		\coordinate (r2) at (5,5,0);
		\coordinate (r3) at (6,2,0);
		
		\begin{scope}[every coordinate/.style={shift={(0,-8,0)}}]
		
		\draw[->] ([c]O)--([c]l2);
		\draw[->] ([c]O)--([c]l1);
		\draw[green] ([c]O)--([c]r1);
		\draw[orange]([c]O)--([c]r2);
		\draw[red] ([c]O)--([c]r3);
		
		 \node at (.7,-4.5,0) [above] {\scriptsize{$\sigma_1$}};
		 \node at (3,-2,0) [above] {\begin{color}{green}\scriptsize{$\tau_1$}\end{color}\scriptsize{, $l_2=2l_1$}};
		 \node at (2.3,-5,0) [above] {\scriptsize{$\sigma_2$}};
		 \node at (5,-3,0) [above right] {\begin{color}{orange}\scriptsize{$\tau_2$}\end{color}\scriptsize{, $l_2=l_1$}};
		 \node at (3.4,-6.2,0) [above] {\scriptsize{$\sigma_3$}};
		 \node at (6,-6,0) [above right] {\begin{color}{red}\scriptsize{$\tau_3$}\end{color}\scriptsize{, $l_1=3l_2$}};
		 \node at (3.8,-7.7,0) [above] {\scriptsize{$\sigma_4$}};
		 \node at ([c]l2) [right] {\scriptsize{$l_1$}};
		  \node at ([c]l1) [left] {\scriptsize{$l_2$}};
		  \end{scope}

		   \coordinate (O2) at (-.5,.5,0);
		   \coordinate (B1) at (-.7,-1.5,0);
		   \coordinate (C1) at (-2.9,-1.5,0);
		  \coordinate (B2) at (1,-1.5,0);
		  \coordinate (B2l) at (2.3,-1.5,0);
		  \coordinate(C2) at (3.5,-1.5,0);
		  \coordinate (O22) at (-.5,.8,0);
		  \coordinate (B12) at (-.9,-.8,0);
		  \coordinate (C12) at (-2.3,-1.5,0);
		  
		   \coordinate (O3) at (0,.9,0);
		   \coordinate (D2) at (.5,-.7,0);

		  \coordinate (Bt1) at (-1,-1.5,0);
	      
	          \coordinate (D22) at (.5,-1.5,0);
	      
	     \coordinate (y1) at (6,0,0);
	   \coordinate (y2) at (6,-3,0);
	    \coordinate (x1) at (4.5,-1.5,0);
	   \coordinate (x2) at (7.5,-1.5,0);
	   \coordinate (oo) at (6,-1.5,0);
	    \coordinate (m1) at (7,-.5,0);
	   \coordinate (m2) at (5,-.5,0);
	    \coordinate (b1) at (5.7,-.4,0);
	    \coordinate (b2) at (7,.2,0);
	    
	       \coordinate (y1e) at (6.2,0,0);
	   \coordinate (y2e) at (6.2,-3,0);
	   \coordinate (b1e) at (5.7,-2.7,0);
	    \coordinate (b2e) at (7,-2.2,0);
	       \coordinate (m3) at (7,-2.5,0);
		  
		  \begin{scope}[every coordinate/.style={shift={(15,12,0)}}]
		  
		  \node at (12,12,0) [left] {{$\sigma_1:$}};

          \draw ([c]O2)-- ([c]B2);
          \draw ([c]O2)-- ([c]B1);
          \draw ([c]B1)-- ([c]C1);
          \draw ([c]B2)-- ([c]B2l);
          \draw ([c]B2l)-- ([c]C2);
          
          \node at (14.5,11.7,0) [left] {{\begin{color}{blue}\scriptsize{$2$}\end{color}}};
           \node at (15.3,11.7,0) [right] {{\begin{color}{blue}\scriptsize{$1$}\end{color}}};
          
          \draw [decorate,decoration={brace,amplitude=5pt,mirror},xshift=0.4pt,yshift=-0.4pt] ([c]O2)--([c]B2) node[black,midway,xshift=-0.1cm,yshift=-0.8cm] {\scriptsize{$2l_1$}};
           \draw [decorate,decoration={brace,amplitude=5pt,mirror},xshift=0.4pt,yshift=-0.4pt] ([c]B2)--([c]B2l) node[black,midway,xshift=0.2cm,yshift=-0.4cm] {\scriptsize{$l_2-2l_1$}};
           
           \draw[gray] ([c]y1)-- ([c]y2);
            \draw[thick] ([c]x1)-- ([c]x2);
              \draw[thick] ([c]b1)-- ([c]b2);
              \draw[gray] ([c]oo) parabola ([c]m1);
               \draw[gray] ([c]oo) parabola ([c]m2);
            
            \node at ([c]m1) [right] {{\scriptsize{$\overline{T}$}}};
             \node at ([c]x2) [right] {{\scriptsize{$T$}}};
              \node at ([c]b2) [right] {{\scriptsize{$W$}}};
              
                
		  \draw[fill=white]  ([c]O2) circle (3pt);

          \fill ([c]B2l) circle (3pt);
            \draw[fill=red]([c]B1) circle (3pt);
            
          \draw ([c]B2) node[cross]{};
              
            \end{scope}

		  \begin{scope}[every coordinate/.style={shift={(15,8,0)}}]
		  
		  \node at (12,8,0) [left] {{\begin{color}{green} $\tau_1:$\end{color}}};

          \draw ([c]O2)-- ([c]B1);
          \draw ([c]O2)-- ([c]B2);
          \draw ([c]B1)-- ([c]C1);
          \draw ([c]B2)-- ([c]B2l);

          \node at (14.5,7.7,0) [left] {{\begin{color}{blue}\scriptsize{$2$}\end{color}}};
           \node at (15.3,7.7,0) [right] {{\begin{color}{blue}\scriptsize{$1$}\end{color}}};

           \draw[ thick] ([c]y1)-- ([c]y2);
            \draw[thick] ([c]x1)-- ([c]x2);
            
              \draw [gray]([c]oo) parabola ([c]m1);
               \draw[gray] ([c]oo) parabola ([c]m2); 
            
            \node at ([c]m1) [right] {{\scriptsize{$\overline{T}$}}};
             \node at ([c]x2) [right] {{\scriptsize{$T$}}};
              \node at ([c]y2) [right] {{\scriptsize{$W$}}};
              
                
		  \draw[fill=white]  ([c]O2) circle (3pt);

          \fill ([c]B2) circle (3pt);
           \draw[fill=red] ([c]B1) circle (3pt);
           
            \end{scope}

		 \begin{scope}[every coordinate/.style={shift={(15,4,0)}}]
		  
		  \node at (12,4,0) [left] {{$\sigma_2:$}};

          \draw ([c]O22)-- ([c]B2l);
          \draw ([c]O22)-- ([c]B12);
          \draw ([c]B12)-- ([c]C12);
           \draw ([c]C12)-- ([c]C1);
          \draw ([c]B2l)-- ([c]C2);
        

            \node at (14.4,4,0) [left] {{\begin{color}{blue}\scriptsize{$2$}\end{color}}};
              \node at (13.9,3.3,0) [left] {{\begin{color}{blue}\scriptsize{$1$}\end{color}}};
              \node at (15.5,4,0) [right] {{\begin{color}{blue}\scriptsize{$1$}\end{color}}};

          \draw [decorate,decoration={brace,amplitude=5pt,mirror},xshift=-0.4pt] ([c]B12)--([c]O22) node[black,midway,xshift=0.4cm,yshift=-0.85cm] {\scriptsize{$l_2-l_1$}};
           \draw [decorate,decoration={brace,amplitude=5pt,mirror},xshift=-0.4pt,yshift=-0.4pt] ([c]C12)--([c]B12) node[black,midway,xshift=-0.2cm,yshift=-0.4cm] {\scriptsize{$2l_1-l_2$}};

           \draw[gray] ([c]y1)-- ([c]y1e)--([c]y2e)--([c]y2)--([c]y1);
            \draw[thick] ([c]b1)-- ([c]b2);
            \draw[thick] ([c]b1e)-- ([c]b2e);

        
             \node at ([c]b2) [right] {{\scriptsize{$T$}}};
              \node at ([c]b2e) [right] {{\scriptsize{$W$}}};
              
                
		  \draw[fill=white]  ([c]O22) circle (3pt);
		  
			\foreach \x in {C12,B2l} 
          \fill ([c]\x) circle (3pt);
          
         \draw ([c]B12) node[red,cross]{};

            \end{scope}

		\begin{scope}[every coordinate/.style={shift={(15,0,0)}}]
		  
		 \node at (12,0,0) [left] {{\begin{color}{orange} $\tau_2:$\end{color}}};

           \draw ([c]O)-- ([c]B2);
          \draw ([c]O)-- ([c]Bt1);
          \draw ([c]Bt1)-- ([c]C1);
       \draw ([c]B2)-- ([c]C2);
        

            \node at (14.5,-.7,0) [left] {{\begin{color}{blue}\scriptsize{$1$}\end{color}}};
            \node at (15.5,-.7,0) [right] {{\begin{color}{blue}\scriptsize{$1$}\end{color}}};

           \draw[gray] ([c]y1)-- ([c]y1e)--([c]y2e)--([c]y2)--([c]y1);
            \draw[thick] ([c]b1)-- ([c]b2);
            \draw[thick] ([c]b1e)-- ([c]b2e);

        
             \node at ([c]b2) [right] {{\scriptsize{$T$}}};
              \node at ([c]b2e) [right] {{\scriptsize{$W$}}};
              
              \draw[fill=red]  ([c]O) circle (3pt);
		  
			\foreach \x in {Bt1,B2} 
          \fill ([c]\x) circle (3pt);

            \end{scope}
		  \begin{scope}[every coordinate/.style={shift={(15,-4,0)}}]
		  
		  \node at (12,-4,0) [left] {{$\sigma_3:$}};

          \draw ([c]O3)-- ([c]Bt1);
          \draw ([c]O3)-- ([c]D2);
          \draw ([c]Bt1)-- ([c]C12);
          \draw ([c]D2)-- ([c]B2l);
          \draw ([c]B2l)-- ([c]C2);
          
          \node at (14.5,-4.5,0) [left] {{\begin{color}{blue}\scriptsize{$1$}\end{color}}};
           \node at (15.6,-3.5,0) {{\begin{color}{blue}\scriptsize{$3$}\end{color}}};
            \node at (15.9,-5.05,0) [above right] {{\begin{color}{blue}\scriptsize{$1$}\end{color}}};
          
          \draw [decorate,decoration={brace,amplitude=5pt},xshift=0.9pt,yshift=-0.4pt] ([c]D2)--([c]O3) node[black,midway,xshift=-0.1cm,yshift=-0.95cm] {\scriptsize{$\frac{l_1-l_2}{2}$}};
           \draw [decorate,decoration={brace,amplitude=5pt,mirror},xshift=0.2pt,yshift=-0.3pt] ([c]D2)--([c]B2l) node[black,midway,xshift=0.15cm,yshift=-0.47cm] {\scriptsize{$\frac{3l_2
           -l_1}{2}$}};

           \draw[gray] ([c]y1)-- ([c]y1e)--([c]y2e)--([c]y2)--([c]y1);
            \draw[thick] ([c]b1)-- ([c]b2);
            \draw[thick] ([c]b1e)-- ([c]b2e);

        
             \node at ([c]b2) [right] {{\scriptsize{$T$}}};
              \node at ([c]b2e) [right] {{\scriptsize{$W$}}};
              
             
		  \draw[fill=white]  ([c]O3) circle (3pt);

          \fill ([c]Bt1) circle (3pt);
           \fill ([c]B2l) circle (3pt);

             \draw ([c]D2) node[cross,red]{};

            \end{scope}
	
		  \begin{scope}[every coordinate/.style={shift={(15,-8,0)}}]
		  
		  \node at (12,-8,0) [left] {{\begin{color}{red} $\tau_3:$\end{color}}};

          \draw ([c]O3)-- ([c]Bt1);
          \draw ([c]O3)-- ([c]D22);
          \draw ([c]Bt1)-- ([c]C1);
          \draw ([c]D22)-- ([c]B2l);

          \node at (14.5,-8.3,0) [left] {{\begin{color}{blue}\scriptsize{$1$}\end{color}}};
           \node at (15.3,-8.3,0) [right] {{\begin{color}{blue}\scriptsize{$3$}\end{color}}};

           \draw[ thick] ([c]y1)-- ([c]y2);
            
              \draw[thick] ([c]oo) parabola ([c]m1);
               \draw[thick] ([c]oo) parabola ([c]m3); 
            
            \node at ([c]m1) [right] {{\scriptsize{$W$}}};
          
              \node at ([c]y2) [right] {{\scriptsize{$T$}}};
              
		  \draw[fill=white]  ([c]O3) circle (3pt);

          \fill ([c]Bt1) circle (3pt);
           \draw[fill=red] ([c]D22) circle (3pt);

            \end{scope}
            
		  \begin{scope}[every coordinate/.style={shift={(15,-12,0)}}]
		  
		  \node at (12,-12,0) [left] {{$\sigma_4:$}};

          \draw ([c]O3)-- ([c]Bt1);
          \draw ([c]O3)-- ([c]D22);
          \draw ([c]Bt1)-- ([c]C12);
           \draw ([c]C12)-- ([c]C1);
          \draw ([c]D22)-- ([c]B2l);

          \node at (14.5,-12.2,0) [left] {{\begin{color}{blue}\scriptsize{$1$}\end{color}}};
           
           \node at (15.3,-12.2,0) [right] {{\begin{color}{blue}\scriptsize{$3$}\end{color}}};

          \draw [decorate,decoration={brace,amplitude=5pt,mirror},xshift=0.2pt,yshift=-0.4pt] ([c]Bt1)--([c]O3) node[black,midway,xshift=0.2cm,yshift=-.88cm] {\scriptsize{$3  l_2$}};
           \draw [decorate,decoration={brace,amplitude=5pt},xshift=0.2pt,yshift=-0.3pt] ([c]Bt1)--([c]C12) node[black,midway,xshift=-.2cm,yshift=-0.4cm] {\scriptsize{$l_1 - 3l_2$}};

           \draw[gray] ([c]y1)-- ([c]y2);
            \draw[thick] ([c]b1)-- ([c]b2);
            
              \draw[thick] ([c]oo) parabola ([c]m1);
               \draw[thick] ([c]oo) parabola ([c]m3); 
            
            \node at ([c]m1) [right] {{\scriptsize{$W$}}};
          
              \node at ([c]b2) [right] {{\scriptsize{$T$}}};
          
           \fill ([c]C12) circle (3pt);
           \draw[fill=red] ([c]D22) circle (3pt);
           
          \draw[fill=white]  ([c]O3) circle (3pt);
         \draw ([c]Bt1) node[cross]{};

          \end{scope}

	\end{tikzpicture}
	\caption{Tropical admissible cover (top left); the subdivision of $\operatorname{trop}(\mathcal A_2^{\text{wt}})$ (bottom left);  $\lambda_T$  and the  Gorenstein singularities (right).}
 \label{esesub1}
\end{figure}
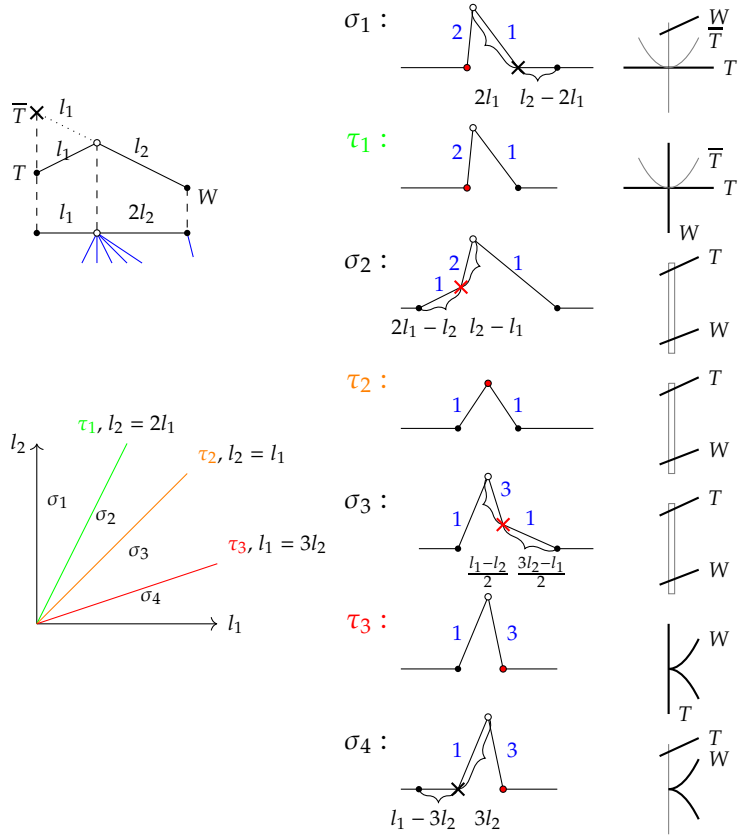
\end{exa}

\begin{exa}
The stabilisation of $C$ has a smooth core of weight zero, and three rational tails attached to general points. Aligning with respect to all the admissible functions produces the non-simplicial subdivision (a); the subdivisions (b) and (c) are the coarsening in case of high-weight, respectively weight-two, tails.
We also represent $\lambda_T$ and the associated singularity on some cones of the subdivision; the other ones can be derived by symmetry. In case (c) the singularity over $\sigma_2$ is replaced by a (sprouted) ribbon. See Figure \ref{esesub2}.

\begin{figure}[h]
\begin{tikzpicture} [scale=.4]
\tikzstyle{every node}=[font=\normalsize]
		\tikzset{arrow/.style={latex-latex}}
		
	\tikzset{cross/.style={cross out, draw, thick,
         minimum size=2*(#1-\pgflinewidth), 
         inner sep=1.2pt, outer sep=1.2pt}}


		\coordinate (Oc) at (0,3,0);
		\coordinate (T1) at (-2,4,0);
		\coordinate (bT1) at (-2,2,0);
		\coordinate (T2) at (4,1.5,0);
		\coordinate (bT2) at (4,2.5,0);
		\coordinate (T3) at (5,6.5,5);
		\coordinate (bT3) at (5,7.5,5);
		\coordinate (Ot) at (0,0,0);
		\coordinate (T1t) at (-2,0,0);
		\coordinate (T2t) at (4,0,0);
		\coordinate (T3t) at (5,0,5);
		
		\coordinate (P1) at (-.5,-1,0);
		\coordinate (P2) at (0,-1,0);
		\coordinate (P3) at (.5,-1,0);
		\coordinate (P4) at (1,-1,0);
		\coordinate (P5) at (1.8,-.7,0);
    	\coordinate (P6) at (-1,-.8,0);
		
		
		\coordinate (l1) at (-1.2,2.3,0);
		\coordinate (bl1) at (-.9,4.4,0);
		\coordinate (l2) at (2.5,1.8,0);
		\coordinate (bl2) at (2.5,2.5,0);
		\coordinate (l3) at (3.6,4.4,3);
		\coordinate (bl3) at (3,5.8,3);
		\coordinate (Ot) at (0,0,0);
		\coordinate (l1t) at (-1,0.2,0);
		\coordinate (l2t) at (1.5,0,0);
		\coordinate (l3t) at (3.5,0,3);
			
		\begin{scope}[every coordinate/.style={shift={(-6,10,0)}}]

     	
     	\draw[dotted] ([c]Oc)-- ([c]T1);
     	\draw ([c]Oc)-- ([c]bT1);
     	\draw ([c]Oc)-- ([c]T2);
     	\draw[dotted] ([c]Oc)-- ([c]bT2);
     	\draw ([c]Oc)-- ([c]T3);
     	\draw[dotted] ([c]Oc)-- ([c]bT3);
     	
     	\draw ([c]Ot)-- ([c]T3t);
     	\draw ([c]Ot)-- ([c]T2t);
     	\draw ([c]Ot)-- ([c]T1t);
     	
     	\draw[blue] ([c]Ot)--([c]P1);
     	\draw[blue] ([c]Ot)--([c]P2);
     	\draw[blue] ([c]Ot)--([c]P3);
     	\draw[blue] ([c]Ot)--([c]P4);
     	\draw[blue] ([c]Ot)--([c]P5);
     	\draw[blue] ([c]Ot)--([c]P6);

       
       \node at ([c]l1) [above] {\scriptsize{$l_1$}};
        \node at ([c]l2) [above] {\scriptsize{$l_2$}};
           \node at ([c]l3) [above=.4cm] {\scriptsize{$l_3$}};
        
      \node at ([c]l1t) [above] {\scriptsize{$l_1$}};
       \node at ([c]l2t) [above] {\scriptsize{$l_2$}};
     
           \node at ([c]l3t) [above] {\scriptsize{$l_3$}};
        
      \node at ([c]T1) [left] {\scriptsize{$\overline{T}_1$}};
        \node at ([c]bT1) [left] {\scriptsize{$T_1$}};
         
        \node at ([c]bT2) [right] {\scriptsize{$\overline{T}_2$}};
        \node at ([c]T2) [right] {\scriptsize{$T_2$}};
         
         \node at ([c]bT3) [right] {\scriptsize{$\overline{T}_3$}};
        \node at ([c]T3) [right] {\scriptsize{$T_3$}};
       
       \draw[dashed] ([c]Oc)--([c]Ot);
       \draw[dashed] ([c]T1)--([c]T1t);
       \draw[dashed] ([c]bT2)--([c]T2t);
       \draw[dashed] ([c]bT3)--([c]T3t);
      
		\draw[fill=white]  ([c]Oc) circle (3pt);
			\draw  ([c]T1) node [cross]{};
			\draw  ([c]bT2) node [cross]{};
    \draw  ([c]bT3) node [cross]{};

    			\draw  ([c]Ot) circle (3pt);

						\foreach \x in {bT1,T2,T3,T1t,T2t,T3t} 
          \fill ([c]\x) circle (3pt);
      
         \end{scope}
         \coordinate(label) at (0,-.5,0);
         \coordinate (A) at (-5,0,0);
         \coordinate (A1) at (-2.5,0,0);
         \coordinate (A2) at (0,0,0);
        \coordinate (A3) at (2.5,0,0);
        \coordinate (B) at (5,0,0);
         \coordinate (B1) at (2.5,4.3,0);
         \coordinate (B2) at (3.75,2.15,0);
         \coordinate (B3) at (1.25,6.45,0);
         \coordinate (C) at (0,8.6,0);
         \coordinate (C1) at (-2.5,4.3,0);
         \coordinate (C2) at (-3.75,2.15,0);
         \coordinate (C3) at (-1.25,6.45,0);
         \coordinate (M1) at (0,2.87,0);
           \coordinate (M2) at (0,1.23,0);
          \coordinate (M3) at (1.43,3.69,0);
        \coordinate (M3b) at (-1.43,3.69,0);
        
         \coordinate (MM) at (0,2,0);
          \coordinate (MA1) at (-2.6,1.5,0);
         \coordinate (MA2) at (-1.8,.8,0);
         \coordinate (CoB) at (-2,1.6,0);
         \coordinate (CoA) at (-.3,.8,0);
         
       \begin{scope} [every coordinate/.style={shift={(-20,-2,0)}}]
       
        \node at ([c]label)[below]  {(a)};
        
        \draw ([c]A)--([c]B);
         \draw ([c]A)--([c]C);
          \draw ([c]C)--([c]B);
          
         \draw([c]C)--([c]A1); 
         \draw([c]C)--([c]A2); 
        \draw([c]C)--([c]A3); 
        
        \draw([c]A)--([c]B1); 
         \draw([c]A)--([c]B2); 
        \draw([c]A)--([c]B3); 
        
        \draw ([c]B)--([c]C1); 
         \draw([c]B)--([c]C2); 
        \draw([c]B)--([c]C3); 
        
        \node at ([c]C) [above] {\scriptsize{$l_3$}};
        \node at ([c]A) [left] {\scriptsize{$l_1$}};
        \node at ([c]B) [right] {\scriptsize{$l_2$}};
        
          \fill ([c]M1) circle (2pt);
           \fill ([c]M2) circle (2pt);
             \fill ([c]M3) circle (2pt);
              \fill ([c]M3b) circle (2pt);
        
        \end{scope}
        
        \begin{scope} [every coordinate/.style={shift={(-6,-2,0)}}]
        
        \node at ([c]label)[below]  {(b)};
        
        \draw ([c]A)--([c]B);
         \draw ([c]A)--([c]C);
          \draw ([c]C)--([c]B);
          
          \draw[red]([c]M1)--([c]A);
           \draw[]([c]M1)--([c]B);
            \draw[orange]([c]M2)--([c]A);
           \draw[]([c]M2)--([c]B);
          \draw[green] ([c]M2)--([c]M1);
           
         \draw([c]C)--([c]M2); 
         \draw([c]M2)--([c]A); 
        \draw([c]M2)--([c]B); 
        
        \draw([c]A)--([c]M3); 
         \draw([c]M3)--([c]B); 
        \draw([c]M3)--([c]C); 
        
        \draw ([c]B)--([c]M3b); 
         \draw([c]M3b)--([c]A); 
        \draw([c]M3b)--([c]C); 
        
        \node at ([c]C) [above] {\scriptsize{$l_3$}};
        \node at ([c]A) [left] {\scriptsize{$l_1$}};
        \node at ([c]B) [right] {\scriptsize{$l_2$}};

             \draw[thick] ([c]A)--([c]B);
           \fill[red] ([c]M1) circle (2pt);
           \fill[orange] ([c]M2) circle (2pt);
           \draw[thick]([c]M1)--([c]A);
           \draw[thick]([c]M1)--([c]B);
            \draw[blue]([c]M2)--([c]A);
           \draw[]([c]M2)--([c]B);
          \draw[green] ([c]M2)--([c]M1);
         
          \node at ([c]M1) [above right] {\begin{color}{red}\scriptsize{$r_1$}\end{color}};
         \node at ([c]M2) [above right] {\begin{color}{orange}\scriptsize{$r_2$}\end{color}};
          \node at ([c]MM) [left=-.1cm] {\begin{color}{green}\scriptsize{$\tau_2$}\end{color}};
                \node at ([c]MA2) [below] {\begin{color}{blue}\scriptsize{$\tau_3$}\end{color}};
          \node at ([c]CoB) [right] {\scriptsize{$\sigma_1$}};
           \node at ([c]CoA) [below] {\scriptsize{$\sigma_2$}};
      
              \fill[red] ([c]M1) circle (3pt);
           \fill[orange] ([c]M2) circle (3pt);
             \fill ([c]M3) circle (2pt);
              \fill ([c]M3b) circle (2pt);

        \end{scope}
        
        
        \begin{scope} [every coordinate/.style={shift={(8,-2,0)}}]
        
         \node at ([c]label)[below]  {(c)};
        
        \draw ([c]A)--([c]B);
         \draw ([c]A)--([c]C);
          \draw ([c]C)--([c]B);
          
         \draw([c]C)--([c]M2); 
         \draw([c]M2)--([c]A); 
        \draw([c]M2)--([c]B); 
            \draw([c]M2)--([c]A2); 
        
        \draw([c]A)--([c]M3); 
         \draw([c]M3)--([c]B); 
        \draw([c]M3)--([c]C); 
           \draw([c]M3)--([c]B1); 
        
        \draw ([c]B)--([c]M3b); 
         \draw([c]M3b)--([c]A); 
        \draw([c]M3b)--([c]C); 
         \draw([c]M3b)--([c]C1); 
        
        \node at ([c]C) [above] {\scriptsize{$l_3$}};
        \node at ([c]A) [left] {\scriptsize{$l_1$}};
        \node at ([c]B) [right] {\scriptsize{$l_2$}};
        
          \fill ([c]M1) circle (2pt);
           \fill ([c]M2) circle (2pt);
             \fill ([c]M3) circle (2pt);
              \fill ([c]M3b) circle (2pt);

        \end{scope}

         \coordinate(cone1) at (2.8,2.2,0);
        \coordinate(cone) at (3.2,2.5,0);
         \coordinate (o) at (0,0,0);
          \coordinate (t3) at (-1,-3.3,0);
           \coordinate (t2) at (1,-3.3,0);
          \coordinate (t1) at (2,-3,-1);
          \coordinate (t3p) at (-2.5,-3.3,0);
          \coordinate (t3pp) at (-2.5,-3.3,0);
           \coordinate (t2p) at (2.5,-3.3,0);
          \coordinate (t1p) at (3.5,-3,-1);
          \coordinate (t1plong) at (4.5,-3,-1);
          \coordinate (t3halfrip) at (-.5,-2.5,-1);
          \coordinate (t3rip) at (-.6,-3.3,0);
          \coordinate (s3) at (-.5,-1.8,0);
           \coordinate (s3rip) at (-.2,-1.8,0);
          \coordinate (s2) at (.5,-1.8,0);
          \coordinate (s1) at (1,-1.4,-.5);
            \coordinate (s31) at (-.2,-2.5,0);
           \coordinate (s3halfrip) at (-.3,-.7,0);
            
            \coordinate (t2plong) at (3.5,-3.3,0);
            
      \coordinate (y1) at (5.5,0,0);
	   \coordinate (y2) at (5.5,-3.4,0);
	     \coordinate (b1) at (5.2,-.4,0);
	    \coordinate (b2) at (6.8,.2,0);
	     \coordinate (b1e) at (6.2,.3,0);
	    \coordinate (b2e) at (7.8,-.5,0);
	     \coordinate (y1e) at (5.7,0,0);
	   \coordinate (y2e) at (5.7,-3.4,0);
	   \coordinate (c1) at (5.2,-3,0);
	    \coordinate (c2) at (6.8,-2.4,0);
	    \coordinate (d1) at (5.2,-1.6,0);
	    \coordinate (d2) at (6.8,-1,0);
            
            
          \coordinate (x1) at (5,-1.5,0);
	   \coordinate (x2) at (8,-1.5,0);
	   \coordinate (oo) at(6.2,-1.5,0);
	    \coordinate (m1) at (8.4,-1.8,2);
	   \coordinate (m2) at (5.7,-1.8,2);
	   \coordinate (z1) at (5.4,-.2,0);
	   \coordinate (z2) at (7,-3,0);
	   \coordinate (w1) at (7,0.2,0);
	   \coordinate (w2) at (5.4,-3,0);
	      
	        \coordinate (u1) at (6.4,-.1,0);
	   \coordinate (u2) at (8.4,-.5,0);
	   
	    \coordinate (v1) at (6.3,-2.5,0);
	   \coordinate (v2) at (8.5,-3,0);
          
          \begin{scope}[every coordinate/.style={shift={(-18,-7,0)}}]
          
           \node at ([c]cone) [left] {\begin{color}{red}$r_1:$\end{color}}; 
             \node at ([c]cone1) [below] {\scriptsize{$l_1=l_2=l_3$}};

            \draw ([c]o)--([c]t1) node[above]{\scriptsize{$T_2$}};
            \draw ([c]o)--([c]t2) node[below]{\scriptsize{$T_1$}};
            \draw ([c]o)--([c]t3) node[below]{\scriptsize{$T_3$}};
              
              \draw ([c]t1)--([c]t1p);
              \draw ([c]t2)--([c]t2p);
              \draw ([c]t3)--([c]t3p);
             
         \node at ([c]s1) [above] {\begin{color}{blue}\scriptsize{$1$}\end{color}};   
          \node at ([c]s2) [right] {\begin{color}{blue}\scriptsize{$1$}\end{color}}; 
           \node at ([c]s3) [left] {\begin{color}{blue}\scriptsize{$1$}\end{color}};

           \draw[gray] ([c]y1)-- ([c]y1e)--([c]y2e)--([c]y2)--([c]y1);
            \draw[thick] ([c]b1)-- ([c]b2);
            \draw[thick] ([c]c1)-- ([c]c2);
              \draw[thick] ([c]d1)-- ([c]d2);
              
            \node at ([c]b2) [right] {\scriptsize{$T_1$}};
            \node at ([c]c2) [right] {\scriptsize{$T_3$}};   
            \node at ([c]d2) [right] {\scriptsize{$T_2$}};   
            
            \draw[fill=red] ([c]o) circle (3pt);
           \fill ([c]t1) circle (3pt);
             \fill ([c]t2) circle (3pt);
              \fill ([c]t3) circle (3pt);
              
          \end{scope}
          
           \begin{scope}[every coordinate/.style={shift={(-2,-7,0)}}]
          
           \node at ([c]cone) [left] {\begin{color}{orange}$r_2:$\end{color}}; 
             \node at ([c]cone1) [below] {\scriptsize{$l_1=l_2=2l_3$}};

            \draw ([c]o)--([c]t1);
            \draw ([c]o)--([c]t2);
            \draw ([c]o)--([c]t3rip);
              
              \draw ([c]t1)--([c]t1p);
              \draw ([c]t2)--([c]t2p);
              \draw ([c]t3rip)--([c]t3pp);
             
         \node at ([c]s1) [above] {\begin{color}{blue}\scriptsize{$1$}\end{color}};   
          \node at ([c]s2) [right] {\begin{color}{blue}\scriptsize{$1$}\end{color}}; 
           \node at ([c]s3rip) [left] {\begin{color}{blue}\scriptsize{$2$}\end{color}};

             \draw[fill=white] ([c]o) circle (3pt);
           \fill ([c]t1) circle (3pt);
             \fill ([c]t2) circle (3pt);
              \draw[fill=red] ([c]t3rip) circle (3pt);
              
           \end{scope}
           
            \begin{scope}[every coordinate/.style={shift={(0,-7,0)}}]
           
           
            \draw[thick] ([c]x1)--([c]x2);
                \draw[thick] ([c]z1)--([c]z2);
                \draw[thick] ([c]w1)--([c]w2);
              \draw[gray] ([c]oo) parabola ([c]m1);
               \draw[gray] ([c]oo) parabola ([c]m2); 
              
            \node at ([c]w1) [right] {\scriptsize{$T_1$}};
            \node at ([c]x2) [right] {\scriptsize{$T_3$}};   
            \node at ([c]z1) [left] {\scriptsize{$T_2$}};

          \end{scope}
         
          \begin{scope}[every coordinate/.style={shift={(-18,-15,0)}}]
          
           \node at ([c]cone) [left] {\begin{color}{green}$\tau_2:$\end{color}}; 
             \node at ([c]cone1) [below] {\scriptsize{$l_3\leq l_1=l_2\leq 2l_3$}};

            \draw ([c]o)--([c]t1);
            \draw ([c]o)--([c]t2);
            \draw ([c]o)--([c]t3halfrip);
            \draw ([c]t3)--([c]t3halfrip);
              
              \draw ([c]t1)--([c]t1p);
              \draw ([c]t2)--([c]t2p);
              \draw ([c]t3)--([c]t3p);
             
         \node at ([c]s1) [above] {\begin{color}{blue}\scriptsize{$1$}\end{color}};   
          \node at ([c]s2) [right] {\begin{color}{blue}\scriptsize{$1$}\end{color}}; 
           \node at ([c]s3halfrip) [left] {\begin{color}{blue}\scriptsize{$2$}\end{color}}; 
            \node at ([c]s31) [below] {\begin{color}{blue}\scriptsize{$1$}\end{color}};

             \draw [decorate,decoration={brace,amplitude=5pt,mirror},xshift=0.2pt,yshift=-0.4pt] ([c]o)--([c]t3halfrip) node[black,midway,xshift=-0.6cm,yshift=-0.2cm] {\scriptsize{$l_2-l_3$}};
           \draw [decorate,decoration={brace,amplitude=5pt,mirror},xshift=0.2pt,yshift=-0.3pt] ([c]t3halfrip)--([c]t3) node[black,midway,xshift=-0.7cm,yshift=0.1cm] {\scriptsize{$2l_3-l_2$}};

           \draw[gray] ([c]y1)-- ([c]y1e)--([c]y2e)--([c]y2)--([c]y1);
            \draw[thick] ([c]b1)-- ([c]b2);
            \draw[thick] ([c]c1)-- ([c]c2);
              \draw[thick] ([c]d1)-- ([c]d2);
              
            \node at ([c]b2) [right] {\scriptsize{$T_1$}};
            \node at ([c]c2) [right] {\scriptsize{$T_3$}};   
            \node at ([c]d2) [right] {\scriptsize{$T_2$}};   
           
             \draw[fill=white] ([c]o) circle (3pt);
           \fill ([c]t1) circle (3pt);
             \fill ([c]t2) circle (3pt);
              \fill ([c]t3) circle (3pt);
            
             \draw ([c]t3halfrip) node[red,cross]{};

          \end{scope}

        \begin{scope}[every coordinate/.style={shift={(-2,-15,0)}}]
          
           \node at ([c]cone) [left] {\begin{color}{blue}$\tau_3:$\end{color}}; 
             \node at ([c]cone1) [below] {\scriptsize{$l_1\geq l_2=2l_3$}};

            \draw ([c]o)--([c]t1);
            \draw ([c]o)--([c]t2);
            \draw ([c]o)--([c]t3rip);
              
              \draw ([c]t1)--([c]t1p);
              \draw ([c]t2)--([c]t2p);
               \draw ([c]t2p)--([c]t2plong);
              \draw ([c]t3rip)--([c]t3pp);
             
         \node at ([c]s1) [above] {\begin{color}{blue}\scriptsize{$1$}\end{color}};   
          \node at ([c]s2) [right] {\begin{color}{blue}\scriptsize{$1$}\end{color}}; 
           \node at ([c]s3rip) [left] {\begin{color}{blue}\scriptsize{$2$}\end{color}};

            \draw [decorate,decoration={brace,amplitude=5pt,mirror},xshift=0.2pt,yshift=-0.4pt] ([c]o)--([c]t2) node[black,midway,xshift=-0.1cm,yshift=-0.4cm] {\scriptsize{$l_2$}};
           \draw [decorate,decoration={brace,amplitude=5pt,mirror},xshift=0.2pt,yshift=-0.3pt] ([c]t2)--([c]t2p) node[black,midway,xshift=0cm,yshift=-0.3cm] {\scriptsize{$l_1-l_2$}};

             \draw[fill=white] ([c]o) circle (3pt);
           \fill ([c]t1) circle (3pt);
             \fill ([c]t2p) circle (3pt);
              \draw[fill=red] ([c]t3rip) circle (3pt);
              
               \draw ([c]t2) node[cross]{};
           \end{scope}
           \begin{scope}[every coordinate/.style={shift={(0,-15,0)}}]
           
           
            \draw[thick] ([c]x1)--([c]x2);
                \draw[thick] ([c]z1)--([c]z2);
                 \draw[thick] ([c]u1)--([c]u2);
                \draw[gray] ([c]w1)--([c]w2);
              \draw[gray] ([c]oo) parabola ([c]m1);
               \draw[gray] ([c]oo) parabola ([c]m2); 
              
            \node at ([c]w1) [right] {\scriptsize{$T_1$}};
            \node at ([c]x2) [right] {\scriptsize{$T_3$}};   
            \node at ([c]z1) [left] {\scriptsize{$T_2$}};

          \end{scope}
         
        \begin{scope}[every coordinate/.style={shift={(-18,-23,0)}}]

           \node at ([c]cone) [left] {$\sigma_1$:}; 
             \node at ([c]cone1) [below] {\scriptsize{$l_3\leq l_2\leq l_1, 2l_3$}};

            \draw ([c]o)--([c]t1);
            \draw ([c]o)--([c]t2);
            \draw ([c]o)--([c]t3halfrip);
            \draw ([c]t3)--([c]t3halfrip);
              
              \draw ([c]t1)--([c]t1p);
              \draw ([c]t2)--([c]t2p);
               \draw ([c]t2p)--([c]t2plong);
              \draw ([c]t3)--([c]t3p);
             
         \node at ([c]s1) [above] {\begin{color}{blue}\scriptsize{$1$}\end{color}};   
          \node at ([c]s2) [right] {\begin{color}{blue}\scriptsize{$1$}\end{color}}; 
           \node at ([c]s3halfrip) [left] {\begin{color}{blue}\scriptsize{$2$}\end{color}}; 
            \node at ([c]s31) [below] {\begin{color}{blue}\scriptsize{$1$}\end{color}};

             \draw [decorate,decoration={brace,amplitude=5pt,mirror},xshift=0.2pt,yshift=-0.4pt] ([c]o)--([c]t3halfrip) node[black,midway,xshift=-0.6cm,yshift=-0.2cm] {\scriptsize{$l_2-l_3$}};
           \draw [decorate,decoration={brace,amplitude=5pt,mirror},xshift=0.2pt,yshift=-0.3pt] ([c]t3halfrip)--([c]t3) node[black,midway,xshift=-0.7cm,yshift=0.1cm] {\scriptsize{$2l_3-l_2$}};

            \draw [decorate,decoration={brace,amplitude=5pt,mirror},xshift=0.2pt,yshift=-0.4pt] ([c]o)--([c]t2) node[black,midway,xshift=-0.1cm,yshift=-0.4cm] {\scriptsize{$l_2$}};
           \draw [decorate,decoration={brace,amplitude=5pt,mirror},xshift=0.2pt,yshift=-0.3pt] ([c]t2)--([c]t2p) node[black,midway,xshift=0cm,yshift=-0.3cm] {\scriptsize{$l_1-l_2$}};
            
           \draw[fill=white] ([c]o) circle (3pt);
           \fill ([c]t1) circle (3pt);
             \fill ([c]t2p) circle (3pt);
              \fill ([c]t3) circle (3pt);
            
             \draw ([c]t3halfrip) node[red,cross]{};
               \draw ([c]t2) node[cross]{};

           \draw[gray] ([c]y1)-- ([c]y1e)--([c]y2e)--([c]y2)--([c]y1);
            \draw[thick] ([c]b1e)-- ([c]b2e);
                 \draw[gray] ([c]b1)-- ([c]b2);
            \draw[thick] ([c]c1)-- ([c]c2);
              \draw[thick] ([c]d1)-- ([c]d2);
              
            \node at ([c]b2e) [right] {\scriptsize{$T_1$}};
            \node at ([c]c2) [right] {\scriptsize{$T_3$}};   
            \node at ([c]d2) [right] {\scriptsize{$T_2$}};

          \end{scope}

        \begin{scope}[every coordinate/.style={shift={(-2,-23,0)}}]
          
           \node at ([c]cone) [left] {$\sigma_2$:}; 
             \node at ([c]cone1) [below] {\scriptsize{$l_1,l_2\geq 2l_3$}};

            \draw ([c]o)--([c]t1);
            \draw ([c]o)--([c]t2);
            \draw ([c]o)--([c]t3rip);
              
              \draw ([c]t1)--([c]t1p);
              \draw ([c]t2)--([c]t2p);
               \draw ([c]t2p)--([c]t2plong);
                \draw ([c]t1p)--([c]t1plong);
              \draw ([c]t3rip)--([c]t3pp);
             
         \node at ([c]s1) [above] {\begin{color}{blue}\scriptsize{$1$}\end{color}};   
          \node at ([c]s2) [right] {\begin{color}{blue}\scriptsize{$1$}\end{color}}; 
           \node at ([c]s3rip) [left] {\begin{color}{blue}\scriptsize{$2$}\end{color}};

            \draw [decorate,decoration={brace,amplitude=5pt,mirror},xshift=0.2pt,yshift=-0.4pt] ([c]o)--([c]t2) node[black,midway,xshift=-0.15cm,yshift=-0.4cm] {\scriptsize{$2l_3$}};
           \draw [decorate,decoration={brace,amplitude=5pt,mirror},xshift=0.2pt,yshift=-0.3pt] ([c]t2)--([c]t2p) node[black,midway,xshift=0cm,yshift=-0.3cm] {\scriptsize{$l_1-2l_3$}};
           
           \draw [decorate,decoration={brace,amplitude=5pt,mirror},xshift=0.2pt,yshift=-0.4pt] ([c]t1)--([c]o) node[black,midway,xshift=0.4cm,yshift=0.2cm] {\scriptsize{$2l_3$}};
           \draw [decorate,decoration={brace,amplitude=5pt,mirror},xshift=0.2pt,yshift=-0.3pt] ([c]t1p)--([c]t1) node[black,midway,xshift=0.2cm,yshift=0.4cm] {\scriptsize{$l_2-2l_3$}};
         
           \draw[fill=white] ([c]o) circle (3pt);
           
             \fill ([c]t2p) circle (3pt);
              \fill ([c]t1p) circle (3pt);
              \draw[fill=red] ([c]t3rip) circle (3pt);
              
               \draw ([c]t2) node[cross]{};
            \draw ([c]t1) node[cross]{};
            
            \end{scope}
             \begin{scope}[every coordinate/.style={shift={(0,-23,0)}}]
           
            \draw[thick] ([c]x1)--([c]x2);
                \draw[gray] ([c]z1)--([c]z2);
                 \draw[thick] ([c]u1)--([c]u2);
                   \draw[thick] ([c]v1)--([c]v2);
                \draw[gray] ([c]w1)--([c]w2);
              \draw[gray] ([c]oo) parabola ([c]m1);
               \draw[gray] ([c]oo) parabola ([c]m2); 
              
            \node at ([c]w1) [right] {\scriptsize{$T_1$}};
            \node at ([c]x2) [right] {\scriptsize{$T_3$}};   
            \node at ([c]v2) [right] {\scriptsize{$T_2$}};   
              
          \end{scope}
       
       \end{tikzpicture}
	\caption{(a) Subdivision $\Sigma$ (non-simplicial). (b) Subdivision $\Sigma'$  for high-degree tails. (c) Subdivision $\Sigma'$  for degree-two tails.}
 \label{esesub2}
\end{figure}
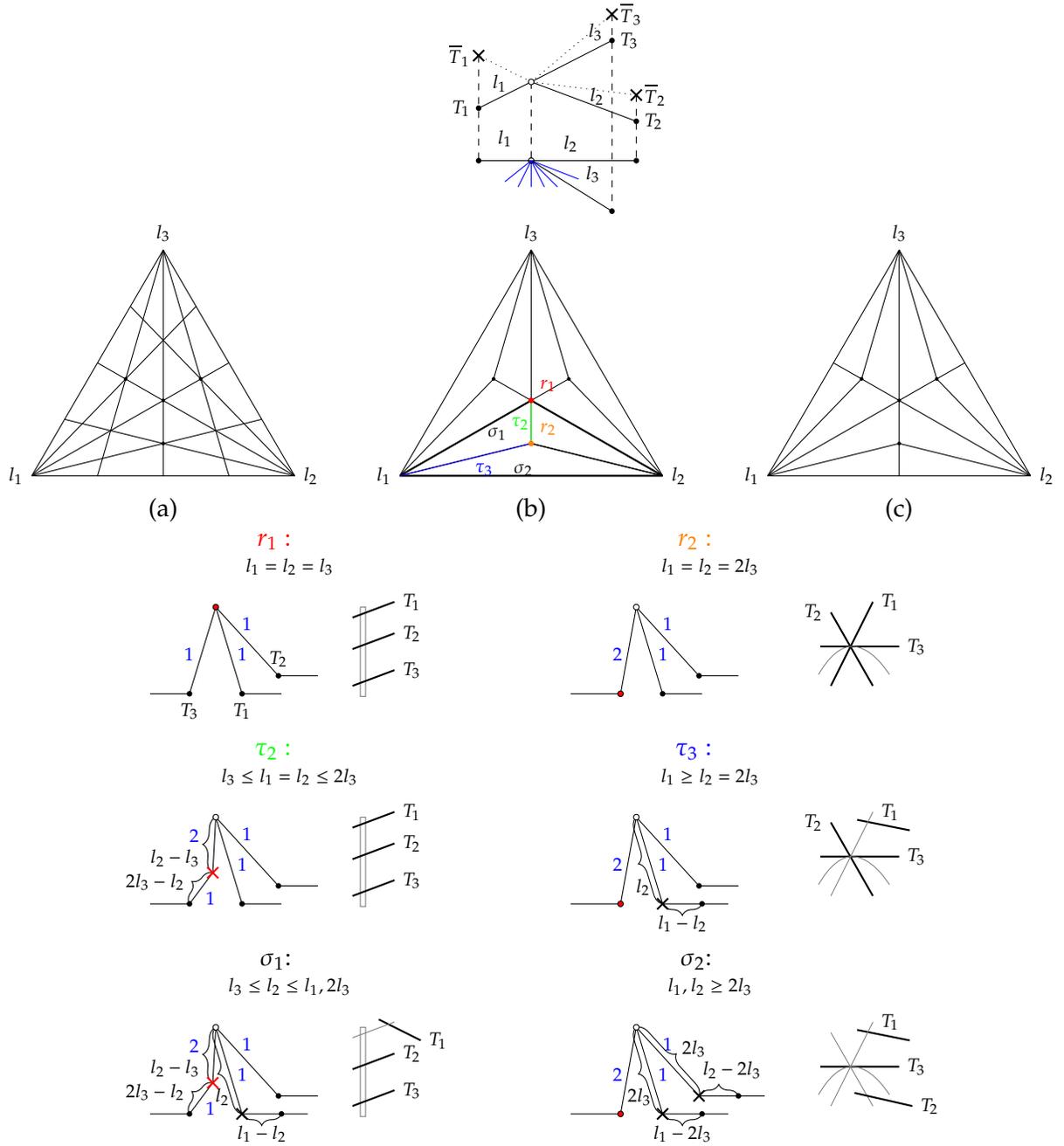
\end{exa}

\begin{exa}\label{exa:reduciblecore}
The core of $C$ consists of two elliptic curves of weight zero, meeting in a node; each elliptic curve is attached to a high-degree rational tail. Notice that the central cone $\sigma_2$ is not smooth; again, we need a Kummer extension (half-lengths).
We also remark that the hyperplane $l_1+l_2=m$ does not come from the alignment, but from the procedure of Definition \ref{def:semistablereduction}. See Figure \ref{esesub3}.
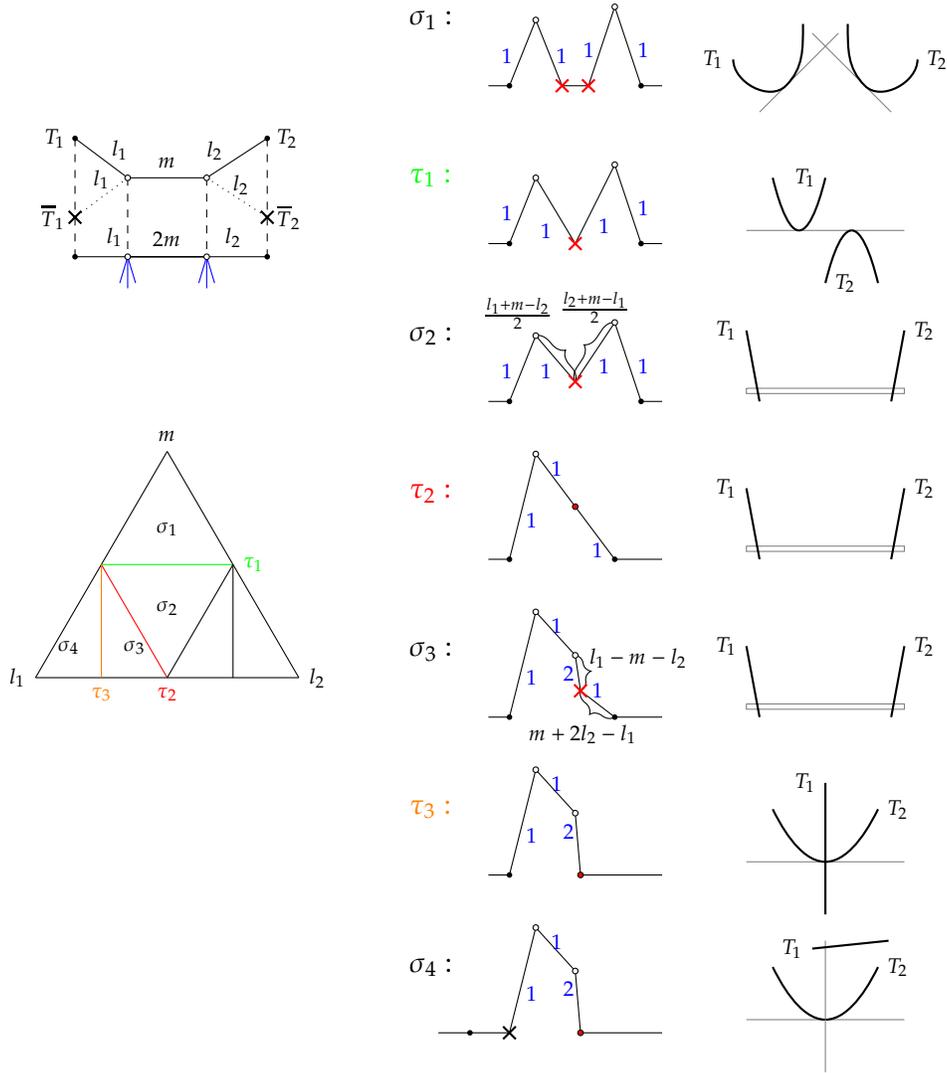
\begin{figure}[h]
\begin{tikzpicture} [scale=.35]
\tikzstyle{every node}=[font=\normalsize]
		\tikzset{arrow/.style={latex-latex}}
		
	\tikzset{cross/.style={cross out, draw, thick,
         minimum size=2*(#1-\pgflinewidth), 
         inner sep=1.2pt, outer sep=1.2pt}}


		\coordinate (E1) at (-1.5,3,0);
		\coordinate (E2) at (1.5,3,0);
		
		\coordinate (T1) at (-3.5,4.5,0);
		\coordinate (bT1) at (-3.5,1.5,0);
		
		\coordinate (T2) at (3.8,4.5,0);
		\coordinate (bT2) at (3.8,1.5,0);
		
		\coordinate (E1t) at (-1.5,0,0);
		\coordinate (E2t) at (1.5,0,0);
		\coordinate (T1t) at (-3.5,0,0);
		\coordinate (T2t) at (3.8,0,0);
		
		\coordinate (P1) at (-1.8,-1,0);
		\coordinate (P2) at (-1.5,-1.2,0);
		\coordinate (P3) at (-1.2,-1,0);
		\coordinate (P4) at (1.8,-1,0);
		\coordinate (P5) at (1.5,-1.2,0);
    	\coordinate (P6) at (1.2,-1,0);
		
		
		\coordinate (l1) at (-1.8,3.3,0);
		\coordinate (bl1) at (-2.5,2.2,0);
		\coordinate (l2) at (1.8,3.3,0);
		\coordinate (bl2) at (2.8,2,0);
		\coordinate (m) at (0,3,0);
		
		\coordinate (l1t) at (-2,0.,0);
		\coordinate (l2t) at (2.5,0.,0);
		\coordinate (mt) at (0,0,0);
		
		\begin{scope}[every coordinate/.style={shift={(-6,6,0)}}]

     	
     	\draw ([c]E1)-- ([c]T1);
     	\draw[dotted] ([c]E1)-- ([c]bT1);
     	\draw ([c]E2)-- ([c]T2);
     	\draw[dotted] ([c]E2)-- ([c]bT2);
     	\draw ([c]E1)-- ([c]E2);

     	\draw ([c]E1t)-- ([c]E2t);
     	\draw ([c]E1t)-- ([c]T2t);
     	\draw ([c]E2t)-- ([c]T1t);
     	
     	\draw[blue] ([c]E1t)--([c]P1);
     	\draw[blue] ([c]E1t)--([c]P2);
     	\draw[blue] ([c]E1t)--([c]P3);
     	\draw[blue] ([c]E2t)--([c]P4);
     	\draw[blue] ([c]E2t)--([c]P5);
     	\draw[blue] ([c]E2t)--([c]P6);

       
       \node at ([c]l1) [above] {\scriptsize{$l_1$}};
        \node at ([c]bl1) [above] {\scriptsize{$l_1$}};
        \node at ([c]l2) [above] {\scriptsize{$l_2$}};
         \node at ([c]bl2) [above] {\scriptsize{$l_2$}};
           \node at ([c]m) [above] {\scriptsize{$m$}};

      \node at ([c]l1t) [above] {\scriptsize{$l_1$}};
       \node at ([c]l2t) [above] {\scriptsize{$l_2$}};
      \node at ([c]mt) [above] {\scriptsize{$2m$}};
        
      \node at ([c]bT1) [left] {\scriptsize{$\overline{T}_1$}};
        \node at ([c]T1) [left] {\scriptsize{$T_1$}};
         
        \node at ([c]bT2) [right] {\scriptsize{$\overline{T}_2$}};
        \node at ([c]T2) [right] {\scriptsize{$T_2$}};

       \draw[dashed] ([c]E1)--([c]E1t);
       \draw[dashed] ([c]T1)--([c]T1t);
       \draw[dashed] ([c]T2)--([c]T2t);
       \draw[dashed] ([c]E2)--([c]E2t);
      
		\draw[fill=white]  ([c]E1) circle (3pt);
			\draw[fill=white]  ([c]E2) circle (3pt);
			\draw[fill=white]   ([c]E1t) circle (3pt);
			\draw[fill=white]   ([c]E2t) circle (3pt);
			\draw  ([c]bT1) node [cross] {};
			\draw  ([c]bT2) node [cross] {};
			
			\foreach \x in {T1,T2,T1t,T2t} 
          \fill ([c]\x) circle (3pt);
      
         \end{scope}
        
      
         \coordinate (A) at (-5,0,0);
         \coordinate (A1) at (-2.5,0,0);
         \coordinate (A2) at (0,0,0);
        \coordinate (A3) at (2.5,0,0);
        \coordinate (B) at (5,0,0);
         \coordinate (B1) at (2.5,4.3,0);
         \coordinate (C) at (0,8.6,0);
         \coordinate (C1) at (-2.5,4.3,0);
         
        
         \coordinate (s1) at (0,5,0);
          \coordinate (s2) at (0,2,0);
           \coordinate (s3) at (-1.25,0.5,0);
            \coordinate (s4) at (-3.75,0.5,0);

         \begin{scope} [every coordinate/.style={shift={(-6,-10,0)}}]

        \draw ([c]A)--([c]B);
         \draw ([c]A)--([c]C);
          \draw ([c]C)--([c]B);
          
         \draw[green]([c]C1)--([c]B1); 
         \draw[red]([c]C1)--([c]A2); 
        \draw ([c]B1)--([c]A2); 
        
        \draw[orange]([c]A1)--([c]C1); 
         \draw([c]A3)--([c]B1);

        \node at ([c]C) [above] {\scriptsize{$m$}};
        \node at ([c]A) [left] {\scriptsize{$l_1$}};
        \node at ([c]B) [right] {\scriptsize{$l_2$}};
        
        \node at ([c]B1) [right] {\begin{color}{green}\scriptsize{$\tau_1$}\end{color}};
         \node at ([c]A1) [below] {\begin{color}{orange}\scriptsize{$\tau_3$}\end{color}};
          \node at ([c]A2) [below] {\begin{color}{red}\scriptsize{$\tau_2$}\end{color}};

          \node at ([c]s1) [above] {\scriptsize{$\sigma_1$}};
          \node at ([c]s2) [above] {\scriptsize{$\sigma_2$}};
          \node at ([c]s3) [above] {\scriptsize{$\sigma_3$}};
          \node at ([c]s4) [above] {\scriptsize{$\sigma_4$}};
        
         \end{scope}
      
      
       \coordinate (label)  at (-4,1,0);
      \coordinate (o1) at  (0,1,0);
       \coordinate (t1m) at  (-1,-1.5,0);
       \coordinate (t1c) at  (-1.8,-1.5,0);
        \coordinate (t1p) at  (1,-1.5,0);
        \coordinate (t2m) at  (2,-1.5,0);
         \coordinate (o2) at  (3,1.5,0);
         \coordinate (t2p) at  (4,-1.5,0);
         \coordinate (t2c) at  (4.8,-1.5,0);
        
         \coordinate (mm) at (1.5,-1.5,0);
         
           \coordinate (nn) at (1.5,-.75,0);
           \coordinate (mn) at (1.5,1.25,0);
        
          \coordinate (o1bis) at  (0,2.5,0);
            \coordinate (o2bis) at  (1.5,.5,0);
            \coordinate (t2bis) at  (3,-1.5,0);
          
              \coordinate (t2tris) at  (1.7,-.5,0);
               \coordinate (o2tris) at  (1.5,.85,0);
               
                
                \coordinate (t2rip) at  (1.7,-1.5,0);
                
                 \coordinate (t1mbis) at  (-2.5,-1.5,0);
                  \coordinate (t1mtris) at  (-3.7,-1.5,0);
         
         \coordinate (s1) at (-.5,-.35,0);
          \coordinate (s2) at (.4,-.35,0);
          \coordinate (s3) at (2.6,-.15,0);
          \coordinate (s4) at (3.5,-.15,0);
        
          
          \coordinate (ss1) at (-.75,0,0);
          \coordinate (ss2) at (.8,1.3,0);
          \coordinate (ss3) at  (2.35,-1.15,0);
           \coordinate (ss4) at  (1.25,.-.5,0);
        
            \coordinate (x1) at  (-2.5,-2.5,0);
         \coordinate (x2) at  (.5,.5,0);
         \coordinate(oo1) at (-1.5,-1.5,0);
         \coordinate (m1) at  (-.85,.85,0);
         \coordinate (m2) at  (-3.5,-.5,0);
          \coordinate(oo2) at (1.5,-1.5,0);
         \coordinate (n1) at  (.85,.85,0);
         \coordinate (n2) at  (3.5,-.5,0);
         \coordinate (y1) at  (2.5,-2.5,0);
         \coordinate (y2) at  (-.5,.5,0);
         
           \coordinate (z1) at (-3,-1,0);
            \coordinate (z2) at (3,-1,0);
            \coordinate (l1) at (-2,1,0);
            \coordinate (ol) at (-1,-1,0);
            \coordinate (l2) at (0,1,0);
            \coordinate (j1) at (2,-3,0);
            \coordinate (oj) at (1,-1,0);
            \coordinate (j2) at (0,-3,0);
            
            \coordinate (z1e) at (-3,-1.2,0);
            \coordinate (z2e) at (3,-1.2,0);
             \coordinate (x1e) at (-2.5,-1.5,0);
            \coordinate (x2e) at (-3,1.2,0);
             \coordinate (y1e) at (2.5,-1.5,0);
            \coordinate (y2e) at (3,1.2,0);
         
            \coordinate (g1) at (-2,1,0);
            \coordinate (on) at (0,-1,0);
            \coordinate (g2) at (2,1,0);
             \coordinate (v1) at (0,2,0);
            \coordinate (v2) at (0,-3,0);
            
                  \coordinate (w1) at (-.5,1.7,0);
            \coordinate (w2) at (2.4,2,0);
         
           \begin{scope} [every coordinate/.style={shift={(8,14,0)}}]
           \node at ([c]label) {$\sigma_1:$};
          
            \draw ([c]o1)--([c]t1m);
             \draw ([c]o1)--([c]t1p);
            \draw ([c]o2)--([c]t2p);
            \draw ([c]o2)--([c]t2m);
              \draw ([c]t1p)--([c]t2m);
              \draw ([c]t1c)--([c]t1m);
              \draw ([c]t2c)--([c]t2p);
               
               
               \node at ([c]s1)  [left] {\begin{color}{blue}\scriptsize{$1$}\end{color}};
             \node at ([c]s2) [right]  {\begin{color}{blue}\scriptsize{$1$}\end{color}};
             \node at ([c]s3)  [left] {\begin{color}{blue}\scriptsize{$1$}\end{color}};
             \node at ([c]s4) [right]  {\begin{color}{blue}\scriptsize{$1$}\end{color}};
               
            \draw[fill=white] ([c]o1) circle (3pt);
           \draw[fill=white] ([c]o2) circle (3pt);
             \fill ([c]t2p) circle (3pt);
              \fill ([c]t1m) circle (3pt);

              \draw ([c]t2m) node[cross,red]{};
            \draw ([c]t1p) node[cross,red]{};
            
           \end{scope}
            \begin{scope} [every coordinate/.style={shift={(19,14,0)}}]
            \draw[gray] ([c]x1)--([c]x2);
                \draw[gray] ([c]y1)--([c]y2);
                \draw[thick] ([c]oo1) to[out=45,in=270] ([c]m1);
               \draw[thick] ([c]oo1) to[out=225,in=270] ([c]m2); 
                \draw[thick] ([c]oo2) to[out=135,in=270] ([c]n1);
               \draw[thick] ([c]oo2) to[out=-45,in=270] ([c]n2);
               
                          \node at ([c]m2)  [left] {\scriptsize{$T_1$}};;
                        \node at ([c]n2) [right]   {\scriptsize{$T_2$}};
        \end{scope}

         \begin{scope} [every coordinate/.style={shift={(8,8,0)}}]
           \node at ([c]label) {\begin{color}{green}$\tau_1:$\end{color}};

            \draw ([c]o1)--([c]mm);
             \draw ([c]o1)--([c]t1m);
            \draw ([c]o2)--([c]t2p);
            \draw ([c]o2)--([c]mm);
             \draw ([c]t1c)--([c]t1m);
              \draw ([c]t2c)--([c]t2p);
               
               
               \node at ([c]s1)  [left] {\begin{color}{blue}\scriptsize{$1$}\end{color}};
             \node at ([c]s2) [below]  {\begin{color}{blue}\scriptsize{$1$}\end{color}};
             \node at ([c]s3)  [below] {\begin{color}{blue}\scriptsize{$1$}\end{color}};
             \node at ([c]s4) [right]  {\begin{color}{blue}\scriptsize{$1$}\end{color}};
               
                  \draw[fill=white] ([c]o1) circle (3pt);
           \draw[fill=white] ([c]o2) circle (3pt);
             \fill ([c]t2p) circle (3pt);
              \fill ([c]t1m) circle (3pt);

              \draw ([c]mm) node[cross,red]{};

           \end{scope}

         \begin{scope} [every coordinate/.style={shift={(19,8,0)}}]
            \draw[gray] ([c]z1)--([c]z2);
              
                \draw[thick] ([c]ol) parabola ([c]l1);
               \draw[thick] ([c]ol) parabola ([c]l2); 
                \draw[thick] ([c]oj) parabola ([c]j1);
               \draw[thick] ([c]oj) parabola ([c]j2);
               
                          \node at ([c]l2)  [left] {\scriptsize{$T_1$}};
                        \node at ([c]j2) [right]   {\scriptsize{$T_2$}};
        \end{scope}
        
         \begin{scope} [every coordinate/.style={shift={(8,2,0)}}]
           \node at ([c]label) {$\sigma_2:$};

            \draw ([c]o1)--([c]nn);
             \draw ([c]o1)--([c]t1m);
            \draw ([c]o2)--([c]t2p);
            \draw ([c]o2)--([c]nn);
             \draw ([c]t1c)--([c]t1m);
              \draw ([c]t2c)--([c]t2p);
               
               
               \node at ([c]s1)  [left] {\begin{color}{blue}\scriptsize{$1$}\end{color}};
             \node at ([c]s2)  {\begin{color}{blue}\scriptsize{$1$}\end{color}};
             \node at ([c]s3) {\begin{color}{blue}\scriptsize{$1$}\end{color}};
             \node at ([c]s4) [right]  {\begin{color}{blue}\scriptsize{$1$}\end{color}};
             
              \draw [decorate,decoration={brace,amplitude=5pt,mirror},xshift=0.2pt,yshift=-0.4pt] ([c]nn)--([c]o1) node[black,midway,xshift=-0.5cm,yshift=0.6cm] {\scriptsize{$\frac{l_1+m-l_2}{2}$}};
           \draw [decorate,decoration={brace,amplitude=5pt,mirror},xshift=0.2pt,yshift=-0.3pt] ([c]o2)--([c]nn) node[black,midway,yshift=0.6cm] {\scriptsize{$\frac{l_2+m-l_1}{2}$}};
               
               \draw[fill=white] ([c]o1) circle (3pt);
           \draw[fill=white] ([c]o2) circle (3pt);
             \fill ([c]t2p) circle (3pt);
              \fill ([c]t1m) circle (3pt);

              \draw ([c]nn) node[cross,red]{};
               
           \end{scope}

         \begin{scope} [every coordinate/.style={shift={(19,2,0)}}]
            \draw[gray] ([c]z1)--([c]z2)--([c]z2e)--([c]z1e)--([c]z1);
             \draw[thick] ([c]x1e)--([c]x2e); 
              \draw[thick] ([c]y1e)--([c]y2e);  
               
                          \node at ([c]x2e)  [left] {\scriptsize{$T_1$}};
                        \node at ([c]y2e) [right]   {\scriptsize{$T_2$}};
        \end{scope}

        \begin{scope} [every coordinate/.style={shift={(8,-4,0)}}]
           \node at ([c]label) {\begin{color}{red}$\tau_2:$\end{color}};

            \draw ([c]o1bis)--([c]o2bis);
             \draw ([c]o1bis)--([c]t1m);
            \draw ([c]o2bis)--([c]t2bis);
            \draw ([c]t2bis)--([c]t2c);
             \draw ([c]t1c)--([c]t1m);

               
               \node at ([c]ss1)  [right] {\begin{color}{blue}\scriptsize{$1$}\end{color}};
             \node at ([c]ss2) [above]  {\begin{color}{blue}\scriptsize{$1$}\end{color}};
             \node at ([c]ss3)   {\begin{color}{blue}\scriptsize{$1$}\end{color}};
            
               \draw[fill=red] ([c]o2bis) circle (3pt);
           \draw[fill=white] ([c]o1bis) circle (3pt);
             \fill ([c]t2bis) circle (3pt);
              \fill ([c]t1m) circle (3pt);

           \end{scope}

        \begin{scope} [every coordinate/.style={shift={(19,-4,0)}}]
            \draw[gray] ([c]z1)--([c]z2)--([c]z2e)--([c]z1e)--([c]z1);
             \draw[thick] ([c]x1e)--([c]x2e); 
              \draw[thick] ([c]y1e)--([c]y2e);  
               
                          \node at ([c]x2e)  [left] {\scriptsize{$T_1$}};
                        \node at ([c]y2e) [right]   {\scriptsize{$T_2$}};
        \end{scope}
        
         \begin{scope} [every coordinate/.style={shift={(8,-10,0)}}]
           \node at ([c]label) {$\sigma_3:$};

            \draw ([c]o1bis)--([c]o2tris);
             \draw ([c]o1bis)--([c]t1m);
            \draw ([c]o2tris)--([c]t2tris);
             \draw ([c]t2bis)--([c]t2tris);
            \draw ([c]t2bis)--([c]t2c);
             \draw ([c]t1c)--([c]t1m);

               
               \node at ([c]ss1)  [right] {\begin{color}{blue}\scriptsize{$1$}\end{color}};
             \node at ([c]ss2) [above]  {\begin{color}{blue}\scriptsize{$1$}\end{color}};
             \node at ([c]ss3)  [above] {\begin{color}{blue}\scriptsize{$1$}\end{color}};
             \node at ([c]ss4)  [above] {\begin{color}{blue}\scriptsize{$2$}\end{color}};

               \draw [decorate,decoration={brace,amplitude=5pt,mirror},xshift=0.2pt,yshift=-0.2pt] ([c]t2tris)--([c]o2tris) node[black,midway,xshift=0.8cm,yshift=0.2cm] {\scriptsize{$l_1-m-l_2$}};
           \draw [decorate,decoration={brace,amplitude=5pt,mirror},xshift=0.2pt,yshift=-0.5pt] ([c]t2tris)--([c]t2bis) node[black,midway,xshift=-0.2cm,yshift=-0.4cm] {\scriptsize{$m+2l_2-l_1$}};

              \draw[fill=white]([c]o1bis) circle (3pt);
           \draw[fill=white] ([c]o2tris) circle (3pt);
             \fill ([c]t2bis) circle (3pt);
              \fill ([c]t1m) circle (3pt);
              
               \draw ([c]t2tris) node[cross,red]{};

           \end{scope}

        \begin{scope} [every coordinate/.style={shift={(19,-10,0)}}]
            \draw[gray] ([c]z1)--([c]z2)--([c]z2e)--([c]z1e)--([c]z1);
             \draw[thick] ([c]x1e)--([c]x2e); 
              \draw[thick] ([c]y1e)--([c]y2e);  
               
                          \node at ([c]x2e)  [left] {\scriptsize{$T_1$}};
                        \node at ([c]y2e) [right]   {\scriptsize{$T_2$}};
        \end{scope}
 
  \begin{scope} [every coordinate/.style={shift={(8,-16,0)}}]
           \node at ([c]label) {\begin{color}{orange}$\tau_3:$\end{color}};

            \draw ([c]o1bis)--([c]o2tris);
             \draw ([c]o1bis)--([c]t1m);
            \draw ([c]o2tris)--([c]t2rip);
           
            \draw ([c]t2rip)--([c]t2c);
             \draw ([c]t1c)--([c]t1m);

               
               \node at ([c]ss1)  [right] {\begin{color}{blue}\scriptsize{$1$}\end{color}};
             \node at ([c]ss2) [above]  {\begin{color}{blue}\scriptsize{$1$}\end{color}};
             
             \node at ([c]ss4)  [above] {\begin{color}{blue}\scriptsize{$2$}\end{color}};
            
             \draw[fill=white]([c]o1bis) circle (3pt);
           \draw[fill=white] ([c]o2tris) circle (3pt);
             \draw[fill=red] ([c]t2rip) circle (3pt);
              \fill ([c]t1m) circle (3pt);
              
               \end{scope}

               \begin{scope} [every coordinate/.style={shift={(19,-16,0)}}]
            \draw[gray] ([c]z1)--([c]z2);
              
                \draw[thick] ([c]on) parabola ([c]g1);
               \draw[thick] ([c]on) parabola ([c]g2); 
                \draw[thick] ([c]v1) -- ([c]v2);

                          \node at ([c]v1)  [left] {\scriptsize{$T_1$}};
                        \node at ([c]g2) [right]   {\scriptsize{$T_2$}};
        \end{scope}
        
        \begin{scope} [every coordinate/.style={shift={(8,-22,0)}}]
           \node at ([c]label) {$\sigma_4:$};

            \draw ([c]o1bis)--([c]o2tris);
             \draw ([c]o1bis)--([c]t1m);
            \draw ([c]o2tris)--([c]t2rip);
           
            \draw ([c]t2rip)--([c]t2c);
             \draw ([c]t1mtris)--([c]t1mbis);
             \draw ([c]t1mbis)--([c]t1m);
               
               
               \node at ([c]ss1)  [right] {\begin{color}{blue}\scriptsize{$1$}\end{color}};
             \node at ([c]ss2) [above]  {\begin{color}{blue}\scriptsize{$1$}\end{color}};
             
             \node at ([c]ss4)  [above] {\begin{color}{blue}\scriptsize{$2$}\end{color}};
            
               \draw[fill=white]([c]o1bis) circle (3pt);
           \draw[fill=white] ([c]o2tris) circle (3pt);
             \draw[fill=red] ([c]t2rip) circle (3pt);
              \fill ([c]t1mbis) circle (3pt);
              
                 \draw ([c]t1m) node[cross]{};
                 
               \end{scope}

               \begin{scope} [every coordinate/.style={shift={(19,-22,0)}}]
            \draw[gray] ([c]z1)--([c]z2);
              
                \draw[thick] ([c]on) parabola ([c]g1);
               \draw[thick] ([c]on) parabola ([c]g2); 
                \draw[gray] ([c]v1) -- ([c]v2);
                \draw[thick] ([c]w1) -- ([c]w2);

                          \node at ([c]w1)  [left] {\scriptsize{$T_1$}};
                        \node at ([c]g2) [right]   {\scriptsize{$T_2$}};
        \end{scope}

  \end{tikzpicture}
	\caption{Tropical admissible cover (top left); the subdivision of $\operatorname{trop}(\mathcal A_2^{\text{wt}})$ (bottom left);  $\lambda_T$  and the  Gorenstein singularities (right).}
 \label{esesub3}
\end{figure}
\end{exa}

\newpage

\subsection{The primary construction}\label{sec:main}
We construct a universal morphism $\widetilde{\mathcal C}\to \overline{\mathcal C}$ over $\TAt$, where $\overline{\mathcal C}$ is a family of Gorenstein singularities (both isolated and ribbons) with core of positive weight. We do so in two steps: first, a contraction informed by $\lambda$, producing a possibly non-Gorenstein curve. The image in $\TAt$ of the non-Gorenstein locus is contained in a divisor, which we name $\Do$ below. We complete the construction of $\overline{\mathcal C}$ by gluing in a portion of $\psi$ over this locus, thus producing a non-reduced structure along the fibres.

\begin{dfn}
 Let $\rho_\text{max}$ denote the maximum value attained by $\lambda$ on $\tropC$, and let $\pazocal D\subseteq\TAt$ be the Cartier divisor determined by the ideal sheaf $\OO_{\widetilde{\mathcal A}}(-\rho_\text{max})\hookrightarrow\OO_{\widetilde{\mathcal A}}$.
\end{dfn} 
\begin{dfn}
Recall from Definition \ref{def:rho1} that $\rho_1$ denotes the value attained by $\lambda$ at the vertex(ices) supporting $K_{\tropD}+\operatorname{div}(\lambda)$. Let $\Do\subseteq\TAt$ be the Cartier divisor cut by the ideal sheaf $\OO_{\widetilde{\mathcal A}}(-\rho_1)\hookrightarrow\OO_{\widetilde{\mathcal A}}$. 
\end{dfn} Note that $\rho_1\leq\rho_\text{max}$ implies $\Do\subseteq \pazocal D$. 
By Lemma \ref{lem:ilcorestasopraD}, the subcurve $\tropC_{\geq\rho_1}$ always contains the core, so in particular $p_a(\tropD^\circ)=2$ over $\Do$. Ribbons will appear in $\overline{\mathcal C}$ precisely over $\Do$.
\begin{dfn}
 Let $\Dt$ be cut by $\rho_\text{max}-\rho_1$, so that we have an exact sequence:
 \[0\to\OO_{\Dt}(-\rho_1)\to\OO_{\pazocal D}\to\OO_{\Do}\to0\]
 \end{dfn}
 The Gorenstein curve $\overline{\mathcal C}$ will have only isolated singularities over $\Dt$. They may have genus one or two according to $p_a(\tropD^\circ)$. 
 \begin{dfn}
 Let $\Z$ be the Cartier divisor on $\TC$, supported over $\pazocal D$, which is determined by the inclusion $\OO_{\TC}(-\lambda)\hookrightarrow\OO_{\TC}$.
\end{dfn}

\begin{rem}\label{rem:localsections}

Locally, we construct a line bundle on $\TC$ trivial on $\pazocal Z$ except when it contains the special component, and relatively ample elsewhere, as follows.
We pick smooth disjoint sections $p_1,\ldots,p_{2d}$ of $\TC$ according to the weight function, namely so that $\#\{i|p_i\in \tC_v\}=2w(v)$.
These sections exist only locally, but we will show that our construction does not depend on this choice, therefore it glues on the whole of $\TAt$.
\end{rem}

\noindent \emph{Notation.} Let $\sigma\colon\TC\to\TC$ be the hyperelliptic involution. Let $\mathbf{p}_i$ denote the multi-section $p_i+\sigma(p_i)$, and $\mathbf{p}=\sum\mathbf{p}_i$. We denote by $\Li$ the line bundle $\omega_{\TC/\widetilde{\mathcal A}}(\mathbf{p})(\lambda)$ on $\TC$.

\begin{thm}\label{thm:contraction}
 The line bundle $\Li$ is $\pi$-semiample. In the diagram
 \bcd
 \TC\ar[rr,"\phi"]\ar[dr,"\pi" below=.3cm] & &\Cp:=\underline{\operatorname{Proj}}_{\TAt}\left(\pi_*\bigoplus_{k\geq 0}\Li^{\otimes k}\right)\ar[dl,"\pi^\prime"] \\
 & \TAt &
 \ecd
 the morphism $\pi^\prime$ is a flat family of reduced, projective, Cohen-Macaulay curves of arithmetic genus two, with Gorenstein fibres outside $\Do$. Moreover, we can perform a parallel contraction $\mathcal T\to \Tp$ so that $\psi^\prime\colon \Cp\to \Tp$ remains finite. Neither $\Cp$ nor $\Tp$ depend on the choice of sections respecting the weight function, as per Remark \ref{rem:localsections}.
\end{thm}

\begin{rem}\label{rem:tailcontraction}
 Notice that $\Li$ is symmetric under the hyperelliptic involution, i.e. it is pulled back from the target of the admissible cover. Any weight-zero branch appearing as the conjugate of a rational tail (or bridge) of positive weight is therefore not contracted under $\phi$. This ensures that $\Cp$ remains ``hyperelliptic'', i.e. it admits a double cover to a rational curve, possibly with ordinary $m$-tuple points.
 As far as the factorisation of the stable map is concerned, all those unstable components of weight zero (different from the special branches) do not interfere. We could have as well carried out the construction without symmetrising the line bundle.
\end{rem}

\begin{thm}\label{thm:pushout}
 Let $\Zp$ denote the image of $\Z$ under $\phi$. Over $\Do$, it is a flat family of Gorenstein curves of arithmetic genus two, and the image of $\Zp$ under $\psi^\prime$ is a rational curve $\Tp_{\pazocal Z}$. 
 Let $\OC$ be obtained as the pushout of:
 \bcd
 \Zp_{\Do}\ar[d,"\psi^\prime"]\ar[r,hook] &\Cp \\
 \Tp_{\pazocal Z} &
 \ecd
 Then $\OC\to\TAt$ is a flat family of projective, Gorenstein curves of arithmetic genus two. Moreover, the weight of any subcurve of positive genus is at least one.
\end{thm}

The proof of the two theorems above occupies the next two sections.

\subsection{First step: the contraction}

\begin{lem}\label{lem:basechange}
 For $k\geq2$, $R^1\pi_*\Li^{\otimes k}$ is supported along $\Dt$. Moreover, it admits a two-term resolution that remains such after pullback to a sufficiently generic base: in particular, if $f\colon S\to\TAt$ is a morphism such that $\OO_S(-\rho_\text{max})\to\OO_S$ remains injective, then
 \[f^*\pi_*(\Li^{\otimes k})\to(\pi_S)_*(\Li_S^{\otimes k})\]
 is an isomorphism.
\end{lem}

\begin{proof}
\emph{Support:} Since $\CL^{\otimes k}$ is flat on the base, by cohomology and base change it is enough to show the vanishing of  $H^1(\TC_s,\CL^{\otimes k})$ for $s\in \TAt\setminus\Dt.$

By weighted stability we have that $\Li$ is relatively ample outside of $\Z$, and in particular over $\TAt\setminus \pazocal D$.

For $s\in \pazocal D\setminus\Dt$ we have $0\neq \rho_1=\rho_{\max}$.
We note that $\CL\geq0$ and, by Corollary \ref{cor:rho1=rhomax}, it has degree two on every positive genus subcurve of $\TC_s$. Then, it is clear by degree reasons that $h^1(\CL^{\otimes k})=0$ for $k\geq 2$.

\noindent\emph{Resolution:}
Note first that the rank of $R^1\pi_*\Li^{\otimes k}$ is not constant along $\Dt$.

Indeed, $\Dt$ has two types of irreducible components:
\begin{itemize}[leftmargin=.5cm]
 \item $\pazocal D_{2,1}$, where \emph{generically} $\tropD^\circ$ has genus one, and $R^1\pi_*\Li^{\otimes k}$ has rank one;
 \item $\pazocal D_{2,2}$, where $\tropD^\circ$ has genus two, and $R^1\pi_*\Li^{\otimes k}$ has rank two.
\end{itemize}
Away from their intersection, it is easy to find the desired \emph{local} resolution:
\[0\to\CO_U^{\oplus 2}\xrightarrow{\left(\substack{e^{\rho_{\max}}\;\; 0\\ 0\;\;e^{\rho_{\max}}}\right)}\CO_U^{\oplus 2}\to \CO_{\pazocal D_{2,2}\cap U}^{\oplus 2}\to 0\]
on $U\subseteq\TAt$ a neighbourhood of a generic point of $\pazocal D_{2,2}$ (assuming $\rho_1=0$), where $e^{\rho_{\max}}$ denotes a local equation for $\pazocal D$. Similarly,
\[0\to\CO_U\xrightarrow{e^{\rho_{\max}}}\CO_U\to \CO_{\pazocal D_{2,1}\cap U}\to 0\]
around a generic point in $\pazocal D_{2,1}.$

The intersection of $\pazocal D_{2,1}$ with $\pazocal D_{2,2}$ has two types of irreducible components, see Figure \ref{intersection with D2}. In order to obtain a local resolution, we adapt an argument of \cite{HLN}. We sketch it here for the reader's benefit. 

\begin{figure}
	\centering
		\begin{tikzpicture} [scale=.6]
		\tikzstyle{every node}=[font=\normalsize]
		\tikzset{arrow/.style={latex-latex}}

\coordinate (E1) at (-1,4,0);

\coordinate (E1l) at (-1,4.2,0);
\coordinate (E2) at (2,2.5,0);

\coordinate (E12) at (0.5,3.25,0);
\coordinate (E2l) at (2,2.7,0);

\coordinate (E2W) at (2.75,1.25);
\coordinate (W) at (3.5,0,0);

\coordinate (A1) at (-2.5,4.3,0);
\coordinate (B1) at (-2.5,3.7,0);

\coordinate (A2) at (.5,2.8,0);
\coordinate (B2) at (.5,2.2,0);

\coordinate (T1) at (-2.5,0,0);
\coordinate (T2) at (-2,0,0);
\coordinate (T3) at (-1.5,0,0);
\coordinate (T4) at (0.7,0,0);
\coordinate (T5) at (1.3,0,0);

\draw[->-] (W) to (E2);
\draw[->-] (E2) to (E1);

\draw (T4) to (E2);
\draw (T5) to (E2);

\draw (T1) to (E1);
\draw (T2) to (E1);
\draw (T3) to (E1);

\draw[gray] (A1) to (E1);
\draw[gray] (B1) to (E1);
\draw[gray] (B2) to (E2);
\draw[dashed] (-3,0,0)--(4,0,0);
\foreach \x in {T1,T2,T3,T4,T5} 
\fill (\x) circle (3pt);
\fill[red] (W) circle (3pt);

\draw[fill=white] (E1) circle (3.5pt);
\draw[fill=white]  (E2) circle (3.5pt);


\node at (E1l) [above] {\scriptsize{$E_1$}};
\node at (E1l) [right] {\scriptsize{$g=1$}};
\node at (E2l) [above] {\scriptsize{$E_2$}};
\node at (E2l) [right] {\scriptsize{$g=1$}};

\node at (E12) [above] {\color{blue}\scriptsize{$1$}};
\node at (E12) [below] {\scriptsize{$\zeta_1$}};

\node at (E2W) [right] {\color{blue}\scriptsize{$2$ or $3$}};
\node at (E2W) [left] {\scriptsize{$\zeta_2$}};

\node at (A1) [left] {\begin{color}{gray}\scriptsize{$A_1$}\end{color}};
\node at (B1) [left] {\begin{color}{gray}\scriptsize{$B_1$}\end{color}};
\node at (B2) [left] {\begin{color}{gray}\scriptsize{$A_2$}\end{color}};

\node at (T3) [above right] {\color{blue}\scriptsize{$1$}};
\node at (T5) [above right] {\color{blue}\scriptsize{$1$}};

\node at (0.5,-.5,0) [below]  {\scriptsize{$\pazocal D^{EE}$}};

\begin{scope}[every coordinate/.style={shift={(10,0,0)}}]

\draw[->-] ([c]W) to ([c]E2);
\draw[->-] ([c]E2)  to [out=80,in=80]([c]E1);
\draw[->-] ([c]E2)  to [out=-100,in=-100]([c]E1);

\draw ([c]T4) to ([c]E2);
\draw ([c]T5) to ([c]E2);

\draw ([c]T1) to ([c]E1);
\draw ([c]T2) to ([c]E1);
\draw ([c]T3) to ([c]E1);

\draw[gray] ([c]A1) to ([c]E1);
\draw[gray] ([c]B1) to ([c]E1);
\draw[gray] ([c]A2) to ([c]E2);

\draw[dashed] (7,0,0)--(14,0,0);
\foreach \x in {E2,T1,T2,T3,T4,T5} 
\fill ([c]\x) circle (3pt);
\fill[red] ([c]W) circle (3pt);

\draw[fill=white] ([c]E1) circle (3.5pt);


\node at (8.9,4.1,0) [above] {\scriptsize{$E_1$}};
\node at (9,4,0) [right] {\scriptsize{$g=1$}};
\node at([c]E2) [right] {\scriptsize{$E_2$}};

\node at (11,4,0) [above] {\color{blue}\scriptsize{$1$}};
\node at (10.25,3,0)  [below] {\color{blue}\scriptsize{$1$}};
\node at (10.5,2.35,0) [below] {\scriptsize{$\zeta_1$}};

\node at ([c]E2W) [right] {\color{blue}\scriptsize{$2$ or $3$}};
\node at ([c]E2W) [left] {\scriptsize{$\zeta_2$}};

\node at ([c]A1) [left] {\begin{color}{gray}\scriptsize{$A_1$}\end{color}};
\node at ([c]B1) [left] {\begin{color}{gray}\scriptsize{$B_1$}\end{color}};
\node at ([c]A2) [above left] {\begin{color}{gray}\scriptsize{$A_2$}\end{color}};

\node at ([c]T3) [above right] {\color{blue}\scriptsize{$1$}};
\node at ([c]T5) [above right] {\color{blue}\scriptsize{$1$}};

\node at (10.5,-.5,0) [below]  {\scriptsize{$\pazocal D^{B}$}};

\end{scope}

\end{tikzpicture}

 \caption{The generic points of $\pazocal D_{2,1}\cap\pazocal D_{2,2}$.}
 \label{intersection with D2}
\end{figure}
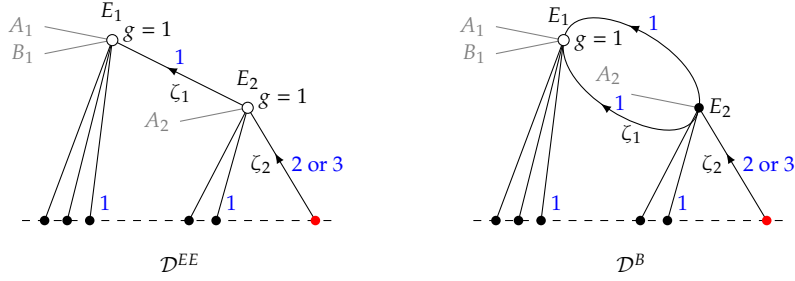


Locally on the base, there is a section $p_{2d+1}$ such that, setting $\mathbf{p}_{2d+1}=p_{2d+1}+\sigma(p_{2d+1})$, we have $\omega_C(\lambda)=\OO_C(\mathbf{p}_{2d+1})$, and therefore
$\mathcal L=\OO_C(\sum_{i=1}^{2d+1} \mathbf{p}_i)$.

Again locally on the base, we may find disjoint generic sections $A_1$, $A_2$, and $B$, such that $A_1$ and $B$ pass through $E_1$, and $A_2$ passes through $E_2$. Then $\Li(A_1+A_2-B)$ has vanishing $h^1$ on fibres, and therefore, by cohomology and base-change, $\pi_*\Li(A_1+A_2-B)$ is a vector bundle. Locally, we can write $\pi_*\Li\simeq\OO_U\oplus\pi_*\Li(-B)$, and the second factor is the kernel of the evaluation map:
\[\pi_*\Li(A_1+A_2-B)\to\pi_*(\OO_{A_1}(A_1)\oplus\OO_{A_2}(A_2)).\]
The evaluation map can be studied fibrewise, since both sheaves in question have vanishing $h^1$; moreover, the source can be decomposed into:
\[\pi_*\Li(A_1+A_2-B)\cong\bigoplus_{i=1}^{2d+2}\pi_*\OO(\mathbf p_i+A_1+A_2-B),\]
and the evaluation map can be studied componentwise:
\[\ev_{i,j}\colon \pi_*\OO(\mathbf p_i+A_1+A_2-B)\to\pi_*(\OO_{A_j}(A_j))\quad i=1,\ldots,2d+2;j=1,2. \]
The cokernel of the latter is $H^1(C,\OO_C(\mathbf p_i+A_{2-j}-B))$. See \cite[\S 4.2]{HL}.

It follows from an argument similar to that of \cite[Proposition 4.13]{HL} that the latter is non-zero precisely when $p_i$ stays away from $E_j$, therefore, in some local trivialisation of the line bundles involved,
\[\ev_{i,j}=c_{i,j}\prod_{q\in[p_i,A_j]}\zeta_q\]
where $c_{i,j}\in\OO_{U}^*$, we denote by $[p_i,A_j]$ the set of nodes separating $p_i$ from $A_j$, and $\zeta_q\in\OO_{U}$ is the smoothing parameter of the node $q$. Thanks to the alignment,
\begin{itemize}[leftmargin=.5cm]
 \item the smoothing parameter $\zeta_2$ of the node separating $E_2$ from the component supporting $D$ divides all the expressions of the form $\ev_{i,2},\ i=1,\ldots,2d+2$;
 \item if $\zeta_1$ denotes the smoothing parameter of the node separating $E_1$ from $E_2$, the product $\zeta_1\zeta_2$ divides all the expressions of the form $\ev_{i,1},\ i=1,\ldots,2d+2$.
\end{itemize}
(Notice that $\zeta_1$ and $\zeta_2$ should be replaced by some products of smoothing parameters when the curve degenerates.)
We can therefore use the column of the evaluation matrix associated to a marking $p_i$ (up to relabelling, $i=1$) on the component supporting $D$ in order to put the matrix in triangular form. In order to diagonalise it, we need a more refined information that we borrow from \cite[\S 2.6-7]{HLN}, and we restate here in streamlined form:

\smallskip

\noindent\textbf{Proposition}[Hu-Li-Niu] \emph{The determinant of the matrix:}
\[\begin{pmatrix}
   c_{1,1} & c_{j,1} \\
   c_{1,2} & c_{j,2}
  \end{pmatrix}\]
\emph{is invertible when the markings $p_1$ and $p_j$ are not conjugate under $\psi$. }

\smallskip

Since the component supporting $D$ contains at least three markings ($p_1,p_{2d+1}$, and $p_{2d+2}$), we can find two non-conjugate ones.

Summing up, the evaluation matrix can be put in the following form:
\[\begin{pmatrix}
   \zeta_1 & 0 \\
   0 & \zeta_1\zeta_2
  \end{pmatrix}\]
(the remaining columns are zero). Noticing that $\pazocal D_{2,i}=\{\zeta_i=0\},\ i=1,2$, around the given point, we have thus found the desired local resolution of $\R^1\pi_*\Li^k$. In particular, it follows that $\pi_*\Li^{\otimes\geq2}$ is a vector bundle on $\TAt$.

\medskip

\noindent\emph{Base-change:} see the proof of \cite[Lemma 3.7.5.3]{RSPW1}. Note that the evaluation matrix above is an explicit instance of the Grothendieck-Mumford's complex for cohomology and base-change, and $\zeta_1\zeta_2=e^{\rho_{\text{max}}}$ is a local equation for $\pazocal D$.

\end{proof}

\begin{proof}[Proof of Theorem~\ref{thm:contraction}] Analogous to \cite[Lemma 2.13]{SMY1} and \cite[Proposition 3.7.6.1]{RSPW1}. We recap for the reader's convenience.

\noindent\emph{Flatness}: is equivalent to $\pi_*\Li^{\otimes k}$ being a vector bundle for $k\gg0$ \cite[\href{https://stacks.math.columbia.edu/tag/0D4D}{Tag 0D4D}]{stacks-project}.

\noindent\emph{Basepoint freeness}: i.e. existence of the morphism $\phi$. This is clear outside $\pazocal D$ where  $\CL^{\otimes k}$ is $\pi-$ ample. Even on $\pazocal D$, we have that $\CL^{\otimes k}$ is $\pi-$ample outside of $\pazocal Z$. Moreover, from the short exact sequence:
\[0\to\CL(-\lambda)^{\otimes k}\to \CL^{\otimes k}\to\CL_{\rvert k\pazocal Z}^{\otimes k}\to 0\]
and the vanishing of $R^1\pi_*\CL(-\lambda)^{\otimes k}$ (stability), it is enough to show that for any $x\in\pazocal Z$ there exists a section of $\pi_*\CL_{\rvert k\pazocal Z}^{\otimes k}$ around $\pi(x)$ which does not vanish in $x$. By definition of $\lambda$, the line bundle $\CL\rvert_{\pazocal Z}$ is the pullback along $\psi$ of a line bundle of degree one and non-negative multidegree on $\mathcal T_{\pazocal Z}$; the latter has enough sections.

\noindent\emph{Properties of the fibres}: can be studied after base-change to a generic trait $\dvr$. We may assume that the generic point corresponds to a smooth curve, and the closed point maps to $\pazocal D=\Do\cup\Dt$. Thus
\[\phi_\dvr\colon C:=\mathcal C_\dvr\to \mathcal C_\dvr^\prime=:C'\]
is a birational contraction satisfying $\phi_*\OO_C=\OO_{C^\prime}$. It follows that $C^\prime$ is a normal surface.
In particular, the central fibre is $S1$; it is also generically reduced, being birational to that of $C$, thus it is reduced (and Cohen-Macaulay).

Finally, we want to argue that the fibres are Gorenstein outside $\Do$. We may assume that the special point of $\dvr$ maps to $\Dt\setminus\Do$. Then $\Li$ is trivial along $Z$, which is therefore contracted to a codimension two locus $Z^\prime$ of $C^\prime$. Outside of $Z$, $\phi$ restricts to an isomorphism. The equality of line bundles:
\[\OO_{C^\prime}(1)(-\sum \mathbf p_i)_{|C^\prime\setminus Z^\prime}=\omega_{C/\dvr|C\setminus Z}=\omega_{C^\prime/\dvr|C^\prime\setminus Z^\prime},\]
together with the fact that $\omega_{C^\prime/\dvr}$ is an $S2$ sheaf, shows  that the latter coincides with the line bundle $\OO_{C^\prime}(1)(-\sum \mathbf p_i)$ on the whole of $C^\prime$ (Hartogs' Theorem). Thus, $C^\prime$ is Gorenstein over $\dvr$.

\noindent\emph{Compatible contraction of $\mathcal T$}: can be performed using the line bundle \[\omega_\mathcal T(D_B)(\lambda_T)(\psi(\mathbf p)),\]
whose pullback to $\TC$ is $\mathcal L$.

\noindent\emph{Well-posedness}: the construction of $\Cp$ and $\Tp$ is independent of the choice of markings respecting the weight function by birational rigidity \cite[Lemma 1.15]{Debarre}.
\end{proof}

\subsection{Second step: the pushout}
Ribbons can be used to interpolate between isolated singularities over $\Do$ - this is completely natural from the point of view of piecewise-linear functions on the tropical side, and serves as a correction of the failure of $\Cp$ at being Gorenstein.

Recall that $\pazocal Z'$ was defined as the image of $\pazocal Z$ under $\phi$. First, we prove that the definition of $\pazocal Z'$ commutes with base-change to a generic trait; then, we will show that the pushout construction commutes with such a base-change, and thus we may reduce to the case of surfaces in order to study the singularities of the fibres.

\begin{prop}\label{prop:basechangeZ'}
Let $\pazocal Z'$ be the subscheme of $\Cp$ defined by the ideal sheaf $\phi_*\CO_{\mathcal C}(-\lambda)$ and supported on $\pazocal D.$
Then we have:
\begin{enumerate}
    \item  $R^1\phi_*\CO_{\mathcal C}(-\lambda)=0$, in particular $\phi_*\CO_{\mathcal C}(-\lambda)=\operatorname{Fitt}(\phi_*\CO_{\pazocal Z})$;
    \item For every generic trait $\dvr\xrightarrow{\iota}\widetilde{\pazocal A}$ with generic point mapping to the smooth locus, the definition of  $\pazocal Z'$ commutes with base change, i.e., $\iota^*\phi_*\CO_{\mathcal C}(-\lambda)=\phi_{\dvr *}\CO_{\mathcal C_{\dvr}}(-\lambda)$.   
\end{enumerate}
The analogous statements about $\mathcal T^\prime_{\pazocal Z}$ hold as well.
\end{prop}
\begin{proof} The discussion has been somewhat inspired by \cite[\S 1]{Teissier}.
\begin{enumerate}[leftmargin=0cm,itemindent=0.8cm,labelwidth=\itemindent,labelsep=0cm,align=left]
    \item Let $C_s$ be a fiber on which $\phi_s$ is not an isomorphism; in particular $\lambda_s\neq 0$. 
    Working locally on the base, we can choose smooth and disjoint sections $p_1\ldots,p_d$ of $\mathcal C$ respecting the weight function.
    It is enough to prove that $R^1\phi_*\CO_{\mathcal C}(-\lambda)(k\mathbf{p})=0$ for $k\gg0.$ Indeed, once we know the latter vanishing, since $\phi$ is an isomorphism around $\mathbf{p},$ we have that $\CO_{\mathcal C}(k\mathbf{p})=\phi^*\CO_{\mathcal C'}(k\mathbf{p}')$ and thus
     \begin{align*}
     0&=R^1\phi_*\CO_{\mathcal C}(-\lambda)(k\mathbf{p})=H^1(R\phi_*\CO_{\mathcal C}(-\lambda)(k\mathbf{p}))\\
     &=H^1(R\phi_*\CO_{\mathcal C}(-\lambda)\otimes\phi^*\CO_{\mathcal C'}(k\mathbf{p}'))=R^1\phi_*\CO_{\mathcal C}(-\lambda)\otimes\CO_{\mathcal C'}(k\mathbf{p}')
     \end{align*}
     implying the desired vanishing.
     
    From the spectral sequence computing $R\pi'_*\circ R\phi_*$, the five-terms exact sequence:
    \[0\to R^1\pi'_*\phi_*\CO_{\mathcal C}(-\lambda)(k\mathbf{p})\to R^1\pi_*\CO_{\mathcal C}(-\lambda)(k\mathbf{p})\to\pi'_*R^1\phi_*\CO_{\mathcal C}(-\lambda)(k\mathbf{p})\]
    ends there, because the next term would involve an $\R^2\pi^\prime_*(-)$, that vanishes as the fibre dimension is bounded by one.
    
    Notice that if $ R^1\pi_*\CO_{\mathcal C}(-\lambda)(k\mathbf{p})\otimes k(s)=0$ for those $s$ where $\phi_s$ is not an isomorphism,
    so is $\pi'_*R^1\phi_*\CO_{\mathcal C}(-\lambda)(k\mathbf{p})\otimes k(s).$ 
    
    On the other hand, $R^1\phi_*\CO_{\mathcal C}(-\lambda)(k\mathbf{p})$ 
is a coherent sheaf supported on the closed locus $V'$ in $\mathcal C'$ over 
    which the fiber of $\phi$ has positive dimension. The latter is finite over $\widetilde{\pazocal A}$, and thus $\pi'_*R^1\phi_*\CO_{\mathcal C}(-\lambda)(k\mathbf{p})\otimes k(s)=0$ implies $R^1\phi_*\CO_{\mathcal C}(-\lambda)(k\mathbf{p})\otimes k(x)=0$ for each $x\in V'_s.$

    Since $R^1\pi_*\CO_{\mathcal C}(-\lambda)(k\mathbf{p})$ satisfies cohomology and base change, the vanishing can be checked after restricting to a fibre. Let $Z_s$ denote the support of $\lambda_s$, $C_i$ the trees rooted at the vertices of $\partial\tropD$, and $e_i$ the first edges encountered in $\tropD^\circ$. Taking the normalisation of $C_s$ at the nodes corresponding to the edges $e_i$ we get
\[0\to \CO_{C_s}(-\lambda)(k\mathbf{p})\to \CO_{Z_s}(-\lambda)\oplus\bigoplus_i \CO_{C_i}(-\lambda)(k\mathbf{p})\to\bigoplus_i\mathbb C_{e_i} \to 0.\]
The evaluation map on the nodes $e_i$ is clearly surjective at the level of $H^0$, as the line bundle restricted to $C_i$ is very ample for $k$ big enough, and the desired vanishing follows from that of $H^1(\CO_{Z_s}(-\lambda_s))$, which can be argued by the definition of $\lambda$ and Serre duality.

\item We know by Lemma \ref{lem:basechange} that the construction of $\mathcal C'$ commutes with the given base change, so the following diagram is Cartesian:
\[\begin{tikzcd}
\mathcal C_{\dvr}\ar[r,"\phi_{\dvr}"]\arrow[d,"\iota"]\ar[dr,phantom,"\Box"] & \mathcal C'_{\dvr}\arrow[d,"\iota' "]\\
\mathcal C\ar[r,"\phi"] & \mathcal C'.
\end{tikzcd}\]
Furthermore, since the source and target are smooth, $\dvr\hookrightarrow\widetilde{\pazocal A}$ is an l.c.i morphism, and so are $\iota$ and $\iota^\prime$; it thus follows from \cite[Corollary~2.27]{kuznetsov2006hyperplane} that
\[L\iota'^*R\phi_*\CO_{\mathcal C}(-\lambda)=R\phi_{\dvr *}L\iota^*\CO_{\mathcal C}(-\lambda).\]
On the other hand, the higher pushforward vanishes by the previous point, so
\[R\phi_*\CO_{\mathcal C}(-\lambda)=\phi_*\CO_{\mathcal C}(-\lambda),\] and  $L\iota^*\CO_{\mathcal C}(-\lambda)=\CO_{\mathcal C_{\dvr}}(-\lambda)$ because it is a line bundle, hence the derived statement is equivalent to what we want.

\end{enumerate}
\end{proof}

\begin{prop}\label{prop:Goreinstein+flatness.}
The restriction of $\pazocal Z^\prime$ to $\Do$ is a flat family of Gorenstein curves of genus two. Similarly, $\mathcal T^\prime_{\pazocal Z}$ is a flat family of Gorenstein curves of genus zero (i.e. at worst nodal).
\end{prop}

Notice that $\Cp$ is not always Gorenstein over $\Do$, in particular it can be the decomposable union of $\pazocal Z^\prime$ with some lines. For the proof we need the following:

\begin{lem}\label{lem:Dreduced}
$\Do$ is a reduced divisor.
\end{lem}
\begin{proof}
$\Do$ is a Cartier divisor in a smooth ambient space, so it is enough to check that it is generically reduced.

The generic point of $\Do$ looks like in Figure \ref{fig:genericD1}.
\begin{figure}[h]
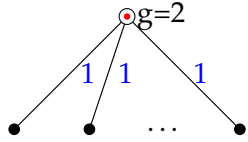

 \tikz{
 \draw[fill=black] (-1,0) circle (2pt) --node[right]{\color{blue}{1}} (.5,1.5) (0,0) circle (2pt) --node[right]{\color{blue}{1}} (.5,1.5) (1,0) node{\ldots}  (2,0) circle (2pt) --node[right]{\color{blue}{1}} (.5,1.5);
 \draw[fill=white] (.5,1.5) circle (3pt) node[right]{g=2};
 \draw[red,fill=red] (.5,1.5) circle (1pt);
 }
 \caption{$\lambda$ at the generic point of $\Do$.}
 \label{fig:genericD1}
\end{figure}

Therefore, a generic trait with uniformiser $t$ will intersect $\Do$ in $(t)$.
 \end{proof}

\begin{proof}[Proof of Proposition~\ref{prop:Goreinstein+flatness.}]

We may change the base to a generic trait $\dvr$ with closed point mapping to the given point of $\Do$:
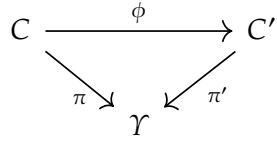
\[\begin{tikzcd}
 C\ar[rr,"\phi"] \arrow[rd,"\pi" below=.1cm] & & C'\ar[dl,"\pi' "]\\
 &\dvr
\end{tikzcd}\]
so $C$ is a smooth surface, $\phi$ is a birational contraction, which is an isomorphism outside the divisor $Z$ defined by $\mathcal I_{Z}=\CO_{C}(-\lambda)$.
Notice that  there is a subcurve of $Z$ on which the line bundle defining the contraction is ample, therefore $Z^\prime\subseteq C^\prime$ has pure codimension one. We want to show that $\omega_{Z'}$
is a line bundle and that $\chi(\omega_{Z'})=1$ (equivalently, $p_a(Z')=2$).

Recall that, by Grothendieck duality:
\[\omega_{Z'}=\mathcal E xt^1_{C'}(\OO_{Z'},\omega_{C'})=\mathcal Hom_{C}(\mathcal I_{Z'},\omega_{C'})\rvert_{Z'}.\]

By adjunction for the Cartier divisor $Z\subseteq C$, there is a short exact sequence:
\begin{equation}\label{eqn:adjunctionZ}
0\to\omega_C\to\omega_C(\lambda)\to\omega_Z\to0,
\end{equation}
which stays exact after pushforward along $\phi$ by Grauert-Riemenschneider vanishing \cite[Corollary 2.68]{kollar2008birational}. By applying $\mathcal Hom_{\OO_{C'}}(-,\omega_{C'})$ to the exact sequence:
\[0\to \phi_*\OO_C(-\lambda)\to\OO_{C'}\to\OO_{Z'}\to 0,\]
we obtain the bottom row of the following diagram:
\bcd
0\ar[r] & \phi_*\omega_C\ar[r]\ar[d]&\phi_*\omega_C(\lambda)\ar[r]\ar[d] &\phi_*\omega_Z\ar[d]\ar[r] & 0\\
0\ar[r] & \omega_{C'}\ar[r] & \mathcal Hom_{C}(\phi_*\OO_C(-\lambda),\omega_{C'})\ar[r] &\omega_{Z'}\ar[r] & 0
\ecd
By Grothendieck duality ($\phi_*\dashv\phi^!$ and $\phi^!\omega_{C^\prime}=\omega_C$) and the snake lemma, we conclude that the vertical arrows are isomorphisms. This implies that $\omega_{Z'}$ is a line bundle if $\phi_*\omega_C(\lambda)$ is.
Since the sections $\mathbf{p}$ are away from $\operatorname{Exc}(\phi)$, by construction \[\phi_*\omega_C(\lambda)=\CO_{C'}(1)(-\phi(\mathbf{p}))\]
is a line bundle. We have thus proved that $Z'$ is Gorenstein.

Moreover, $\phi_*\omega_Z=\omega_{Z'}$. Now, to prove that $\chi(\omega_{Z'})=1$ it is enough to prove that $\chi(\omega_{Z'})=\chi(\omega_Z)$. Indeed, over $\Do$, the definition of $\lambda$ and adjunction \eqref{eqn:adjunctionZ} show that $\omega_C(\lambda)$ restricts to a line bundle of degree two on $Z$, which is therefore a curve of genus two. Since $\phi_*\OO_Z=\OO_{Z'}$, it is enough to prove the vanishing of $\R^1\phi_*\OO_Z$. By Proposition \ref{prop:basechangeZ'}, $\R^1\phi_*\OO_Z\simeq\R^1\phi_*\OO_C$, therefore the desired statement is equivalent to the fact the $C'$ has rational singularities. This follows from $\phi_*\omega_C=\omega_{C'}$, see for example \cite{kovacs}. We have thus proved that $Z'$ has genus two.

In order to prove the flatness of $\pazocal Z'_{\Do}\to\Do$, it is sufficient to show that \[\pi'_*(\CO_{\pazocal Z'_{\Do}}\otimes\CO_{\mathcal C'}(n))\]  is a vector bundle on $\Do$ for $n$ large enough. Since $\Do$ is a reduced divisor (Lemma \ref{lem:Dreduced}), we only have to show that its rank is constant along $\Do$. It is easy to see that:
\[\CO_{\mathcal C'}(1)_{\rvert\pazocal Z'_{\pazocal D_1}}\cong\omega_{\pazocal Z'_{\pazocal D_1}}.\]

Since $\pazocal Z^\prime_{\Do}$ is a curve of genus two, it follows from Riemann-Roch that
\[H^0(Z'_s,\omega_{Z'_s}^{\otimes n})=2n-1,\]
so it is enough to show that $\CO_{\mathcal C'}(1)_{\rvert\pazocal Z'_{\pazocal D_1}}$ satisfies cohomology and base-change. 
For this, we observe the short exact sequence:
\bcd
0\ar[r] &\pi'_*(\phi_*\CO_{\mathcal C}(-\lambda)\otimes\CO_{\mathcal C'}(n))\ar[d,"\cong"]\ar[r] &\pi'_*\CO_{\mathcal C'}(n)\ar[r] &\pi'_*\CO_{\mathcal C'}(n)\rvert_{\pazocal Z'}\ar[r] & 0\\
& \pi_*(\omega_{\pi}^{\otimes n}(n\mathbf{p})\otimes\CO_{\mathcal C}((n-1)\lambda))
\ecd
As $\mathcal C$ and $\mathcal C^\prime$ are flat over the base, and both $\omega_\pi(\mathbf{p})$ and $\OO_{\mathcal C^\prime}(1)$ are relatively ample, the first two bundles satisfy cohomology and base-change. It follows from a diagram-chase that so does the third. We have thus proved that $\pazocal Z^\prime$ is flat on $\Do$.

The statement about $\mathcal T^\prime_{\pazocal Z}$ can be proven in an analogous (but easier) fashion.
\end{proof}

\begin{proof}[Proof of Theorem \ref{thm:pushout}]
Noticing that $\Zp\to\mathcal T^\prime_{\pazocal Z}$ is finite, the existence of the pushout as a scheme over $\TAt$ follows from results of D. Ferrand \cite{Ferrand}.
We have already proved that the construction of $\Cp,\ \Zp,$ and $\mathcal T^\prime_Z$ commutes with pullback to a generic trait. The pushout does as well in virtue of \cite[\href{https://stacks.math.columbia.edu/tag/0ECK}{Tag 0ECK}]{stacks-project}. So, in order to prove that the fibres of $\OC$ are Gorenstein, we may work with fibred surfaces, in which case we may apply some results of M. Reid \cite{reid1994nonnormal}.

Following \cite[\S~2.1]{reid1994nonnormal}, $\bar\phi\colon C^\prime\to\overline{C}$ is the normalisation, with conductor:
\[\operatorname{Ann}(\bar\phi_*\OO_{C'}/\OO_{\overline C})=\operatorname{Ann}(\psi_*\OO_{Z'}/\OO_{T'_Z}).\]

Since $\psi$ is a double cover, it is in particular flat, and $\psi_*\OO_{Z'}/\OO_{T'_Z}$ is a line bundle; it follows from \cite[Proposition~2.2]{reid1994nonnormal} that $\overline{C}$ is $S2$.

Moreover, we have seen in the proof of Proposition \ref{prop:Goreinstein+flatness.} that \[\omega_{C'}(Z')=\mathcal Hom_{\mathcal C'}(\mathcal I_{\pazocal Z'},\omega_{\mathcal C'})\] is a line bundle on $C'$. Since $\psi\colon Z'\to T'_Z$ is a double cover of a rational curve, it follows that the kernel of the canonical map $\psi_*\omega_{Z'}\to\omega_{T'_Z}$ is a line bundle as well. The criterion of \cite[Corollary~2.8~(iv)]{reid1994nonnormal} allows us to conclude that $\overline C$ is Gorenstein.

Finally, the statement about weights is obvious from the construction of $\lambda$ in case its support is connected. If there are two components of $\tropD^\circ$, each of them necessarily of genus one by Lemma \ref{lem:supp_pos_gen}, the statement follows from Lemma \ref{lem:two_genus_one_pos_wt}.
\end{proof}

\subsection{The secondary construction}
We restrict now to $\doublewidetilde{\mathcal A_2}$ defined in \ref{def:doubletildeA}. Thanks to the alignment, we have the following

\begin{dfn}
 Let $\Do^2\subseteq \doublewidetilde{\mathcal A_2}$ be the logarithmic divisor where the weight of the subcurve $\tropD$ determined by $\lambda$ is two.
\end{dfn}

The curve $\overline{\mathcal C}_{|\doublewidetilde{\mathcal A_2}}$ will be manipulated further to insert ribbons over $\Do^2$. We could descend the line bundle $\OO_{\mathcal C'}(1)$ to $\overline{\mathcal C}$ using the results of \cite[Theorem 2.2]{Ferrand} and proceed from there. For simplicity we chose instead to work on $\TTC$ from Theorem \ref{thm:doubletildeA}, so that the phrasing of the following two theorems is analogous to that of the main theorems in \S\ref{sec:main}. The proof goes along the same lines too, so we leave it to the reader to figure out the details.

Let $\LLi$ denote the line bundle $\omega_{\TTC/\doublewidetilde{\mathcal A_2}}(\mathbf{p})(\widetilde\lambda)$ on $\TTC$, where $\mathbf{p}$ is a local multi-section compatible with the weight function, see the notation at the beginning of \S\ref{sec:main}.

\begin{thm}\label{thm:contraction2}
 The line bundle $\LLi$ is $\tilde\pi$-semiample. In the diagram
 \bcd
 \TTC\ar[rr,"\tilde \phi"]\ar[dr,"\tilde\pi" below=.3cm] & &\Cpp:=\underline{\operatorname{Proj}}_{\TTA}\left(\tilde\pi_*\bigoplus_{k\geq 0}\LLi^{\otimes k}\right)\ar[dl,"\pi^{\prime\prime}"] \\
 & \TTA &
 \ecd
 the morphism $\pi^{\prime\prime}$ is a flat family of reduced, projective, Cohen-Macaulay curves of arithmetic genus two, with Gorenstein fibres outside $\Do\cup\Do^2$. Moreover, we can perform a parallel contraction $\mathcal T\to \Tpp$ so that $\psi^{\prime\prime}\colon \Cpp\to \Tpp$ remains finite. Neither $\Cpp$ nor $\Tpp$ depend on the choice of local sections $\mathbf{p}$.
\end{thm}

Let $\tZ$ be the logarithmic divisor on $\TTC$ defined by $\OO_{\TTC}(-\widetilde\lambda)\hookrightarrow\OO_{\TTC}$.

\begin{thm}\label{thm:pushout2}
 Let $\Zpp$ denote the image of $\tZ$ under $\tilde\phi$. Over $\Do\cup\Do^2$, it is a flat family of Gorenstein curves of arithmetic genus two, and the image of $\Zpp$ under $\psi^{\prime\prime}$ is a rational curve $\Tpp_{\pazocal Z}$. 
 Let $\OOC$ be obtained as the pushout of:
 \bcd
 \Zpp_{\Do\cup\Do^2}\ar[d,"\psi^{\prime\prime}"]\ar[r,hook] &\Cpp \\
 \Tpp_{\pazocal Z} &
 \ecd
 Then $\OOC\to\TTA$ is a flat family of projective, Gorenstein curves of arithmetic genus two. Moreover, the weight of the minimal subcurve of genus two is at least three.
\end{thm}

\begin{rem}
Notice that $\Do$ and $\Do^2$ do not intersect, thus the flatness of $\Zpp$ can be checked independently on the two components: to this goal, we notice that $\mathcal O_{\Cpp}(1)$ restricts to $\omega_{\Zpp_{\Do}/\Do}$ on $\Zpp_{\Do}$ (as seen in the previous section) and to  $\omega^{\otimes 3}_{\Zpp_{\Do^2}/\Do^2}$  on $\Zpp_{\Do^2}$. The push-out construction can be carried out independently on these two loci.
\end{rem}

\begin{rem}
 We could have extended $\lambda$ to $\widetilde\lambda$ over the whole $\TAt$ by choosing a different cutoff level whenever $w(\tropD)\leq 2$. In this case, though, we would have had to deal with the cases that the weight is $1$, or $1+1$, or $2$ but supported on a different vertex than $v$. In these cases, the construction above produces singularities worse than tailed ribbons. When there are two vertices in $\tropD^\circ$ on which $\Li$ has positive degree, the singularity looks like a chain of two ribbons on an underlying node, with local equations $\mathbb C[\![x,y]\!]/(x^2y^2)$. When there are three, if the contraction $\phi$ acts non-trivially, there may even be three ribbons with underlying curve an ordinary $3$-fold point. We chose to keep the singularities under control by discarding the bad locus in Definition \ref{def:niceopen}. We shall see below that this forces an intermediate step on us (which can actually have some interest of its own), but it does not affect the end result.
\end{rem}

\subsection{Markings}\label{sec:admmarkings}
To avoid  overloading the notation and the exposition,  we have so far considered only weighted admissible covers without markings. However, with the application to stable maps in mind, markings are necessary to impose cohomological constraints using the evaluation maps.
Our construction extends to the marked version essentially unchanged. We can  consider:
\begin{dfn}
  A weighted admissible cover with markings consists of: \[\left(\psi\colon(C,D_R,\mathbf{x})\to(T,D_B,\mathbf{y}=\psi(\mathbf{x})),w\colon V(\tropC)\to\mathbb N\right)\]  such that $D_R$ and $\mathbf{x}$ are separately disjoint (multi-)sections of $C$, and $T$ is weighted-stable, i.e every weight-zero component has at least three special (marked, branch, or nodal) points. We denote the moduli space of weighted admissible covers with $n$ markings by $\mathcal A_{2,n}^{\text{wt}}$.
\end{dfn}
Markings are represented by infinitely long legs on $\tropC$. They play no role in the alignment: admissible functions will have slope $1$ along them, and the infinite legs will be subdivided accordingly. In particular, if a marking is supported on a vertex of $\tropD^\circ$, then we may \emph{sprout} (blow-up) the marking as many times as it is necessary for its strict transform to be supported on a vertex of $\partial\tropD$. (In fact, only one blow-up is necessary if it is allowed to be weighted.)
\begin{thm}\label{thm:markings_sm}
There exists a logarithmically \'etale modification $\widetilde {\mathcal A}_{2,n}\to\mathcal A_{2,n}^{\text{wt}}$ (resp. $\doublewidetilde{\mathcal A_{2,n}}$) parametrising weighted admissible covers with markings and a primary (resp. secondary) alignment. The moduli space $\widetilde {\mathcal A}_{2,n}$ (resp. $\doublewidetilde{\mathcal A_{2,n}}$) is a logarithmically smooth stack with locally free logarithmic structure (and therefore smooth).
\end{thm}

In particular, the strict transform of the markings never touches the singularity.
\begin{thm}\label{thm:markings_prim}
 There is a diagram of flat families of projective, Gorenstein curves of arithmetic genus two, with $n$ corresponding disjoint sections of the smooth locus:
 \bcd
 & \widetilde{\mathcal C}\ar[dr]\ar[dd]\ar[dl] & \\
 \mathcal C\ar[dr] & & \OC\ar[dl] \\
 & \TAtn\ar[uu,bend right,"\tilde{\mathbf x}"] &
 \ecd
 Moreover, the weight of every subcurve of positive genus in $\OC$ is at least one.
\end{thm}

\begin{thm}\label{thm:markings_sec}
 There is a diagram of flat families of projective, Gorenstein curves of arithmetic genus two, with $n$ corresponding disjoint sections of the smooth locus:
 \bcd
 & \TTC\ar[dr]\ar[dd]\ar[dl] & \\
 \mathcal C\ar[dr] & & \OOC\ar[dl] \\
 & \doublewidetilde{\mathcal A_{2,n}}\ar[uu,bend right,"\tilde{\tilde{\mathbf x}}"] &
 \ecd
 Moreover, the weight of the core of $\OOC$ is at least three.
\end{thm}

\section{A modular desingularisation of $\M{2}{n}{\PP^r}{d}^\text{main}$}\label{sec:maps}

\subsection{Irreducible components}\label{sec:irrcompos}
We draw the weighted dual graph of the general member of all possible irreducible components of $\overline{\pazocal{M}}_2(\PP^r,d)$.
Our running convention is that a white vertex corresponds to a contracted component, a gray one to a genus two subcurve covering a line two-to-one, and a black vertex corresponds to a non-special subcurve. Vertices are labelled with their genus and weight.
\begin{enumerate}[leftmargin=.6cm]
 \item \emph{main} is the closure of the locus of maps from a smooth curve of genus two;
 \item $\pazocal D^{(d_1,\ldots,d_k)}=\left\{\Dk\right\}$
 \item ${}^{\rm{hyp}}\!\pazocal D^{(d_1,\ldots,d_k)}=\left\{\hypell\right\}$
 \item ${}^{d_0}\!\pazocal E^{(d_1,\ldots,d_k)}=\left\{\Ek\right\}$
 \item ${}^{(d_{1,1},\ldots,d_{1,k_1})}\!\pazocal E^{d_0}\pazocal E^{(d_{2,1},\ldots,d_{2,k_2})}=\left\{\Ekk\right\}$
 \item ${}^{\rm{br}=d_0}\!\pazocal E^{(d_1,\ldots,d_k)}=\left\{\pacman\right\}$
\end{enumerate}
This is taken from the first author's PhD thesis \cite{Battistella-thesis}, and is implicit in \cite{HLN}.

\subsection{Factorisation through a Gorenstein curve}

\begin{dfn}
 Let $(X,\OO_X(1))$ be a polarised variety, and $\beta\in H^+_2(X,\ZZ)$ an effective curve class. The moduli space $\mathcal A_{2,n}(X,\beta)$ of \emph{admissible maps} to $X$ is defined by the fibre diagram:
 \bcd
 \mathcal A_{2,n}(X,\beta)\ar[d]\ar[r]\ar[dr,phantom,"\Box"] & \overline{\pazocal M}_{2,n}(X,\beta)\ar[d]\\
 \mathcal A_{2,n}^{\text{wt}}\ar[r] & \mathfrak{M}_{2,n}^{\text{wt}}.
 \ecd
 
 \noindent It can be described as a space of maps and admissible covers with the same source:
 \bcd
 (C,D_R,\mathbf{x})\ar[d,"\psi"]\ar[r,"f"] & \PP^r\\
 (T,D_B,\psi(\mathbf{x}))
 \ecd
 subject to the stability condition that, if $\sigma\colon C\to C$ is the hyperelliptic involution,
 \[\omega_C(\mathbf{x})\otimes f^*\OO_X(2)\otimes\sigma^*f^*\OO_X(2)\]
 is relatively ample. 
\end{dfn}

If we set $d=\OO_X(1)\cdot\beta$, we can restrict to the component of the base where the weight is $d$.

\begin{rem}\label{rem:mainbirational}
 On the locus of maps from a smooth curve, $\mathcal A_{2,n}(X,\beta)\to\overline{\pazocal M}_{2,n}(X,\beta)$ is an isomorphism; therefore, for $X=\PP^r$, the main components are birational.
\end{rem}

\begin{dfn}
 Let $\TAtn(X,\beta)$ denote the fibre product:
 \bcd
 \TAtn(X,\beta)\ar[d]\ar[r]\ar[dr,phantom,"\Box"] & \mathcal A_{2,n}(X,\beta)\ar[d]\\
 \TAtn\ar[r] & \mathcal A_{2,n}^{\text{wt}}.
 \ecd
 We call it the moduli space of \emph{aligned admissible maps}.
\end{dfn}
\begin{rem}\label{rem:logetale}
 $\TAtn(X,\beta)$ is logarithmically \'{e}tale over $\mathcal A_{2,n}(X,\beta)$. It comes with universal structures:
 \bcd[cramped]
  & \TC\ar[dl,"p"]\ar[dr]\ar[drr,"\tilde f"] & & \\
  \OC & & \mathcal C\ar[r,"f" below]\ar[d,"\psi"] & X \\
  & & \mathcal T &
 \ecd
\end{rem}
\begin{dfn}
 Let $\TAtn(X,\beta)^{\text{fact}}\subseteq\TAtn(X,\beta)$ be the locus of maps such that $\tilde f\colon\TC\to X$ factors through a map $\bar f\colon\OC\to X$. We call it the moduli space of aligned admissible maps satisfying the \emph{first factorisation property}.
\end{dfn}
\begin{rem}
 The map $f$ need \emph{not} factor through the admissible cover $\psi$.
\end{rem}

\begin{prop}
 $\TAtn(X,\beta)^{\text{fact}}\subseteq\TAtn(X,\beta)$ is a closed substack. If $X$ is smooth, there is a perfect obstruction theory:
 \[(\R^\bullet\bar\pi_*\bar f^*T_X)^\vee\to\mathbb L^\bullet_{\TAtn(X,\beta)^{\text{fact}}/\TAtn}\]
 endowing $\TAtn(X,\beta)^{\text{fact}}$ with a virtual fundamental class in $A_{\operatorname{vdim}}(\TAtn(X,\beta)^{\text{fact}})$, where:
 \[\operatorname{vdim}(\VZ(X,\beta))=3-\dim(X)+n-K_X\cdot\beta.\]
\end{prop}
\begin{proof}
 For the first claim, we refer the reader to \cite[Theorem 4.3]{RSPW1}.
 The second claim goes back to K. Behrend and B. Fantechi \cite[Proposition 6.3]{BF}.
\end{proof}

Although $\TAtn(\PP^r,d)^{\text{fact}}$ is not necessarily smooth, as it may still have a hyperelliptic component, it can already be useful to the end of computing the invariants of a projective complete intersection.

\begin{lem}
 If $\iota\colon X\hookrightarrow\PP^r$ is a complete intersection of degree $(l_1,\ldots,l_h)$ ($l_i\geq2$ for all $i$),
 \[ \TAtn(X,\beta)^{\text{fact}}\subseteq \TAtn(\PP^r,\iota_*\beta)^{\text{fact}}\]
 is cut out by a section of the vector bundle \[\mathcal E_{(l_1,\ldots,l_k)}=\bar\pi_*\bar f^*\left(\bigoplus_{i=1}^h\OO_{\PP^r}(l_i)\right).\]
 In particular, the invariants satisfy the \emph{quantum Lefschetz} hyperplane principle.
\end{lem}
\begin{proof}
 Thanks to the last statement of Theorem \ref{thm:pushout}, the pullback $f^*\OO_{\PP^r}(1)$ has degree at least one on any subcurve of positive genus. Then, the degree has to be at least two on the minimal subcurve of genus two in order for a non-constant map to exist. This implies that $\bar f^*\OO_{\PP^r}(l_i)$ has vanishing $h^1$ along the fibres of $\OC$ for $l_i\geq 2$. Therefore, $\mathcal E_{(l_1,\ldots,l_k)}$ is a vector bundle by ``cohomology and base-change''. The virtual statement follows as in genus zero from \cite{KKP}.
\end{proof}

\begin{lem}
 The projection $\TAtn(\PP^r,d)^{\text{fact}}\to\TAtn$ factors through $\TAtn^\circ$.
\end{lem}
\begin{proof}
 See Definition \ref{def:niceopen} for the notation. The key observation is that the only (non-constant) maps of degree two from a minimal Gorenstein curve of genus two are those which factor through the hyperelliptic cover of a line. Hence, we may discard the locus where the weight of $\tropD$ is two but (partly) supported away from $D$.
\end{proof}

Thanks to this lemma we can make the following

\begin{dfn}
 Let $\doublewidetilde{\mathcal A_{2,n}}(X,\beta)$ be the moduli space of \emph{admissible maps with a secondary alignment} defined by the fibre diagram:
 \bcd
 \doublewidetilde{\mathcal A_{2,n}}(X,\beta)\ar[d]\ar[r]\ar[dr,phantom,"\Box"] & \TAtn(X,\beta)^\text{fact}\ar[d]\\
 \doublewidetilde{\mathcal A_{2,n}}\ar[r] & \TAtn^\circ.
 \ecd
 This space comes with universal structures:
 \bcd[cramped]
  & \TTC\ar[dl,"\tilde p"]\ar[dr]\ar[drr,"\tilde{\tilde{f}}"] & & \\
  \OOC & & \mathcal C\ar[r,"f" below]\ar[d,"\psi"] & X \\
  & & \mathcal T &
 \ecd
Let $\VZ(X,\beta)\subseteq\doublewidetilde{\mathcal A_{2,n}}(X,\beta)$ be the locus of maps such that $\tilde{\tilde{f}}\colon\TTC\to X$ factors through a map $\bar{\bar{f}}\colon\OOC\to X$. We call it the moduli space of aligned admissible maps satisfying the \emph{second factorisation property}.\footnote{The notation is reminiscent of the celebrated desingularisation of $\pazocal{M}_{1,n}(\PP^r,d)^{\text{main}}$ due to R. Vakil and A. Zinger in {\cite{VZ}}.}
\end{dfn}

\begin{dfn}
 Let $\ev_i\colon \VZ(X,\beta)\to X$ denote the evaluation map at the $i$-th marked point, $i=1,\ldots,n$. Let $\alpha_1,\ldots,\alpha_n\in H^*(X)$ be cohomology classes on the target manifold. The \emph{reduced genus two Gromov-Witten invariants} are defined as:
 \[\langle\alpha_1,\ldots,\alpha_n\rangle^{X,\text{red}}_{2,\beta,n}=\int_{[\VZ(X,\beta)]^\text{vir}}\ev_1^*\alpha_1\cup\ldots\cup\ev_n^*\alpha_n.\]
\end{dfn}

\begin{thm}
 For $d\geq 3$, $\VZ(\PP^r,d)$ is a desingularisation of $\overline{\pazocal M}_{2,n}(\PP^r,d)^\text{main}$.
\end{thm}
\begin{proof}
 Consider the factorisation:
 \[ \VZ(\PP^r,d)\to \mathfrak{P}ic_{\TAtn}\to\TAtn. \]
 Obstructions to the first map can be found in $H^1(\overline{C},L)$, where $L=\bar f^*\OO_{\PP^r}(1)$. If there are two disconnected subcurves of genus one, the degree must be positive on each of them by Theorem \ref{thm:markings_prim}. On the other hand, if the core is a minimal subcurve of genus two, it must have degree at least three by Theorem \ref{thm:markings_sec}, and either the special component has positive degree (at least two) by Corollary \ref{cor:Dpositiveweight}, or the special component is a collapsed ribbon with at least two tails of positive degree by Corollary \ref{cor:oneitherside}. It follows from Lemmas \ref{lem:h1van1} and \ref{lem:h1van2} that the obstructions vanish.
 
 Obstructions to the second map lie in $H^2(\overline{C},\OO)$, which vanishes by dimension reasons. The map is therefore unobstructed. The base is smooth by Theorem \ref{thm:doubletildeA}. We conclude that $\VZ(\PP^r,d)$ is smooth as well. Since it is proper and it contains the locus of maps from a smooth curve as an open dense, $\VZ(\PP^r,d)\to \overline{\pazocal M}_{2,n}(\PP^r,d)^\text{main}$ is birational (see Remarks \ref{rem:mainbirational} and \ref{rem:logetale}).
\end{proof}

\begin{rem}
 A posteriori, we note that aligning and the factorisation property do not alter the main component of $\overline{\pazocal M}_2(\PP^r,2)$, which therefore is already smooth.
\end{rem}

\newpage

\newcommand{\etalchar}[1]{$^{#1}$}

\medskip

\noindent Luca Battistella\\
Mathematisches Institut, Ruprecht-Karls-Universit\"at Heidelberg \\
\texttt{lbattistella@mathi.uni-heidelberg.de}\\

\noindent Francesca Carocci\\
Ecole Polytechnique F\'ed\'erale de Lausanne \\
\texttt{francesca.carocci@epfl.ch}

\end{document}